\documentclass[12pt]{article}
\usepackage{amsmath}
\usepackage{amsfonts}
\usepackage{amssymb}

\newtheorem{Def}{Definition}[section]
\newtheorem{Prop}{Proposition}[section]
\newtheorem{Le}{Lemma}[section]
\newtheorem{Cor}{Corollary}[section]
\newtheorem{Rem}{Remark}[section]
\newtheorem{Ex}{Example}[section]

\newcommand\ds{\displaystyle}
\newcommand\N{{\mathbb{N}}}

\newcommand\R{{\mathbb{R}}}
\newcommand\C{{\mathbb{C}}}

\newcommand\Cc{\mathcal{C}}

\newcommand\Fc{\mathcal{F}}

\newcommand\Hc{\mathcal{H}}
\newcommand\Kc{\mathcal{K}}
\newcommand\Lc{\mathcal{L}}
\newcommand\Mc{\mathcal{M}}
\newcommand\Nc{\mathcal{N}}

\newcommand\Rc{\mathcal{R}}
\newcommand\Sc{\mathcal{S}}
\newcommand\Tc{\mathcal{T}}

\newcommand\Vc{\mathcal{V}}

\newcommand\Yc{\mathcal{Y}}

\newcommand\Rcs{\Rc^{*}}
\newcommand\ld{\lambda}
\newcommand\sg{\sigma}
\newcommand\Jl{J_{\ld}}
\newcommand\gl{g_{\ld}}
\newcommand\ol{\omega_{\ld}}
\newcommand\scl{\Sc_{\ld}}
\newcommand\n[2]{\nu_{#1}^{#2}}
\newcommand\sig[2]{\sg_{#1}^{#2}}
\newcommand\Ker{\mathrm{Ker}}
 
\newcommand\Uh{U^{1,0}}
\newcommand\Xh{X^{1,0}}
\newcommand\Yh{Y^{1,0}}
\newcommand\Zh{Z^{1,0}}

\newcommand\Xah{X^{0,1}}
\newcommand\Yah{Y^{0,1}}
\newcommand\Zah{Z^{0,1}}
\newcommand\Wah{W^{0,1}}

\newcommand\Xhl{X_{\ld}^{1,0}}
\newcommand\Yhl{Y_{\ld}^{1,0}}
\newcommand\Zhl{Z_{\ld}^{1,0}}

\newcommand\Xahl{X_{\ld}^{0,1}}
\newcommand\Yahl{Y_{\ld}^{0,1}}
\newcommand\Zahl{Z_{\ld}^{0,1}}
\newcommand\Wahl{W_{\ld}^{0,1}}
\newcommand\So[2]{\Sc^{{#1},{#2}}}
\newcommand\Sl[2]{\Sc_{\ld}^{{#1},{#2}}}
\newcommand\om[1]{\omega_{#1}}
\newcommand\op{\stackrel{\perp}{\oplus}}
\newcommand\tr{\mathrm{Tr}_{\Sc}}
\newcommand\tu[1]{\Tc_{#1}}

\newcommand\del{\partial_{\Sc}}
\newcommand\delb{{\overline{\partial}}_{\Sc}}
\newcommand\sdel{\partial_{\Sc}^{*}}
\newcommand\sdelb{{\overline{\partial}}_{\Sc}^{*}}
\newcommand\dls{{d\ld}_{\Sc}}
\newcommand\dlsl{{d\ld}_{\scl}}
\newcommand\dlsd[1]{{(d\ld)}^{#1}_{\scl}}
\newcommand\dlsn{||{(d\ld)}^{\sharp}_{\Sc}||^{2}}
\newcommand\nab[2]{\nabla_{#1}^{#2}}

\begin{document}

\title{\LARGE\bf Robinson manifolds and the Chern-Robinson connection}
\author{Robert Petit*}
\date{}
\maketitle
\vspace{2cm}
\begin{abstract}
In this article, we define the Chern-Robinson connection on the complexify tangent bundle of an almost Robinson manifold and we study  the curvature associated to. Various Bianchi identities are obtained together with an application to geometry of some Robinson manifolds.
\end{abstract}

\vspace{6cm}\noindent
{\bf Author Keywords:} Robinson manifolds; Witt structures; Chern-Robinson connection; Chern-Robinson curvature tensor; Chern-Moser-Robinson tensor; Bianchi identities

\vspace{2.5cm}\noindent
{\footnotesize *Laboratoire de Math\'ematiques Jean Leray, UMR 6629 CNRS, Nantes Universit\'e,\\
2, rue de la Houssini\`ere BP 92208, 44322 Nantes - France.\\
E-mail address : robert.petit@univ-nantes.fr}

\newpage
\section{Introduction}
The almost Robinson manifolds are even dimensional Lorentzian manifolds endowed with a complex subbundle of the complexify tangent bundle which is maximally totally null with respect to the complexified metric. These manifolds can be viewed as Lorentzian analogues of almost Hermitian manifolds. This geometric structure appears in the context of the general relativity for some solutions of Einstein equations like gravitational waves or black holes (Kerr or Taub-NUT black holes). The above formal definition is quite recent and due to Nurowski and Trautman (\cite{NT}). The interest focused on the Robinson manifolds has increased in the last years. Among the works on the subject, the following \cite{AGS},\cite{AGSS},\cite{FLTC1},\cite{FLTC2},\cite{HLN},\cite{Pa},\cite{TC1},\cite{TC2} deal with relationships between Robinson manifolds and CR manifolds and also about classification of almost Robinson manifolds. For instance, in \cite{FLTC2}, the authors give a classification of almost Robinson manifolds based on the concept of intrinsic torsion. In this article, we continue to investigate the geometry of almost Robinson manifolds (started in \cite{Pe2}) by means of the Witt structures. Recall that a Witt structure on the (complexify) tangent bundle of a pseudo-Riemannian manifold is a decomposition into totally null subbundles (directly inspired by the Witt decomposition of pseudo-euclidean spaces). As described in (\cite{FLTC2},\cite{Pe2}), the maximally totally null complex subbundle of an almost Robinson manifold gives rise to a Witt decomposition on the complexify tangent bundle of the almost Robinson manifold. Also it is possible to equiped the (complexify) tangent bundle together with a metric connection preserving the Witt structure (note that a such connection necessarily possesses torsion). The so called Lichnerowicz connection introduced in (\cite{Pe2}) provides a such example, however, there is no canonical choice for this type of connections. In this article, we introduce an other metric connection preserving the Witt structure called the Chern-Robinson connection with the property to coincide with the $\overline{\partial}$-operator associated to the almost Robinson structure (which is a counterpart of the Chern connection in Hermitian geometry (\cite{Ch},\cite{Ga})). Now, by means of this connection together with its curvature and its torsion, we can investigate the geometry of almost Robinson manifolds from a new perspective. In that follows, we restrict this investigation mainly to so called almost Fefferman-Robinson manifolds (subclass of almost Robinson manifolds realized as line bundles over CR strictly pseudoconvex manifolds). The plan of this article is as follow. In section 2, the definitions of almost optical and almost Robinson manifolds are reviewed thus examples. We also describe the Witt structure on their complexify tangent bundle. Section 3 is devoted to construct the Chern-Robinson connection for an almost Robinson manifold and to derive first Bianchi identity for the curvature of the connection of almost Fefferman-Robinson manifolds (Propositions 3.1 and 3.2). In the last part of this section, we define various type of curvatures, especially, we define the Chern-Robinson curvature tensor and the Chern-Moser-Robinson tensor (which is the counterpart of the Chern-Moser tensor \cite{CM} for CR strictly pseudoconvex manifolds). The Section 4 is devoted to prove (Proposition 4.2 and Corollary 4.2) that the Chern-Moser-Robinson tensor is an optical conformal invariant for the strongly geodetic (in the sense of Definition 2.3 after) almost Fefferman-Robinson manifolds. In Section 5, we derive a second Bianchi identity for the Chern-Robinson curvature tensor and the Chern-Moser-Robinson tensor of a Fefferman-Robinson manifold (Propositions 5.2 and 5.3). As application of the second Bianchi identity, we give a formula expressing explicitly the Chern-Moser-Robinson tensor of an Einstein-Fefferman-Robinson manifold with parallel torsion in terms of the metric and the torsion (Proposition 5.4). This formula is, in particular, valid for a locally symmetric Fefferman-Robinson manifold Einstein-Fefferman (Corollary 5.1) and is similar to those already obtained in \cite{Pe1} for contact locally subsymmetric spaces.

\section{Optical and Robinson manifolds}
Most of definitions given in this section can be found in \cite{AGS},\cite{AGSS},\cite{FLTC1},\cite{FLTC2},\cite{NT},\cite{Pa}.\\

\noindent
In the following, $(M^{2m+2},g)$ is a $(2m+2)$-dimensional Lorentzian manifold.
\begin{Def}
An almost optical structure (cf. \cite{FLTC1}) on $(M^{2m+2},g)$ is a null line subbundle $\Rc$ of $TM$.
\end{Def}
An almost optical structure gives rise to a filtration and a grading of $TM$ :
$$\{0\}\subset\Rc\subset\Rc^{\perp}\subset TM,$$
\begin{equation}\label{e1}
TM=(\Rc\oplus TM/{\Rc}^{\perp})\oplus({\Rc}^{\perp}/\Rc).
\end{equation}
The vector bundle $\Sc:={\Rc}^{\perp}/\Rc$ is called screen bundle associated to $\Rc$. Note that the  Lorentzian metric $g$ induces a riemannian metric $g_{\Sc}$ on $\Sc$.

\begin{Def}
An almost Robinson (or complex almost optical) structure (cf. \cite{FLTC2},\cite{NT},\cite{Pa}) on $(M^{2m+2},g)$ is a rank-$(m+1)$ complex subbundle $\Nc$ of $T^{\C}M$ which is totally null with respect to the complexified metric $g^{\C}$ (i.e $\Nc^{\perp}=\Nc$).
\end{Def}
In this case, $\Nc\cap\overline{\Nc}$ is a complex null line subbundle of $T^{\C}M$ and there exists a real null line subbundle $\Rc$ of $TM$ such that $\Rc^{\C}=\Nc\cap\overline{\Nc}$. ($V^{\C}$ is the complexification of $V$).\\
We obtain the following filtration and grading of $T^{\C}M$ :
$$\{0\}\subset\Rc^{\C}\subset{\Rc^{\C}}^{\perp}=\Nc+\overline{\Nc}\subset T^{\C}M,$$
\begin{equation}\label{e2}
T^{\C}M=(\Rc^{\C}\oplus{(TM/{\Rc}^{\perp})}^{\C})\op (\Nc/\Rc^{\C}\oplus\overline{\Nc}/\Rc^{\C}).
\end{equation}

It follows that $\Rc$ is an almost optical structure associated to the almost Robinson structure $\Nc$ together with $\Rc^{\perp}=\mathrm{Re}\,(\Nc+\overline{\Nc})$. It is possible to define an almost complex structure $J_{\Sc}$ on $\Sc$ such that $\Nc/\Rc^{\C}$ (respectively $\overline{\Nc}/\Rc^{\C}$) is the $\sqrt{-1}$ (respectively $-\sqrt{-1}$)-eigenbundle of $J_{\Sc}$ on $\Sc^{\C}=(\Nc+\overline{\Nc})/(\Nc\cap\overline{\Nc})$ together with a Hermitian scalar product $g_{\Sc}$ on $\Sc$. In the following, the Hermitian form associated to $(g_{\Sc},J_{\Sc})$ will be denoted by $\omega_{\Sc}$ (i.e. $\om{\Sc}(X,Y)=g_{\Sc}(J_{\Sc}X,Y)$, $X,Y\in\Gamma(\Sc)$) and the subbundles $\Nc/\Rc^{\C}$, $\overline{\Nc}/\Rc^{\C}$ will be respectively denoted by $\So{1}{0}$ and $\So{0}{1}$ (and so $\Nc=\Rc^{\C}\oplus\So{1}{0}$). Note that an almost optical structure whose screen bundle is endowed with an almost Hermitian structure give rise to an almost Robinson structure.

\begin{Ex} Basic example of almost Robinson structure. Let $\R^4=\R\times\R\times\C$ with coordinates $(u,v,z)$ and the metric $g=du\odot dv+dz\odot d\overline{z}$ (with the convention $\alpha\odot\beta=\alpha\otimes\beta+\beta\otimes\alpha$ for $1$-forms $\alpha,\beta$). Then $\ds\Nc=span_{\C}\bigl(\frac{\partial}{\partial v},\frac{\partial}{\partial z}\bigr)\subset\C^4$ is an almost Robinson structure on $(\R^4,g)$ together with $\ds\Rc=span\bigl(\frac{\partial}{\partial v}\bigr)$. 
\end{Ex}

\begin{Rem} Note that almost optical and almost Robinson structures also exist on the odd dimensional Lorentzian manifolds but we restrict here to the even dimensional case.
\end{Rem}

\newpage

\begin{Def}
\begin{enumerate}
\item An almost optical structure $\Rc$ on $(M^{2m+2},g)$ is called geodetic if 
$$[\Gamma(\Rc),\Gamma(\Rc^{\perp})]\subset\Gamma(\Rc^{\perp}).$$
\item An almost Robinson structure $\Nc$ on $(M^{2m+2},g)$ is called\\ geodetic if 
$$[\Gamma(\Rc^{\C}),\Gamma(\Nc)]\subset\Gamma(\Nc+\overline{\Nc})$$
strongly geodetic if 
$$[\Gamma(\Rc^{\C}),\Gamma(\Nc)]\subset\Gamma(\Nc)$$
partially integrable if  
$$[\Gamma(\Nc),\Gamma(\Nc)]\subset\Gamma(\Nc+\overline{\Nc})$$
integrable if  
$$[\Gamma(\Nc),\Gamma(\Nc)]\subset\Gamma(\Nc).$$
\end{enumerate}
\end{Def}

\begin{Rem}
\begin{enumerate}
\item Note that the condition $[\Gamma(\Rc^{\C}),\Gamma(\Nc)]\subset\Gamma(\Nc)$ is called nearly Robinson in \cite{FLTC2} but the term strongly geodetic seems more suitable.
\item If $\Nc$ is an almost Robinson structure together with its almost optical structure associated $\Rc$, then $\Nc$ geodetic is equivalent to $\Rc$ geodetic.
\item The integrability condition $[\Gamma(\Nc),\Gamma(\Nc)]\subset\Gamma(\Nc)$ is equivalent to $d\,\Gamma(\Nc^{0})\subset\Gamma(\Nc^{0})\wedge\Gamma({(T^{\C}M)}^{*})$ where $\Nc^{0}\subset {(T^{\C}M)}^{*}$ is the set of linear forms vanishing on $\Nc$ and $d$ is the exterior derivative.
\end{enumerate}
\end{Rem}

For an almost optical structure $\Rc$, the choice of a null line subbundle $\Rcs$ of $TM$ (dual to $\Rc$) such that $TM=\Rcs\oplus{\Rc}^{\perp}$ allows to rewrite grading (\ref{e1}), as : 
\begin{equation}\label{e3}
TM=(\Rc\oplus\Rcs)\op\Sc,
\end{equation}
with $\Sc\simeq{\Rc}^{\perp}\cap{(\Rcs)}^{\perp}$.
Now if $\Nc$ is an almost Robinson structure with almost optical structure $\Rc$ associated to, then grading (\ref{e2}) is :
\begin{equation}\label{e4}
T^{\C}M=(\Rc^{\C}\oplus{(\Rcs)}^{\C})\op(\So{1}{0}\oplus\So{0}{1}),
\end{equation}
with $\So{1}{0}\simeq\Nc\cap\Sc^{\C}$ and $\So{0}{1}\simeq\overline{\Nc}\cap\Sc^{\C}$.

\begin{Def} A Lorentzian manifold $(M^{2m+2},g)$ endowed with an almost optical structure $\Rc$ will be called an almost optical manifold.\\
A Lorentzian manifold $(M^{2m+2},g)$ endowed with an almost Robinson structure $\Nc$ will be called an almost Robinson manifold. A Robinson manifold is an almost Robinson manifold such that $\Nc$ integrable.
\end{Def}

If $Z$ is a vector field on $M$ (respectively a $1$-form $\alpha$ on $M$), then we denote by $Z^{\flat}$ (respectively $\alpha^{\sharp}$) the $1$-form on $M$ given by ${Z}^{\flat}(X)=g(Z,X)$ (respectively the vector field on $M$ given by $\alpha(X)=g(\alpha^{\sharp},X)$). 
\newpage
\begin{Def} Let $(M^{2m+2},g,\Rc)$ be an almost optical manifold.
\begin{enumerate}
\item A local non vanishing section $\n{}{}\in\Gamma(\Rc)$ is called an optical vector field and the $1$-form ${\n{}{}}^{\flat}$ is called optical form associated to.
\item An optical vector field $\n{}{}\in\Gamma(\Rc)$ is called a geodetic vector field if $(\Lc_{\n{}{}}\sig{}{*})=f\sig{}{*}$ with $\sig{}{*}={\n{}{}}^{\flat}$ and $f\in\Cc^{\infty}(M)$ (equivalently, ${(\Lc_{\n{}{}}\sig{}{*})}_{/{\Rc}^{\perp}}=0$ with ${\Rc}^{\perp}=\mathrm{Ker}\,\sig{}{*}$).
\item An optical vector field $\n{}{}\in\Gamma(\Rc)$ is called a shearfree vector field if $\Lc_{\n{}{}}g=\rho g+\sig{}{*}\odot\alpha$, with $\rho\in\Cc^{\infty}(M)$ and $\alpha\in\Omega^{1}(M)$. The function $\rho$ is called expansion.
\end{enumerate}
\end{Def}

Note that $\Rc$ geodetic if and only if any optical vector field $\n{}{}\in\Gamma(\Rc)$ is a geodetic vector field and that a shearfree vector field is always a geodetic vector field.\\

From now, we will always assume that an almost optical manifold $(M^{2m+2},g,\Rc)$ will be endowed with a null line bundle $\Rcs$ such that $TM=\Rcs\oplus{\Rc}^{\perp}$. 

\begin{Def} Let $(M^{2m+2},g,\Rc,\Rcs)$ be an almost optical manifold.
We call $g$-optical pairing a pairing of optical vector fields $(\n{}{},\n{}{*})\in\Gamma(\Rc\oplus\Rcs)$ together with $g(\n{}{},\n{}{*})=1$. The $1$-forms $(\sig{}{}={\n{}{*}}^{\flat},\sig{}{*}={\n{}{}}^{\flat})$ will be called optical forms associated to the $g$-optical pairing $(\n{}{},\n{}{*})$\\(note that $\sig{}{}(\n{}{})=\sig{}{*}(\n{}{*})=1$ and $\sig{}{}(\n{}{*})=\sig{}{*}(\n{}{})=0$).
\end{Def}

Given a $g$-optical pairing $(\n{}{},\n{}{*})\in\Gamma(\Rc\oplus\Rcs)$ together with $(\sig{}{},\sig{}{*})$ optical forms associated, then we have 
$${\Rc}^{\perp}=\mathrm{Ker}\,\sig{}{*},\quad {(\Rcs)}^{\perp}=\mathrm{Ker}\,\sig{}{},\quad \Sc=\mathrm{Ker}\,\sig{}{}\cap\mathrm{Ker}\,\sig{}{*}(={(\Rcs)}^{\perp}\cap{\Rc}^{\perp})$$
and $\Lc_{\n{}{(*)}}\sig{}{(*)}=\n{}{(*)}\lrcorner\,d\sig{}{(*)}$ (by the Cartan formula $\Lc_{X}=d\circ X\lrcorner+X\lrcorner\circ d$ where $X\lrcorner$ denotes the interior product by $X$).\\

If $\Pi_{\Sc}:TM\to\Sc$ is the orthogonal projection on $\Sc$ associated to (\ref{e3}), then we have the decompositions:
\begin{eqnarray}\label{e5}
X\in\Gamma(TM)&=&X_{\Sc}+\sig{}{}(X)\n{}{}+\sig{}{*}(X)\n{}{*},\quad X_{\Sc}=\Pi_{\Sc}(X)\nonumber\\
g&=&g_{\Sc}+\sig{}{*}\odot\sig{}{},\quad g_{\Sc}=g\circ\Pi_{\Sc}\nonumber\\
\alpha\in\Omega^{1}(M)&=&\alpha_{\Sc}+\sig{}{}\otimes\alpha(\n{}{})+\sig{}{*}\otimes\alpha(\n{}{*}),\quad \alpha_{\Sc}=\alpha\circ\Pi_{\Sc}\nonumber\\
\alpha\in\Omega^{p}(M)&=&\alpha_{\Sc}+\sig{}{}\wedge {(\n{}{}\lrcorner\,\alpha)}_{\Sc}+\sig{}{*}\wedge {(\n{}{*}\lrcorner\,\alpha)}_{\Sc}+\sig{}{}\wedge\sig{}{*}\wedge {(\n{}{}\lrcorner\,(\n{}{*}\lrcorner\,\alpha))}_{\Sc}\quad (p\geq 2).
\end{eqnarray}
Let $\alpha_{\Sc}\in\Omega_{\Sc}^{1}(M)=\Gamma({(\Sc^{\C})}^*)$ and $\Pi_{\Sc^{\C}}^{1,0}:\Sc^{\C}\to\So{1}{0},\Pi_{\Sc^{\C}}^{0,1}:\Sc^{\C}\to\So{0}{1}$ the canonical projections associated to (\ref{e4}), then $\alpha_{\Sc}^{1,0}:=\alpha_{\Sc}\circ\Pi_{\Sc^{\C}}^{1,0}$ and $\alpha_{\Sc}^{0,1}:=\alpha_{\Sc}\circ\Pi_{\Sc^{\C}}^{0,1}$.\\
Note that
\begin{equation}\label{e6}
\Lc_{\n{}{}}g={(\Lc_{\n{}{}}g)}_{\Sc}+\Lc_{\n{}{}}\sig{}{*}\odot\sig{}{}+\sig{}{*}\odot(\Lc_{\n{}{}}\sig{}{}-[\n{}{},\n{}{*}]_{\Sc}^{\flat}). 
\end{equation}

\begin{Ex}
Let $(N^{2m+1},H,J_N)$ be an orientable almost $CR$ manifold. Since $N$ is orientable, we can choose a non-vanishing $1$-form $\theta$ and a non-vanishing vector field $\xi$ on $N$ such that $\theta(\xi)=1$ and $\Ker\,\theta=H$ ($\theta$ is a pseudo-Hermitian structure). Note that $TN=H\oplus span(\xi)$. Assume moreover that there exists a symmetric $2$-tensor field $h$ on $TN$ such that $\Ker\,h=span(\xi)$ and ($h_{/H},J_N)$ is an almost hermitian structure on $H$. Let $\pi:M^{2m+2}\to N^{2m+1}$ be a $G$-principal bundle over $N$ with $G=\R,\R^+\,\mathrm{or}\,S^1$. Then we have the splitting of $TM$ into horizontal and vertical spaces :
$$TM=\Vc_M\oplus\Hc_M=\Vc_M\oplus (span(\xi^{*})\oplus\pi_{*}^{-1}(H)),$$ 
with $\xi^{*}=\pi_{*}^{-1}(\xi)$.  We fix a non-vanishing vector field $\n{}{}\in\Gamma(\Vc_M)$ and we consider the $1$-form $\alpha$ on $M$ satisfying $\alpha(\n{}{})=1$ and $\Ker\,\alpha=\Hc_M$. We set :
\begin{equation}\label{e7}
\sig{}{*}=\pi^{*}\theta,\quad\sig{}{}=\alpha+\lambda\pi^{*}\theta,\quad 
\n{}{*}=\xi^{*}-\ld\n{}{},\quad g=\pi^{*}h+\sig{}{*}\odot\sig{}{},
\end{equation}
with $\ld\in\Cc^{\infty}(M)$. Then $(M^{2m+2},g,\Nc,\Rc=span(\n{}{}),\Rcs=span(\n{}{*}))$ with $\Nc=\Rc^{\C}\oplus\pi_{*}^{-1}(H^{1,0})$ 
and $\pi_{*}^{-1}(H^{1,0})=\{X\in\pi_{*}^{-1}(H^{\C});\,(\pi_{*}^{-1}\circ J_N\circ\pi_{*})(X)=\sqrt{-1}X\}$ is an almost Robinson structure and $(\n{}{},\n{}{*})$ is a $g$-optical pairing together with $(\sig{}{},\sig{}{*})$ optical forms associated to.\\

Let $M$ be a $\R^+$-principal bundle over $N$. Denote by $r\in\R^+$ the coordinate on the fiber and by $(u,z^{i},z^{\bar{i}})\in\R\times\C^{2m}$ (with $z^{\bar{i}}=\overline{z^{i}}$) local coordinates on an open set $U\subset N$, then a metric of form (\ref{e7}) on $\R^+\times U$ is given by :
$$g=Ph_{i\bar{j}}dz^{i}\odot dz^{\bar{j}}+\bigl(du+A_{i}dz^{i}+A_{\bar{i}}dz^{\bar{i}}\bigr)\odot\bigl(dr+B_{i}dz^{i}+B_{\bar{i}}dz^{\bar{i}}+Q(du+A_{i}dz^{i}+A_{\bar{i}}dz^{\bar{i}})\bigr),$$
with $A_{i},B_{i},h_{i\bar{j}}$ complex valued functions on $\R^+\times U$ ($h_{i\bar{j}}$ independent of variable $r$) and $P,Q$ real valued functions on $\R^+\times U$ ($P$ function of $r$ only). 
\end{Ex}

\begin{Rem} Note that in dimension $4$, black holes metrics like Taub-Nut metric or Kerr metric have the form described by the above $g$ (cf. \cite{AGS},\cite{HLN}).
\end{Rem}

\begin{Ex}
Let $(\tilde{M}^{2m},\tilde{g},\tilde{\Nc},\tilde{\Rc}=span(\tilde{\nu}),\tilde{\Rc}^{*}=span(\tilde{\nu}^{*}))$ be an almost Robinson manifold and $\pi:M^{2m+2}\to \tilde{M}$ be a $S^2$-principal bundle over $\tilde{M}$. In this case, 
$$TM=\Vc_M\oplus\Hc_M=\Vc_M\oplus (\pi_{*}^{-1}(\tilde{\Rc})\oplus\pi_{*}^{-1}(\tilde{\Rc}^{*})\oplus\pi_{*}^{-1}(\tilde{\Sc})).$$
Setting $$\n{}{}=\pi_{*}^{-1}(\tilde{\nu}),\quad \n{}{*}=\pi_{*}^{-1}(\tilde{\nu}^{*}),\quad\sig{}{*}=\pi^{*}\tilde{\sig{}{}}^*,\quad\sig{}{}=\pi^{*}\tilde{\sig{}{}},\quad g=g_{\Vc_M}+\pi^{*}\tilde{g},$$
and defining an almost complex structure $J$ on $\Vc_M\oplus\pi_{*}^{-1}(\tilde{\Sc})$ by 
$J=J_{\Vc_M}\oplus(\pi_{*}^{-1}\circ\tilde{J}\circ\pi_{*})$,
then $(M^{2m+2},g,\Rc=span(\n{}{}),\Rcs=span(\n{}{*}),\Nc=\Rc^{\C}\oplus(\Vc^{1,0}_{M}\oplus\pi_{*}^{-1}(\tilde{\Sc}^{1,0}))$ is an almost Robinson manifold and $(\n{}{},\n{}{*})$ is a $g$-optical pairing.
\end{Ex}
 
In the following, we simply denote $X^{1,0}$ (respectively $X^{0,1}$) for $X_{\Sc}^{1,0}\in\Gamma({\Sc}^{1,0})$ (respectively $X_{\Sc}^{0,1}\in\Gamma({\Sc}^{0,1})$).

\subsection{Fefferman-Robinson manifold}

\begin{Def} Let $\Rc$ be a geodetic almost optical structure on $(M^{2m+2},g)$ and $\n{}{}\in\Gamma(\Rc)$ a geodetic vector field. Then $\Rc$ is called twistfree if $\sig{}{*}\wedge d\sig{}{*}=0$ (equivalently if ${(d\sig{}{*})}_{/{\Rc}^{\perp}}=0$) and $\Rc$ is called maximally twisted if $\sig{}{*}\wedge (d\sig{}{*})^m\neq 0$ (equivalently if ${(d\sig{}{*})}_{/{\Rc}^{\perp}}$ has rank $2m$).
\end{Def} 

\begin{Rem}
Note that $\Rc$ maximally twisted means that ${\Rc}^{\perp}$ is an even contact structure on $M$ together with characteristic line field equal to $\Rc$.
\end{Rem}

In Proposition 4.29 of \cite{FLTC1}, the authors show that for a maximally twisted almost optical structure $\Rc$ then there exists a $g$-optical pairing $(\n{}{},\n{}{*})\in\Gamma(\Rc\oplus\Rcs)$ such that $(\Lc_{\n{}{}}\sig{}{*})=(\Lc_{\n{}{*}}\sig{}{*})=0$. In this case, we say that the pairing $(\n{}{},\n{}{*})$ is adapted. From now, we will assume that a maximally twisted almost optical structure will be always endowed with an adapted $g$-optical pairing $(\n{}{},\n{}{*})$.

\begin{Def} An almost Fefferman-Robinson manifold is an almost optical manifold $(M^{2m+2},g,\Rc,\Rcs)$ with $\Rc$ maximally twisted and such that $J_{\Sc}=g_{\Sc}^{-1}\circ d{\sig{}{*}/}_{\Sc}$ is an almost Hermitian structure for $g_{\Sc}$.
\end{Def}

If we denote by $\omega$ the extension to $TM$ of the Hermitian form $\omega_{\Sc}$ on $\Sc$ (i.e. $\omega(X,Y)=g(JX,Y)$ where $J$ is the extension of $J_{\Sc}$ to $TM$ by $0$ on $\Rc\oplus\Rcs$), then an almost Fefferman-Robinson manifold has almost Robinson structure $\Nc$ together with $g=d\sig{}{*}(.,J.)+\sig{}{*}\odot\sig{}{}$ and $\omega=d\sig{}{*}$.\\A Fefferman-Robinson manifold is an almost Fefferman-Robinson with $\Nc$ integrable.

\begin{Prop} Let $(M^{2m+2},g,\Nc,\Rc,\Rcs)$ be an almost Fefferman-Robinson manifold together with $(\n{}{},\n{}{*})$ adapted $g$-optical pairing. Then 
\begin{enumerate}
\item $\Nc$ is partially integrable.
\item $\Lc_{\n{}{}}\omega=\Lc_{\n{}{*}}\omega=0$ and $[\n{}{},\n{}{*}]_{\Sc}=0$.
\item $\Lc_{\n{}{}}g={(\Lc_{\n{}{}}g)}_{\Sc}^{(2,0)+(0,2)}+\sig{}{*}\odot\Lc_{\n{}{}}\sig{}{}$.
\end{enumerate}
\end{Prop}

\noindent{Proof.}  Let $(\n{}{},\n{}{*})\in\Gamma(\Rc\oplus\Rcs)$ be an adapted $g$-optical pairing. Since $M$ is an almost Fefferman-Robinson manifold then $\Lc_{\n{}{}}\sig{}{*}=0$ and $\om{\Sc}^{2,0}={(d\sig{}{*})}_{\Sc}^{2,0}=0$. Also 1. directly follows. Now $\omega=d\sig{}{*}$ and $\n{}{}\lrcorner\,\omega=\n{}{*}\lrcorner\,\omega=0$ which implies by the Cartan formula that $\Lc_{\n{}{}}\omega=\Lc_{\n{}{*}}\omega=0$. The relation
$\n{}{*}\lrcorner\,(\Lc_{\n{}{}}\omega)=0$ gives $[\n{}{},\n{}{*}]_{\Sc}=0$ and 2. is obtained.\\
For $X\in\Gamma(\Sc^{\C})=\Xh+\Xah$ and $Y\in\Gamma(\Sc^{\C})=\Yh+\Yah$ we have, by (\ref{e6}) together with the assumptions $\Lc_{\n{}{}}\sig{}{*}=0$ and $[\n{}{},\n{}{*}]_{\Sc}=0$ that : 
\begin{eqnarray*}
 (\Lc_{\n{}{}}g)(X,Y)&=&(\Lc_{\n{}{}}g)(\Xh,\Yh)+(\Lc_{\n{}{}}g)(\Xah,\Yah)
+(\Lc_{\n{}{}}g)(\Xh,\Yah)+(\Lc_{\n{}{}}g)(\Xah,\Yh)\\
&+&(\sig{}{*}\odot\Lc_{\n{}{}}\sig{}{})(X,Y).
\end{eqnarray*}
Since 
$(\Lc_{\n{}{}}g)(\Xh,\Yah)=\sqrt{-1}(\Lc_{\n{}{}}\omega)(\Xh,\Yah)=0$ and 
$(\Lc_{\n{}{}}g)(\Xah,\Yh)=\sqrt{-1}(\Lc_{\n{}{}}\omega)(\Xah,\Yh)=0$, we deduce 3. $\Box$

\begin{Rem}
\begin{enumerate}
 \item It follows from Proposition 2.1 that an almost Fefferman-Robinson manifold together with $(\n{}{},\n{}{*})$ adapted $g$-optical pairing admits a codimension $2m$ foliation $\Fc$ of $M$ such that $T(\Fc)=\Rc\oplus{\Rc}^{*}$ (since $[\n{}{},\n{}{*}]_{\Sc}=0$) and that, for a strongly geodetic almost Fefferman-Robinson manifold, then $\n{}{}$ is a shearfree vector field without expansion (since ${(\Lc_{\n{}{}}g)}_{\Sc}^{(2,0)+(0,2)}=0$).
\item Note that an almost Fefferman-Robinson manifold is an example of conformal twist-induced almost Robinson manifold with non-shearing congruence of null geodesics as defined in Proposition 5.16 of \cite{FLTC1} and Proposition 4.8 of \cite{FLTC2}.
\end{enumerate}

\end{Rem}

\begin{Ex}
Let $(N^{2m+1},H=\Ker\,\theta,\xi,J_N)$ be a partially integrable strictly pseudoconvex $CR$-manifold. Then recall that the Levi form $L_{\theta}(X,Y)=d\theta(X,J_N Y)$ (with $J_N$ extended on $TN$ by $J_N\xi=0$) satisfies $\Ker\,L_{\theta}=span(\xi)$ and $({L_{\theta}}_{/H},J_N)$ is an almost hermitian structure on $H$. If $\pi:M^{2m+2}\to N^{2m+1}$ is a $G$-principal bundle over $N$ with $G=\R,\R^+\,\mathrm{or}\,S^1$ then following construction of example 2.2, we obtain that $(M^{2m+2},g=\pi^{*}L_{\theta}+\sig{}{*}\odot\sig{}{},\Nc=\Rc^{\C}\oplus\pi_{*}^{-1}(H^{1,0}),\Rc=span(\n{}{}),\Rcs=span(\n{}{*}))$ is an almost Fefferman-Robinson manifold. Note that the Fefferman spaces (cf. \cite{Ba},\cite{Fe}) are examples of Fefferman-Robinson manifolds for which the $1$-form $\sig{}{}$ is a connection form.
\end{Ex}

\section{Connection and curvature}
In the following, the $\C$-bilinear extensions of $g$ and $\omega$ to $T^{\C}M$ will be denoted by the same letters.\\
Recall that the torsion and the curvature of a connection $\nabla$ are respectively defined by :  
$$T(X,Y)=\nabla_{X}Y-\nabla_{Y}X-[X,Y]\quad\mathrm{and}\quad R(X,Y)=[\nabla_{X},\nabla_{Y}]-\nabla_{[X,Y]}.$$ 

If $(M^{2m+2},g,\Nc,\Rc,\Rcs)$ is an almost Robinson manifold, then the grading $\ds T^{\C}M=(\Rc^{\C}\oplus{(\Rcs)}^{\C})\op(\So{1}{0}\oplus\So{0}{1})$ gives rise to a complex Witt structure in the sense of definition 2.1 of \cite{Pe2}. In this context, we have defined in lemma 6.1 of \cite{Pe2} a metric connection with torsion preserving the complex Witt structure called the canonical Witt connection. Now, we define a new metric connection with torsion preserving the complex Witt structure called Chern-Robinson connection of the following way :

\begin{Prop} Let $(M^{2m+2},g,\Nc,\Rc,\Rcs)$ be an almost Robinson manifold and \\$T^{\C}M=(\Rc^{\C}\oplus{(\Rcs)}^{\C})\op(\So{1}{0}\oplus \So{0}{1})$ its complex Witt structure. We assume in the following that $(\n{}{},\n{}{*})\in\Gamma(\Rc\oplus\Rcs)$ is a $g$-optical pairing. Then, there exists a unique metric connection $\nabla$ on $T^{\C}M$ preserving each component of the previous complex Witt structure (called the Chern-Robinson connection) with torsion $T$ given by : 
\begin{eqnarray}\label{e8}
T(\Xh,\Yh)&=&-[\Xh,\Yh]^{0,1}+d\sig{}{}(\Xh,\Yh)\n{}{}+d\sig{}{*}(\Xh,\Yh)\n{}{*}+\Tc^{1,0}(\Xh,\Yh)\nonumber\\
T(\Xh,\Yah)&=&d\sig{}{}(\Xh,\Yah)\n{}{}+d\sig{}{*}(\Xh,\Yah)\n{}{*}\nonumber\\
T(\n{}{},\Xh)&=&-[\n{}{},\Xh]^{0,1}+(\Lc_{\n{}{}}\sig{}{*})(\Xh)\n{}{*}-\frac{1}{2}(\Lc_{\Xh}g)(\n{}{},\n{}{*})\n{}{}+\Tc^{1,0}(\n{}{},\Xh)\nonumber\\
T(\n{}{*},\Xh)&=&-[\n{}{*},\Xh]^{0,1}+(\Lc_{\n{}{*}}\sig{}{})(\Xh)\n{}{}-\frac{1}{2}(\Lc_{\Xh}g)(\n{}{},\n{}{*})\n{}{*}+\Tc^{1,0}(\n{}{*},\Xh)\nonumber\\
T(\n{}{},\n{}{*})&=&-[\n{}{},\n{}{*}]_{\Sc},
\end{eqnarray}
with 
\begin{eqnarray*}
g(\Tc^{1,0}(\Xh,\Yh),\Zah)&=&(\Lc_{\Xh}g)(\Zah,\Yh)+d{(\Xh)}^{\flat}(\Zah,\Yh)\\
&=&-\sqrt{-1}(d\omega)(\Xh,\Yh,\Zah)=-(J^{*}d\omega)(\Xh,\Yh,\Zah)\\
g(\Tc^{1,0}(\n{}{(*)},\Xh),\Yah)&=&\frac{1}{2}(\Lc_{\n{}{(*)}}g)(\Xh,\Yah)=-\frac{\sqrt{-1}}{2}(\Lc_{\n{}{(*)}}\omega)(\Xh,\Yah)\\
-g([\n{}{(*)},\Xh]^{0,1},\Yh)&=&-\frac{1}{2}g((J\Lc_{\n{}{(*)}}J)(\Xh),\Yh)\\
&=&\frac{1}{2}(\Lc_{\n{}{(*)}}g)(\Xh,\Yh)+\frac{\sqrt{-1}}{2}(\Lc_{\n{}{(*)}}\omega)(\Xh,\Yh),
\end{eqnarray*}
where $\Xh,\Yh\in\Gamma(\So{1}{0})$ and $\Yah\in\Gamma(\So{0}{1})$.
\end{Prop}
Moreover, we have $\nabla_{\Yah}\Xh={\overline{\partial}}_{\Yah}\Xh$ with
${\overline{\partial}}_{\Yah}\Xh=[\Yah,\Xh]^{1,0}$.\\

\noindent{Proof.} Let $\nabla$ be a metric connection on $T^{\C}M$, then we have, for any $X,Y,Z\in\Gamma(T^{\C}M)$ :
\begin{eqnarray}\label{e9}
(\nab{X}{}Y^{\flat})(Z)&=&\frac{1}{2}\Bigl(Xg(Y,Z)-Zg(X,Y)+Yg(X,Z)-g([X,Z],Y)+g([X,Y],Z)-g([Y,Z],X)\nonumber\\
&+&g(T(X,Y),Z)-g(T(X,Z),Y)-g(T(Y,Z),X)\Bigr)\nonumber\\
&=&\frac{1}{2}\Bigl(dY^{\flat}(X,Z)+(\Lc_{Y}g)(X,Z)\Bigr)+\Kc^{\nabla}(X,Y,Z),
\end{eqnarray}
with $\ds\Kc^{\nabla}(X,Y,Z)=\frac{1}{2}\Bigl(g(T(X,Y),Z)-g(T(X,Z),Y)-g(T(Y,Z),X)\Bigr)$.\\ 
If $\nabla$ preserves each component of the complex Witt structure $T^{\C}M=(\Rc^{\C}\oplus{(\Rcs)}^{\C})\op(\So{1}{0}\oplus \So{0}{1})$, then, for $\Yh\in\So{1}{0}$ and $Z\in{(\So{1}{0})}^{\perp}$, we have $(\nab{X}{}{\Yh}^{\flat})(Z)=0$ and (\ref{e9}) gives the equation :
$$
d{\Yh}^{\flat}(X,Z)+(\Lc_{\Yh}g)(X,Z)+g(T(X,\Yh),Z)-g(T(X,Z),\Yh)-g(T(\Yh,Z),X)=0.
$$
Since $T^{\C}M={(\So{1}{0})}^{\perp}\oplus\So{0}{1}$, then, depending of $X\in{(\So{1}{0})}^{\perp}$ or $\So{0}{1}$, the previous equation becomes
$$g(T(X,\Yh)+[X,\Yh],Z)+g(T(Z,\Yh)+[Z,\Yh],X)-g(T(X,Z)+[X,Z],\Yh)=0,$$
for $X\in{(\So{1}{0})}^{\perp}$ such that $g(X,Z)=0$.
$$g(T(X,\Yh),Z)+g(T(Z,\Yh),X)+(\Lc_{\Yh}g)(X,Z)-g(T(X,Z)+[X,Z],\Yh)=0,$$
for $X\in{(\So{1}{0})}^{\perp}$ such that $g(X,Z)\neq 0$,\\
and,
$$g(T(\Yh,Z),\Xah)+(\Lc_{Z}g)(\Yh,\Xah)+d{Z}^{\flat}(\Xah,\Yh)-g(T(\Xah,\Yh),Z)+g(T(\Xah,Z),\Yh)=0,$$
for $\Xah\in\So{0}{1}$.\\

Now, taking respectively $Z=\Zh$, $Z=\n{}{}$ and $Z=\n{}{*}$, we obtain the following equations : 
\begin{eqnarray}\label{e10}
g({T(\Yh,\Zh)}^{1,0},\Xah)&+&(\Lc_{\Zh}g)(\Xah,\Yh)+d{\Zh}^{\flat}(\Xah,\Yh)\nonumber\\
&-&g({T(\Xah,\Yh)}^{0,1},\Zh)+g({T(\Xah,\Zh)}^{0,1},\Yh)=0\nonumber\\
g({T(\Xh,\Yh)}^{0,1}+{[\Xh,\Yh]}^{0,1},\Zh)&+&g({T(\Zh,\Yh)}^{0,1}+{[\Zh,\Yh]}^{0,1},\Xh)\nonumber\\
&-&g({T(\Xh,\Zh)}^{0,1}+{[\Xh,\Zh]}^{0,1},\Yh)=0\nonumber\\
g({T(\n{}{},\Yh)}^{0,1}+{[\n{}{},\Yh]}^{0,1},\Zh)&-&g({T(\n{}{},\Zh)}^{0,1}+{[\n{}{},\Zh]}^{0,1},\Yh)\nonumber\\
&+&\sig{}{*}(T(\Zh,\Yh))-d\sig{}{*}(\Zh,\Yh)=0\nonumber\\
g({T(\n{}{*},\Yh)}^{0,1}+{[\n{}{*},\Yh]}^{0,1},\Zh)&-&g({T(\n{}{*},\Zh)}^{0,1}+{[\n{}{*},\Zh]}^{0,1},\Yh)\nonumber\\
&+&\sig{}{}(T(\Zh,\Yh))-d\sig{}{}(\Zh,\Yh)=0\nonumber\\
\sig{}{*}(T(\n{}{},\Yh))-(\Lc_{\n{}{}}\sig{}{*})(\Yh)=0&,& \sig{}{}(T(\n{}{*},\Yh))-(\Lc_{\n{}{*}}\sig{}{})(\Yh)=0\nonumber\\
-g({T(\n{}{},\Yh)}^{1,0},\Xah)&-&g({T(\n{}{},\Xah)}^{0,1},\Yh)+(\Lc_{\n{}{}}g)(\Yh,\Xah)\nonumber\\
&+&\sig{}{*}(T(\Yh,\Xah))-d\sig{}{*}(\Yh,\Xah)=0\nonumber\\
-g({T(\n{}{*},\Yh)}^{1,0},\Xah)&-&g({T(\n{}{*},\Xah)}^{0,1},\Yh)+(\Lc_{\n{}{*}}g)(\Yh,\Xah)\nonumber\\
&+&\sig{}{}(T(\Yh,\Xah))-d\sig{}{}(\Yh,\Xah)=0\nonumber\\
\sig{}{}(T(\n{}{},\Yh))&+&\sig{}{*}(T(\n{}{*},\Yh))+(\Lc_{\Yh}g)(\n{}{},\n{}{*})\nonumber\\
&+&g({T(\n{}{},\n{}{*})}^{0,1}+{[\n{}{},\n{}{*}]}^{0,1},\Yh)=0.
\end{eqnarray}
Now we take 
\begin{eqnarray*}
g({T(\Yh,\Zh)}^{1,0},\Xah)&=&-(\Lc_{\Zh}g)(\Xah,\Yh)-d{\Zh}^{\flat}(\Xah,\Yh)\\
{T(\Xah,\Yh)}^{0,1}&=&0,\quad {T(\Xh,\Yh)}^{0,1}=-{[\Xh,\Yh]}^{0,1}\\
{T(\n{}{(*)},\Yh)}^{0,1}&=&-{[\n{}{(*)},\Yh]}^{0,1},\quad{T(\n{}{},\n{}{*})}^{0,1}=-{[\n{}{},\n{}{*}]}^{0,1}\\
g({T(\n{}{(*)},\Yh)}^{1,0},\Xah)&=&g({T(\n{}{(*)},\Xah)}^{0,1},\Yh)=\frac{1}{2}(\Lc_{\n{}{(*)}}g)(\Yh,\Xah)
\end{eqnarray*}
\begin{eqnarray*}
\sig{}{*}(T(\Zh,\Yh))&=&d\sig{}{*}(\Zh,\Yh),\quad \sig{}{}(T(\Zh,\Yh))=d\sig{}{}(\Zh,\Yh)\\
\sig{}{*}(T(\Yh,\Xah))&=&d\sig{}{*}(\Yh,\Xah),\quad \sig{}{}(T(\Yh,\Xah))=d\sig{}{}(\Yh,\Xah)\\
\sig{}{*}(T(\n{}{},\Yh))&=&(\Lc_{\n{}{}}\sig{}{*})(\Yh),\quad \sig{}{}(T(\n{}{*},\Yh))=(\Lc_{\n{}{*}}\sig{}{})(\Yh)\\
\sig{}{*}(T(\n{}{*},\Yh))&=&\sig{}{}(T(\n{}{},\Yh))=-\frac{1}{2}(\Lc_{\Yh}g)(\n{}{},\n{}{*})=\frac{1}{2}\bigl((\Lc_{\n{}{}}\sig{}{})(\Yh)+(\Lc_{\n{}{*}}\sig{}{*})(\Yh)\bigr).
\end{eqnarray*}
Hence equations (\ref{e10}) are satisfy. By an easy calculation, we obtain that 
$$(\Lc_{\Zh}g)(\Xah,\Yh)+d{\Zh}^{\flat}(\Xah,\Yh)=\sqrt{-1}(d\omega)(\Yh,\Zh,\Xah)$$ 
and 
\begin{eqnarray*}
 g([\n{}{(*)},\Yh]^{0,1},\Zh)&=&\frac{1}{2}g((J\Lc_{\n{}{(*)}}J)(\Yh),\Zh)\\
&=&-\frac{1}{2}(\Lc_{\n{}{(*)}}g)(\Yh,\Zh)-\frac{\sqrt{-1}}{2}(\Lc_{\n{}{(*)}}\omega)(\Yh,\Zh).
\end{eqnarray*}
Morever we choose $\sig{}{*}(T(\n{}{},\n{}{*}))=\sig{}{}(T(\n{}{},\n{}{*}))=0$. Note that $T(\Xah,\Yah)$ and $T(\n{}{(*)},\Yah)$ are obtained by conjugating. The formula (\ref{e8}) for $T$ follows from $T(X,Y)={T(X,Y)}^{1,0}+{T(X,Y)}^{0,1}+\sig{}{}(T(X,Y))\n{}{}+\sig{}{*}(T(X,Y))\n{}{*}$. Conversely, by means of formula (\ref{e8}) for $T$ and (\ref{e9}), we define a pseudo-Riemannian connection $\nabla$ with the required properties. $\Box$\\

Now, considering the restriction of the Chern-Robinson connection $\nabla$ to $TM$, we obtain similary to proposition 6.1 of \cite{Pe2} the :

\begin{Cor} 
Let $(M^{2m+2},g,\Nc,\Rc,\Rcs)$ be an almost Robinson manifold. There exists a connection $\nabla$ on $TM$ (called the real Chern-Robinson connection) such that :
\begin{enumerate}
\item $\nabla$ preserves the Witt structure $TM=(\Rc\oplus{\Rc}^{*})\op\Sc$.
\item $\nabla g=\nabla J=0$ (hence $\nabla\omega=0$).
\item The torsion $T$ of $\nabla$ is given by   
\begin{eqnarray*}
T(X,Y)&=&-\frac{1}{4}N_J(X,Y)+\frac{1}{2}\Bigl({(d^{c}\omega)}^{+}-\Mc({(d^{c}\omega)}^{+})\Bigr)^{\sharp}(X,Y)+d\sig{}{}(X,Y)\n{}{}+d\sig{}{*}(X,Y)\n{}{*}\\
T(\n{}{},X)&=&(\Lc_{\n{}{}}\sig{}{*})(X)\n{}{*}-\frac{1}{2}(\Lc_{X}g)(\n{}{},\n{}{*})\n{}{}+\tu{\Sc}(\n{}{},X)+\frac{1}{4}\Bigl(\Lc_{\n{}{}}\omega-\Mc(\Lc_{\n{}{}}\omega)\Bigr)^{\sharp}(JX)\\
T(\n{}{*},X)&=&(\Lc_{\n{}{*}}\sig{}{})(X)\n{}{}-\frac{1}{2}(\Lc_{X}g)(\n{}{},\n{}{*})\n{}{*}+\tu{\Sc}(\n{}{*},X)+\frac{1}{4}\Bigl(\Lc_{\n{}{*}}\omega-\Mc(\Lc_{\n{}{*}}\omega)\Bigr)^{\sharp}(JX)\\
T(\n{}{},\n{}{*})&=&-[\n{}{},\n{}{*}]_{\Sc},
\end{eqnarray*}
with $X,Y\in\Sc$,   
\begin{eqnarray*}
N_J(X,Y)&=&[X,Y]_{\Sc}-[JX,JY]_{\Sc}+J[JX,Y]_{\Sc}+J[X,JY]_{\Sc}\\
(d^{c}\omega)(X,Y,Z)&=&-(J^{*}d\omega)(X,Y,Z)\\
g(\tu{\Sc}(\n{}{(*)},X),Y)&=&\frac{1}{2}(\Lc_{\n{}{(*)}}g)(X,Y),
\end{eqnarray*}
and, for $\alpha\in\Omega^{p}(M)$ and $\beta\in\Omega^{3}(M)$
\begin{eqnarray*}
g(\alpha^{\sharp}(X_1,X_2,\ldots,X_{p-1}),X_{p})&=&\alpha(X_1,X_2,\ldots,X_{p-1},X_{p})\\
\Mc(\alpha)(X_1,X_2,X_3,\ldots,X_{p})&=&\alpha(JX_1,JX_2,X_3,\ldots,X_{p})
\end{eqnarray*}
$$
\beta^{+}(X_1,X_2,X_3)=\frac{1}{4}\Bigl(3\beta(X_1,X_2,X_3)+\beta(X_1,JX_2,JX_3)
+\beta(JX_1,X_2,JX_3)+\beta(JX_1,JX_2,X_3)\Bigr).
$$
\end{enumerate}
\end{Cor}

\begin{Rem} If $d\omega=0$, note that the Chern-Robinson connection and the canonical Witt connection (\cite{Pe2})
coincide.
\end{Rem}

For an almost Fefferman-Robinson manifold, we have, by proposition 2.1, that $\Lc_{\n{}{}}\sig{}{*}=\Lc_{\n{}{*}}\sig{}{*}=0$, $[\n{}{},\n{}{*}]_{\Sc}=0$, ${(d\sig{}{*})}_{\Sc}^{2,0}=0$, $\Tc^{1,0}(\Xh,\Yh)=\Tc^{1,0}(\n{}{(*)},\Xh)=0$ and $(\Lc_{\Xh}g)(\n{}{},\n{}{*})=-(\Lc_{\n{}{}}\sig{}{})(\Xh)$. Also we obtain :
 
\begin{Cor} 
Let $(M^{2m+2},g,\Nc,\Rc,\Rcs)$ be an almost Fefferman-Robinson manifold. Then torsion of 
the Chern-Robinson connection is given by :
\begin{eqnarray}\label{e11}
T(\Xh,\Yh)&=&-[\Xh,\Yh]^{0,1}+d\sig{}{}(\Xh,\Yh)\n{}{}\nonumber\\
T(\Xh,\Yah)&=&d\sig{}{}(\Xh,\Yah)\n{}{}+\omega(\Xh,\Yah)\n{}{*}\nonumber\\
T(\n{}{},\Xh)&=&-[\n{}{},\Xh]^{0,1}+\frac{1}{2}(\Lc_{\n{}{}}\sig{}{})(\Xh)\n{}{},\hspace{1cm} T(\n{}{},\n{}{*})=0\nonumber\\
T(\n{}{*},\Xh)&=&-[\n{}{*},\Xh]^{0,1}+(\Lc_{\n{}{*}}\sig{}{})(\Xh)\n{}{}+\frac{1}{2}(\Lc_{\n{}{}}\sig{}{})(\Xh)\n{}{*},
\end{eqnarray}
with $$-g([\n{}{(*)},\Xh]^{0,1},\Yh)=\frac{1}{2}(\Lc_{\n{}{(*)}}g)(\Xh,\Yh).$$
\end{Cor}

\noindent Curvature and first Bianchi identity.\\
\\
If $E\to M$ is a vector bundle over a manifold $M$, then we denote by 
$\Omega^{k}(M;E)=\Gamma(\wedge^{k}{TM}^*\otimes E)$ the bundle of $E$-valued $k$-forms on $M$. 
For $\varphi\in\Omega^{k}(M;\wedge^{r}E^*)$ and $\psi\in\Omega^{l}(M;\wedge^{s}E^*)$, 
we define $\varphi\wedge\psi\in\Omega^{k+l}(M;\wedge^{r+s}E^*)$ by :
$$(\varphi\wedge\psi)(X_1,\ldots,X_{k+l})=\frac{1}{k!l!}\sum_{\pi\in\mathfrak{S}_{k+l}}\epsilon(\pi)\varphi(X_{\pi(1)},\ldots,X_{\pi(k)})\wedge\psi(X_{\pi(k+1)},\ldots,X_{\pi(k+l)}).$$

Now, assume that $E$ and $TM$ are endowed with connections $\nabla^{E}$ and $\nabla$, then we recall that the covariant derivative and the exterior covariant derivative on $\Omega^{k}(M;E)$ are respectively defined by :

\begin{eqnarray*}
(\nabla_{X}\varphi)(X_{1},\ldots,X_{k})&=&\nabla_{X}^{E}\varphi(X_{1},\ldots,X_{k})-\sum_{i=1}^{k}\varphi(X_{1},\ldots,\nabla_{X}X_{i},\ldots,X_{k})\\
(d^{\nabla^{E}}\varphi)(X_{1},\ldots,X_{k+1})&=&\sum_{i=1}^{k+1}{(-1)}^{i+1}\nabla_{X_{i}}^{E}\varphi(X_{1},\ldots,\hat{X_{i}},\ldots,X_{k+1})\\                                           
                                           &+&\sum_{i<j}{(-1)}^{i+j}\varphi([X_{i},X_{j}],X_{1},\ldots,\hat{X_{i}},\ldots,\hat{X_{j}},\ldots,X_{k+1})\\
                                           &=&\sum_{i=1}^{k+1}{(-1)}^{i+1}(\nabla_{X_{i}}\varphi)(X_{1},\ldots,\hat{X_{i}},\ldots,X_{k+1})\\
                                           &+&\sum_{i<j}{(-1)}^{i+j+1}\varphi(T(X_{i},X_{j}),X_{1},\ldots,\hat{X_{i}},\ldots,\hat{X_{j}},\ldots,X_{k+1}).
\end{eqnarray*}
Recall that the curvature and the torsion of a connection $\nabla$ on $TM$ satisfy the first Bianchi identity $b(R)=d^{\nabla}T$ (i.e. $R(X,Y)Z+R(Z,X)Y+R(Y,Z)X=(d^{\nabla}T)(X,Y,Z)$, $X,Y,Z\in TM$).\\ 

Now, let $(M^{2m+2},g,\Nc,\Rc,\Rcs)$ be an almost Fefferman-Robinson manifold endowed with its Chern-Robinson connection $\nabla$. Then we derive the first Bianchi for the curvature tensor $R$ of $\nabla$.

\begin{Prop} 
Let $(M^{2m+2},g,\Nc,\Rc,\Rcs)$ be an almost Fefferman-Robinson manifold endowed with its Chern-Robinson connection $\nabla$ and $R$ be the curvature tensor of $\nabla$ (i.e. $R(X,Y,Z,W)=g(R(X,Y)Z,W)$). Then we have
\begin{eqnarray*}
R(\Xh,\Yah,\Zh,\Wah)&-&R(\Zh,\Yah,\Xh,\Wah)=-g(N(N(\Xh,\Zh),\Yah),\Wah)\\
                                         &-&d\sig{}{}(\Xh,\Zh)A_{\n{}{}}(\Yah,\Wah)\\
R(\Xah,\Yah,\Zh,\Wah)&=&-g((\nab{\Zh}{}N)(\Xah,\Yah),\Wah)\\
                                         &+&d\sig{}{}(\Zh,\Xah)A_{\n{}{}}(\Yah,\Wah)-d\sig{}{}(\Zh,\Yah)A_{\n{}{}}(\Xah,\Wah)\\
                                         &+&\omega(\Zh,\Xah)A_{\n{}{*}}(\Yah,\Wah)-\omega(\Zh,\Yah)A_{\n{}{*}}(\Xah,\Wah)\\
R(\n{}{(*)},\Yah,\Zh,\Wah)&=&(\nab{\Zh}{}A_{\n{}{(*)}})(\Yah,\Wah)-g(N([\n{}{(*)},\Zh]^{0,1},\Yah),\Wah)\\
                        &-&(\Lc_{\n{}{(*)}}\sig{}{})(\Zh)A_{\n{}{}}(\Yah,\Wah)\\
R(\n{}{},\n{}{*},\Zh,\Wah)&=&A_{\n{}{}}([\n{}{*},\Zh]^{0,1},\Wah)-A_{\n{}{*}}([\n{}{},\Zh]^{0,1},\Wah)\\
R(\Xh,\Yah,\n{}{},\n{}{*})&=&-\frac{1}{2}d(\Lc_{\n{}{}}\sig{}{})(\Xh,\Yah)+(\nab{\n{}{}}{}d\sig{}{})(\Xh,\Yah)\\
&+&d\sig{}{}([\n{}{},\Yah]^{1,0},\Xh)-d\sig{}{}([\n{}{},\Xh]^{0,1},\Yah)+\frac{1}{2}\omega(\Xh,\Yah)(\Lc_{\n{}{}}\sig{}{})(\n{}{*})\\
R(\Xah,\Yah,\n{}{},\n{}{*})&=&-\frac{1}{2}d(\Lc_{\n{}{}}\sig{}{})(\Xah,\Yah)+(\nab{\n{}{}}{}d\sig{}{})(\Xah,\Yah)\\
&+&d\sig{}{}([\n{}{},\Yah]^{1,0},\Xah)-d\sig{}{}([\n{}{},\Xah]^{1,0},\Yah)\\
R(\n{}{*},\Yah,\n{}{},\n{}{*})&=&-\frac{1}{2}(\Lc_{\n{}{*}}(\Lc_{\n{}{}}\sig{}{}))(\Yah)+(\Lc_{\n{}{}}(\Lc_{\n{}{*}}\sig{}{}))(\Yah)+\frac{1}{2}(\Lc_{\n{}{}}\sig{}{})(\Yah)(\Lc_{\n{}{}}\sig{}{})(\n{}{*})\\
R(\n{}{},\Yah,\n{}{},\n{}{*})&=&\frac{1}{2}(\Lc_{\n{}{}}(\Lc_{\n{}{}}\sig{}{}))(\Yah)
\end{eqnarray*}
\begin{eqnarray}\label{e12}
g((\nab{\Xah}{}N)(\Yah,\Zah),\Wah)+g((\nab{\Zah}{}N)(\Xah,\Yah),\Wah)&+&g((\nab{\Yah}{}N)(\Zah,\Xah),\Wah)\nonumber\\
                                         &=&d\sig{}{}(\Xah,\Yah)A_{\n{}{}}(\Zah,\Wah)\nonumber\\&+&d\sig{}{}(\Zah,\Xah)A_{\n{}{}}(\Yah,\Wah)\nonumber\\&+&d\sig{}{}(\Yah,\Zah)A_{\n{}{}}(\Xah,\Wah)\nonumber\\
(\nab{\Yah}{}A_{\n{}{(*)}})(\Zah,\Wah)-(\nab{\Zah}{}A_{\n{}{(*)}})(\Yah,\Wah)
&=&-g((\nab{\n{}{(*)}}{}N)(\Yah,\Zah),\Wah)\nonumber\\ 
&+&(\Lc_{\n{}{(*)}}\sig{}{})(\Yah)A_{\n{}{}}(\Zah,\Wah)\nonumber\\&-&(\Lc_{\n{}{(*)}}\sig{}{})(\Zah)A_{\n{}{}}(\Yah,\Wah)\nonumber\\ 
(\nab{\n{}{}}{}A_{\n{}{*}})(\Zah,\Wah)-(\nab{\n{}{*}}{}A_{\n{}{}})(\Zah,\Wah)&=&-(\Lc_{\n{}{}}\sig{}{})(\n{}{*})A_{\n{}{}}(\Zah,\Wah),
\end{eqnarray}
with
$$
A_{\n{}{(*)}}=\frac{1}{2}{(\Lc_{\n{}{(*)}}g)}_{\Sc^{\C}},\quad N(\Xh,\Yh)=[\Xh,\Yh]^{0,1},\quad N(\Xah,\Yah)=[\Xah,\Yah]^{1,0}.
$$
\end{Prop}

\noindent{Proof.}
For any $X,Y,Z\in T^{\C}M$ and $\Wah\in\So{0}{1}$ we have, by the first Bianchi identity, that :
\begin{eqnarray}\label{e13}
R(X,Y,Z,\Wah)&+&R(Z,X,Y,\Wah)+R(Y,Z,X,\Wah)\\
&=&g((\nab{X}{}T)(Y,Z),\Wah)+g((\nab{Z}{}T)(X,Y),\Wah)+g((\nab{Y}{}T)(Z,X),\Wah)\nonumber\\
&+&g(T(T(X,Y),Z),\Wah)+g(T(T(Z,X),Y),\Wah)+g(T(T(Y,Z),X),\Wah).\nonumber
\end{eqnarray}
Taking $\Xh,\Zh\in\So{1}{0},\Yah\in\So{0}{1}$ and using $R(\Zh,\Xh,\Yah,\Wah)=0$ together with $T$ given by (\ref{e11}), we obtain the first formula by (\ref{e13}). Since the proof for the others formulas is similar, we omit it. $\Box$

\begin{Rem}
Note that it follows from the first Bianchi identity that the components of $R$ excepted $R(\Xh,\Yah,\Zh,\Wah)$ and $R(\n{}{},\n{}{*},\n{}{},\n{}{*})$ are determined by torsion. We note also that, due to torsion, $R$ restricted to $\Sc^{\C}$ is not a Kahler-like curvature tensor.
\end{Rem}

In that follows, we construct from $R$ a Kahler-like curvature tensor on $\Sc^{\C}$  (directly inspired by the construction given by Matsuo in (\cite{Ma}).\\

Let $(V,J)$ be a complex vector space and ${\wedge}^{2,+}V^{*}$ be the space of $J$-invariant antisymmetric $(0,2)$-tensor on $(V,J)$. For $Q\in\bigotimes^{2}({\wedge}^{2}V^{*})$, we define $Q^{+}\in\bigotimes^{2}({\wedge}^{2,+}V^{*})$ by :
$$Q^{+}(X,Y,Z,W)=\frac{1}{4}\bigl(Q(X,Y,Z,W)+Q(JX,JY,Z,W)+Q(X,Y,JZ,JW)+Q(JX,JY,JZ,JW)\bigr).$$

\begin{Def} (Kahler-like curvature tensor) 
Let $(V,J)$ be a complex vector space, then $Q\in\bigotimes^{2}({\wedge}^{2}V^{*})$ will be called a Kahler-like curvature tensor if 
$$b(Q)=0\quad\mathrm{and}\quad Q=Q^{+},$$
with $$b(Q)(X,Y,Z,W)=Q(X,Y,Z,W)+Q(Z,X,Y,W)+Q(Y,Z,X,W).$$
\end{Def}

A Kahler-like curvature tensor will be denoted by $Q^{Kl}$. Note that $Q^{Kl}\in \bigodot^{2}({\wedge}^{2,+}V^{*})\cap Ker\,b$.\\For $Q\in\bigotimes^{2}({\wedge}^{2}V^{*})$ with $(V,J)$ complex vector space then the following formula (cf. \cite{Ma}) gives a Kahler-like curvature tensor $Q^{Kl}$ :
\begin{eqnarray}\label{e14}
Q^{Kl}(X,Y,Z,W)
&=&\frac{1}{8}\Bigl(2Q^{+}(X,Y,Z,W)+2Q^{+}(Z,W,X,Y)+Q^{+}(X,Z,Y,W)-Q^{+}(Y,Z,X,W)\nonumber\\
&&-Q^{+}(X,W,Y,Z)+Q^{+}(Y,W,X,Z)+Q^{+}(JX,Z,JY,W)-Q^{+}(JY,Z,JX,W)\nonumber\\
&&-Q^{+}(JX,W,JY,Z)+Q^{+}(JY,W,JX,Z)\Bigr).
\end{eqnarray}

Let $(M^{2m+2},g,\Nc,\Rc,\Rcs)$ be an almost Robinson manifold endowed with its Chern-Robinson connection $\nabla$ and $R$ be its curvature tensor. Consider $R_{\Sc^{\C}}\in\bigotimes^{2}({\wedge}^{2}{(\Sc^{\C})}^{*})$ and define $R^{\Cc}:=R_{\Sc^{\C}}^{Kl}$ called the Chern-Robinson curvature tensor. It follows from (\ref{e14}) that $R^{\Cc}_{/\So{1}{0}\wedge\So{1}{0}}=R^{\Cc}_{/\So{0}{1}\wedge\So{0}{1}}=0$, so $R^{\Cc}\in\bigodot^{2}(\wedge^{1,1}{(\Sc^{\C})}^*)$ and $R^{\Cc}$ is given, for $\Xh,\Zh\in\Gamma(\So{1}{0})$ and $\Yah,\Wah\in\Gamma(\So{0}{1})$ by :
\begin{eqnarray}\label{e15}
R^{\Cc}(\Xh,\Yah,\Zh,\Wah)&=&\frac{1}{4}\Bigl(R(\Xh,\Yah,\Zh,\Wah)+R(\Zh,\Yah,\Xh,\Wah)\nonumber\\
&+&R(\Xh,\Wah,\Zh,\Yah)+R(\Zh,\Wah,\Xh,\Yah)\Bigr).
\end{eqnarray}

In the following, we call adapted basis of $S^{\C}$ a local basis $(Z_1\,\ldots,Z_m,\overline{Z_1},\ldots,\overline{Z_m})$ of $S^{\C}=\So{1}{0}\oplus\So{0}{1}$ such that :
$$g(Z_i,Z_j)=g(\overline{Z_i},\overline{Z_j})=0\;\mathrm{and}\; g(Z_i,\overline{Z_j})=\delta_{ij}.$$

Now, we define
\begin{eqnarray}\label{e16}
Ric^{\Cc}(\Xh,\Yah)&=&\sum_{i=1}^{m}R^{\Cc}(\Xh,\overline{Z_i},Z_i,\Yah)=\sum_{i=1}^{m}R^{\Cc}(Z_i,\overline{Z_i},\Xh,\Yah)\nonumber\\
s^{\Cc}&=&\tr\,Ric^{\Cc}=\sum_{i=1}^{m}Ric^{\Cc}(Z_i,\overline{Z_i})\nonumber\\
P^{\Cc}(\Xh,\Yah)&=&\frac{1}{m+2}\Bigl(Ric^{\Cc}(\Xh,\Yah)-\frac{s^{\Cc}}{2(m+1)}g(\Xh,\Yah)\Bigr)\nonumber\\
&=&\frac{1}{m+2}\Bigl(Ric^{\Cc}(\Xh,\Yah)+\frac{\sqrt{-1}s^{\Cc}}{2(m+1)}\omega(\Xh,\Yah)\Bigr)\nonumber\\
C^{\Mc}(\Xh,\Yah,\Zh,\Wah)&=&R^{\Cc}(\Xh,\Yah,\Zh,\Wah)\nonumber\\
&-&\Bigl(P^{\Cc}(\Xh,\Yah)g(\Zh,\Wah)+P^{\Cc}(\Zh,\Wah)g(\Xh,\Yah)\nonumber\\
&+&P^{\Cc}(\Zh,\Yah)g(\Xh,\Wah)+P^{\Cc}(\Xh,\Wah)g(\Zh,\Yah)\Bigr)\nonumber\\
&=&R^{\Cc}(\Xh,\Yah,\Zh,\Wah)\nonumber\\
&+&\sqrt{-1}\Bigl(P^{\Cc}(\Xh,\Yah)\omega(\Zh,\Wah)+P^{\Cc}(\Zh,\Wah)\omega(\Xh,\Yah)\nonumber\\
&+&P^{\Cc}(\Zh,\Yah)\omega(\Xh,\Wah)+P^{\Cc}(\Xh,\Wah)\omega(\Zh,\Yah)\Bigr),
\end{eqnarray}
with $(Z_1\,\ldots,Z_m,\overline{Z_1},\ldots,\overline{Z_m})$ adapted basis of $S^{\C}$.\\

The tensors $Ric^{\Cc}$, $s^{\Cc}$, $P^{\Cc}$ and $C^{\Mc}$ are respectively called Chern-Ricci-Robinson tensor, Chern-Robinson scalar curvature, Chern-Robinson-Schouten tensor and Chern-Moser-Robinson tensor. Note that $C^{\Mc}=0$ when $m=2$.

\begin{Def} (Einstein-Robinson manifold)\\
An almost Robinson manifold $(M^{2m+2},g,\Nc,\Rc,\Rcs)$ endowed with its Chern-Robinson connection $\nabla$ will be called an Einstein-Robinson manifold if $Ric^{\Cc}=f\omega$ for $f\in\Cc^{\infty}(M)$.
\end{Def}

Note that if $M$ is Einstein-Robinson then $f=-\dfrac{\sqrt{-1}s^{\Cc}}{m}$ and $P^{\Cc}=-\dfrac{\sqrt{-1}s^{\Cc}}{2m(m+1)}\omega$.

\section{Optical conformal invariant tensor on almost Fefferman-Robinson manifolds}

Let $(M^{2m+2},g,\Nc,\Rc,\Rcs)$ be an almost Robinson manifold and $(\n{}{},\n{}{*})\in\Gamma(\Rc\oplus\Rcs)$ be a $g$-optical pairing together with $(\sig{}{},\sig{}{*})$ optical forms associated. For $\ld\in\Cc^{\infty}(M)$, the vector fields $\n{\ld}{},\n{\ld}{*}$ and the metric $\gl$ on $M$ are defined (cf. Proposition 4.8 of \cite{FLTC2}) by :
\begin{eqnarray}\label{e17}
\n{\ld}{}&=&\n{}{}\nonumber\\
\n{\ld}{*}&=&e^{-2\ld}\bigl(\n{}{*}-2J{(d\ld)}^{\sharp}_{\Sc}-2\dlsn\n{}{}\bigr)\nonumber\\
\gl&=&e^{2\ld}g.
\end{eqnarray}
Note that $(\n{\ld}{},\n{\ld}{*})$ is a $\gl$-optical pairing with optical forms $(\sig{\ld}{},\sig{\ld}{*})$ given by :
\begin{eqnarray}\label{e18}
\sig{\ld}{}&=&{(\n{\ld}{*})}^{\flat_{\ld}}=\sig{}{}+2J^{*}\dls-2\dlsn\sig{}{*}\nonumber\\
\sig{\ld}{*}&=&{(\n{\ld}{})}^{\flat_{\ld}}=e^{2\ld}\sig{}{*},
\end{eqnarray}
with $Y^{\flat_{\ld}}(X)=\gl(Y,X)$ and $J^{*}\dls(X)=\dls(JX)$.\\

Let ${\Rc}_{\ld}=span(\n{\ld}{}),\;{\Rc}^{*}_{\ld}=span(\n{\ld}{*})\;\mathrm{and}\;\scl=\mathrm{Ker}\,\sig{\ld}{}\cap\mathrm{Ker}\,\sig{\ld}{*}$, then
$$TM=({\Rc}_{\ld}\oplus{\Rc}^{*}_{\ld})\stackrel{\perp_{\ld}}{\oplus}\scl.$$
For $X\in\Gamma(TM)$, then (\ref{e5}) yields the splittings $X=X_{\Sc}+\sig{}{}(X)\n{}{}+\sig{}{*}(X)\n{}{*}=X_{\scl}+\sig{\ld}{}(X)\n{\ld}{}+\sig{\ld}{*}(X)\n{\ld}{*}$. It follows from (\ref{e17}) and (\ref{e18}) that $X_{\scl}$ and $X_{\Sc}$ are related by :
$$
X_{\scl}=X_{\Sc}+2\sig{}{*}(X)J{(d\ld)}^{\sharp}_{\Sc}-2\dls\bigl(J(X_{\Sc}+2\sig{}{*}(X)J{(d\ld)}^{\sharp}_{\Sc})\bigr)\n{}{}.
$$
Now $\scl=\{X_{\scl}=X_{\Sc}-2\dls(JX_{\Sc})\n{}{};\;X_{\Sc}\in\Gamma(\Sc)\}$. If $\Jl$ is the endomorphism on $\scl$ given by $\Jl(X_{\scl}=X_{\Sc}-2\dls(JX_{\Sc})\n{}{})=JX_{\Sc}+2\dls(X_{\Sc})\n{}{}$, then $(\gl,\Jl)$ is an almost Hermitian structure on $\scl$ (in the following $\Jl$ is extended by $0$ on ${\Rc}_{\ld}\oplus{\Rc}^{*}_{\ld}$). Let $\Sl{1}{0}$ be the $\sqrt{-1}$-eigenbundle of $\Jl$ on $\scl^{\C}$ and $\Nc_{\ld}:=\Rc_{\ld}^{\C}\oplus\Sl{1}{0}$ then $\Nc_{\ld}$ is an almost Robinson structure on $(M^{2m+2},\gl)$ with gradings
$$T^{\C}M=({\Rc}_{\ld}^{\C}\oplus{({\Rc}^{*}_{\ld})}^{\C})\stackrel{\perp_{\ld}}{\oplus}(\Sl{1}{0}\oplus\Sl{0}{1}).$$
Let $\Xhl:=X_{\scl}^{1,0}\in\Gamma(\Sl{1}{0})$ and $\Xahl:=X_{\scl}^{0,1}\in\Gamma(\Sl{0}{1})$ then we have the following transformation rules (cf. Proposition 4.8 of \cite{FLTC2}) called optical conformal transformations :

\begin{equation}\label{e19}
\begin{pmatrix}\n{\ld}{}\cr\n{\ld}{*}\cr\Xhl\cr\Xahl
\end{pmatrix}=\begin{pmatrix}\n{}{}\cr e^{-2\ld}\bigl(\n{}{*}-2J{(d\ld)}^{\sharp}_{\Sc}-2\dlsn\n{}{}\bigr)\cr\Xh-2\sqrt{-1}\dls(\Xh)\nu{}{}\cr\Xah+2\sqrt{-1}\dls(\Xah)\nu{}{}
\end{pmatrix}.
\end{equation}

For $\alpha\in\Omega^{1}(M)$, then (\ref{e5}) yields $\alpha=\alpha_{\Sc}+\sig{}{}\otimes\alpha(\n{}{})+\sig{}{*}\otimes\alpha(\n{}{*})=
\alpha_{\scl}+\sig{\ld}{}\otimes\alpha(\n{\ld}{})+\sig{\ld}{*}\otimes\alpha(\n{\ld}{*})$ and we have the following relation :
$$
\alpha_{\scl}=\alpha_{\Sc}+\bigl(4\dlsn\sig{}{*}-2J^{*}\dls\bigr)\otimes\alpha(\n{}{})+2\sig{}{*}\otimes\alpha(J{(d\ld)}^{\sharp}_{\Sc}).
$$
Also, for any $\ld\in\Cc^{\infty}(M)$ such that $d\ld(\nu{}{})=0$, then it follows from the previous equation that $\dlsl=\dls$ and $\Jl^{*}\dlsl=J^{*}\dls-2\dlsn\sig{}{*}$.\\

In the following $\Cc_{\nu{}{}}^{\infty}(M)=\{\ld\in\Cc^{\infty}(M);\;d\ld(\nu{}{})=0\}$.

\begin{Le} Let $(M^{2m+2},g,\Nc,\Rc,\Rcs)$ be an almost Robinson manifold, $\ld\in\Cc_{\nu{}{}}^{\infty}(M)$ and $(\n{\ld}{},\n{\ld}{*})$ be a $g_{\ld}$-optical pairing. Let $\ol=e^{2\ld}(\omega+2d\ld\wedge\sig{}{*})$. Then 
$$\n{\ld}{}\lrcorner\,\ol=\n{\ld}{*}\lrcorner\,\ol=0,\quad \ol=\gl\circ\Jl$$
and 
$$d\ol=e^{2\ld}(d\omega+2d\ld\wedge(\omega-d\sig{}{*})),\quad \Lc_{\n{\ld}{}}\ol=e^{2\ld}(\Lc_{\n{}{}}\omega+2d\ld\wedge(\Lc_{\n{}{}}\sig{}{*})).$$
\end{Le}

\noindent{Proof.} 
We have 
$$\n{\ld}{}\lrcorner\,\ol=e^{2\ld}\n{}{}\lrcorner\,(\omega+2d\ld\wedge\sig{}{*})=e^{2\ld}(\n{}{}\lrcorner\,\omega+2d\ld(\n{}{})\sig{}{*}-2\sig{}{*}(\n{}{})d\ld)=0.$$
Since 
$\n{\ld}{*}=e^{-2\ld}\bigl(\n{}{*}-2J{(d\ld)}^{\sharp}_{\Sc}-2\dlsn\n{}{}\bigr)$,
we obtain 
\begin{eqnarray*}
\n{\ld}{*}\lrcorner\,\ol&=&(\n{}{*}-2J{(d\ld)}^{\sharp}_{\Sc}-2\dlsn\n{}{})\lrcorner\,(\omega+2d\ld\wedge\sig{}{*})\\
&=&-2(J{(d\ld)}^{\sharp}_{\Sc})\lrcorner\,\omega+2d\ld(\n{}{*})\sig{}{*}-2
d\ld\\
&=&2\dls-2(d\ld-\sig{}{*}d\ld(\n{}{*}))=2\dls-2\dls=0.
\end{eqnarray*}
Now, let $X_{\ld}=X-2\dls(JX)\n{}{}\in\Gamma(\scl)$ and $Y_{\ld}=Y-2\dls(JY)\n{}{}\in\Gamma(\scl)$ with $X,Y\in\Gamma(\Sc)$. Since $J_{\ld}X_{\ld}=JX+2\dls(X)\n{}{}$, we have :
$$\gl(\Jl X_{\ld},Y_{\ld})=e^{2\ld}g(JX,Y)=e^{2\ld}\omega(X,Y)=\ol(X_{\ld},Y_{\ld}).$$
Using $\Jl\n{\ld}{}=\Jl\n{\ld}{*}=0$ and $\n{\ld}{}\lrcorner\,\ol=\n{\ld}{*}\lrcorner\,\ol=0$, we deduce that $\ol=\gl\circ\Jl$.\\
Now
$$d\ol=2e^{2\ld}d\ld\wedge(\omega+2d\ld\wedge\sig{}{*})+e^{2\ld}(d\omega-2d\ld\wedge d\sig{}{*})=e^{2\ld}(d\omega+2d\ld\wedge(\omega-d\sig{}{*})).$$
Since 
$$\Lc_{\n{\ld}{}}\ol=\n{\ld}{}\lrcorner\,(d\ol)+d(\n{\ld}{}\lrcorner\,\ol)=\n{}{}\lrcorner\,(d\ol),$$
then 
$$\Lc_{\n{\ld}{}}\ol=e^{2\ld}\n{}{}\lrcorner\,(d\omega+2d\ld\wedge(\omega-d\sig{}{*}))=e^{2\ld}(\Lc_{\n{}{}}\omega+2d\ld\wedge(\Lc_{\n{}{}}\sig{}{*})).\quad\Box$$
Note that, if $\omega=d\sig{}{*}$, then $\ol=d\sig{\ld}{*}$.

\begin{Prop} Let $(M^{2m+2},g,\Nc,\Rc,\Rcs)$ be a geodetic almost Robinson manifold and $\ld\in\Cc_{\nu{}{}}^{\infty}(M)$. Then, for any $\Xhl,\Yhl\in\Gamma(\Sl{1}{0})$ and $\Xahl,\Yahl,\Zahl\in\Gamma(\Sl{0}{1})$, we have :
\begin{eqnarray}\label{e20}
(\nab{\Xhl}{\ld}{\Yhl}^{\flat_{\ld}})(\Zahl)&=&e^{2\ld}(\nab{\Xhl}{}{\Yhl}^{\flat})(\Zahl)+2\Bigl(\dlsl(\Xhl)\gl(\Yhl,\Zahl)+\dlsl(\Yhl)d\sig{\ld}{*}(\Xhl,\Jl\Zahl)\Bigr)\nonumber\\
&-&\dlsl(\Xhl)(\Lc_{\n{}{}}\ol)(\Yhl,\Zahl)\nonumber\\
&+&\dlsl(\Jl\Zahl)\Bigl((\Lc_{\n{}{}}\gl)(\Xhl,\Yhl)+\sqrt{-1}(\Lc_{\n{}{}}\ol)(\Xhl,\Yhl)\Bigr)\nonumber\\
(\nab{\Xahl}{\ld}{\Yhl}^{\flat_{\ld}})(\Zahl)&=&e^{2\ld}(\nab{\Xahl}{}{\Yhl}^{\flat})(\Zahl)-2\dlsl(\Zahl)d\sig{\ld}{*}(\Yhl,\Jl\Xahl)\nonumber\\
&-&\dlsl(\Xahl)(\Lc_{\n{}{}}\ol)(\Yhl,\Zahl)\nonumber\\
&-&\dlsl(\Jl\Yhl)\Bigl((\Lc_{\n{}{}}\gl)(\Xahl,\Zahl)-\sqrt{-1}(\Lc_{\n{}{}}\ol)(\Xahl,\Zahl)\Bigr)\nonumber\\
(\nab{\n{\ld}{}}{\ld}{\Yhl}^{\flat_{\ld}})(\Zahl)&=&e^{2\ld}(\nab{\n{}{}}{}{\Yhl}^{\flat})(\Zahl)\nonumber\\
(\nab{\n{\ld}{*}}{\ld}{\Yhl}^{\flat_{\ld}})(\Zahl)&=&e^{2\ld}(\nab{\n{\ld}{*}}{}{\Yhl}^{\flat})(\Zahl)+\gl(\Yhl,\Zahl)d\ld(\n{\ld}{*})\nonumber\\
&-&2\Bigl(\dlsl(\Jl\Yhl)\dlsl(\Zahl)-\dlsl(\Jl\Zahl)\dlsl(\Yhl)\Bigr)\nonumber\\
&-&\dlsl(\Jl\Yhl)\Bigl(\Lc_{\n{\ld}{*}}\sig{\ld}{*}(\Zahl)-\gl([\n{}{},\n{\ld}{*}]^{{(1,0)}_{\ld}}_{\scl},\Zahl)\Bigr)\nonumber\\
&+&\dlsl(\Jl\Zahl)\Bigl(\Lc_{\n{\ld}{*}}\sig{\ld}{*}(\Yhl)-\gl([\n{}{},\n{\ld}{*}]^{{(0,1)}_{\ld}}_{\scl},\Yhl)\Bigr)\nonumber\\
&+&(\Lc_{\dlsd{\sharp}}\omega)(\Yhl,\Zahl)+2\dlsl(\Jl\Yhl)(\Lc_{\dlsd{\sharp}}\omega)(\n{}{},\Zahl)\nonumber\\
&-&2\dlsl(\Jl\Zahl)(\Lc_{\dlsd{\sharp}}\omega)(\n{}{},\Yhl).
\end{eqnarray}
\end{Prop}

For a strongly geodetic almost Fefferman-Robinson manifold, the previous proposition becomes : 

\newpage

\begin{Cor} Let $(M^{2m+2},g,\Nc,\Rc,\Rcs)$ be a strongly geodetic almost Fefferman-Robinson manifold and $\ld\in\Cc_{\nu{}{}}^{\infty}(M)$. Then, for any $\Xhl,\Yhl\in\Gamma(\Sl{1}{0})$ and $\Xahl\in\Gamma(\Sl{0}{1})$, we have :
\begin{eqnarray}\label{e21}
\nab{\Xhl}{\ld}\Yhl&=&{\bigl(\nab{\Xhl}{}\Yhl\bigr)}^{{(1,0)}_{\ld}}_{\scl}+2\Bigl(\dlsl(\Xhl)\Yhl+\dlsl(\Yhl)\Xhl\Bigr)\nonumber\\
\nab{\Xahl}{\ld}\Yhl&=&{\bigl(\nab{\Xahl}{}\Yhl\bigr)}^{{(1,0)}_{\ld}}_{\scl}-2g(\Xahl,\Yhl)\dlsd{\sharp\,{(1,0)}_{\ld}}\nonumber\\
\nab{\n{\ld}{}}{\ld}\Yhl&=&{\bigl(\nab{\n{\ld}{}}{}\Yhl\bigr)}^{{(1,0)}_{\ld}}_{\scl}\nonumber\\
\nab{\n{\ld}{*}}{\ld}\Yhl&=&{\bigl(\nab{\n{\ld}{*}}{}\Yhl\bigr)}^{{(1,0)}_{\ld}}_{\scl}
+2e^{-2\ld}\Bigl(\nab{\Yhl}{}\Jl\dlsd{\sharp\,{(1,0)}_{\ld}}-2\dlsl(\Yhl)\Jl\dlsd{\sharp\,{(1,0)}_{\ld}}\Bigr).
\end{eqnarray}
\end{Cor}
\hspace{0.3cm}\\
Proofs of Proposition 4.1 and Corollary 4.1 will be given in annexe.

\begin{Prop} Let $(M^{2m+2},g,\Nc,\Rc,\Rcs)$ be a strongly geodetic almost Fefferman-Robinson manifold and $\ld\in\Cc_{\nu{}{}}^{\infty}(M)$. Then, for any $\Xhl,\Zhl\in\Gamma(\Sl{1}{0})$ and $\Yahl,\Wahl\in\Gamma(\Sl{0}{1})$, we have :
\begin{eqnarray*}
g(R^{\ld}(\Xhl,\Yahl)\Zhl,\Wahl)&=&g(R(\Xh,\Yah)\Zh,\Wah)\nonumber\\
&-2&\Bigl((\nab{\Yah}{}\dls)(\Xh)g(\Zh,\Wah)+(\nab{\Zh}{}\dls)(\Wah)g(\Xh,\Yah)\nonumber\\
&+&(\nab{\Yah}{}\dls)(\Zh)g(\Xh,\Wah)+(\nab{\Xh}{}\dls)(\Wah)g(\Zh,\Yah)\Bigr)\nonumber\\
&-&2\dlsn\Bigl(g(\Xh,\Yah)g(\Zh,\Wah)+g(\Zh,\Yah)g(\Xh,\Wah)\Bigr)\\
R^{\Cc}_{\ld}(\Xhl,\Yahl,\Zhl,\Wahl)&=&e^{2\ld}\Bigl[R^{\Cc}(\Xh,\Yah,\Zh,\Wah)\\
&-&\Bigl((\nab{\Xh}{sym}\dls)(\Yah)g(\Zh,\Wah)+(\nab{\Zh}{sym}\dls)(\Wah)g(\Xh,\Yah)\\
&+&(\nab{\Yah}{sym}\dls)(\Zh)g(\Xh,\Wah)+(\nab{\Xh}{sym}\dls)(\Wah)g(\Zh,\Yah)\Bigr)\\
&-&2\dlsn\Bigl(g(\Xh,\Yah)g(\Zh,\Wah)+g(\Zh,\Yah)g(\Xh,\Wah)\Bigr)\Bigr]\\
Ric^{\Cc}_{\ld}(\Xhl,\Yahl)&=&Ric^{\Cc}(\Xh,\Yah)-(m+2)(\nab{\Xh}{sym}\dls)(\Yah)\\
&-&g(\Xh,\Yah)\bigl(\tr(\nab{}{sym}\dls)+2(m+1)\dlsn\bigr)\\
s^{\Cc}_{\ld}&=&e^{-2\ld}\Bigl(s^{\Cc}-2(m+1)\tr(\nab{}{sym}\dls)-2m(m+1)\dlsn\Bigr)\\
P^{\Cc}_{\ld}(\Xhl,\Yahl)&=&P^{\Cc}(\Xh,\Yah)-(\nab{\Xh}{sym}\dls)(\Yah)-\dlsn g(\Xh,\Yah),\\
\end{eqnarray*}
and, for $m\geq 2$,
\begin{equation}\label{e22}
C^{\Mc}_{\ld}(\Xhl,\Yahl,\Zhl,\Wahl)=e^{2\ld}C^{\Mc}(\Xh,\Yah,\Zh,\Wah),
\end{equation}
with $(\nab{X}{sym}\dls)(Y)=(\nab{X}{}\dls)(Y)+(\nab{Y}{}\dls)(X)$ and $\ds \tr(\nab{}{sym}\dls)=\sum_{i=1}^{m}(\nab{Z_i}{sym}\dls)(\overline{Z_i})$.
\end{Prop}

It follows from ($\ref{e22}$) that

\begin{Cor} If $(M^{2m+2},g,\Nc,\Rc,\Rcs)$ is a strongly geodetic almost Fefferman-Robinson manifold with $m\geq 2$ then $C^{\Mc}$ is an optical conformal invariant tensor.
\end{Cor}

\noindent{Proof.} We have 
\begin{eqnarray*}
R^{\ld}(\Xhl,\Yahl)\Zhl&=&\Bigl(\nab{\Xhl}{\ld}\nab{\Yahl}{\ld}-\nab{\Yahl}{\ld}\nab{\Xhl}{\ld}-\nab{[\Xhl\Yahl]}{\ld}\Bigr)\Zhl\\
&=&\Bigl(\nab{\Xhl}{\ld}\nab{\Yahl}{\ld}-\nab{\nab{\Xhl}{\ld}\Yahl}{\ld}-\nab{\Yahl}{\ld}\nab{\Xhl}{\ld}+\nab{\nab{\Yahl}{\ld}\Xhl}{\ld}+\nab{T^{\lambda}(\Xhl,\Yahl)}{\ld}\Bigr)\Zhl\\
&=&B^{\ld}(\Xhl,\Yahl)\Zhl+\nab{T^{\ld}(\Xhl,\Yahl)}{\ld}\Zhl.
\end{eqnarray*}
Now, using (\ref{e21}), we have 
\begin{eqnarray}\label{e23}
\nab{\Xhl}{\ld}\nab{\Yahl}{\ld}\Zhl&=&{\Bigl(\nab{\Xhl}{}{\bigl(\nab{\Yahl}{}\Zhl\bigr)}^{{(1,0)}_{\ld}}_{\scl}\Bigr)}^{{(1,0)}_{\ld}}_{\scl}-2g(\Zhl,\Yahl){\bigl(\nab{\Xhl}{}\dlsd{\sharp\,{(1,0)}_{\ld}}\bigr)}^{{(1,0)}_{\ld}}_{\scl}\nonumber\\
&+&2\Bigl(\dlsl(\Xhl){\bigl(\nab{\Yahl}{}\Zhl\bigr)}^{{(1,0)}_{\ld}}_{\scl}
+\dlsl\bigl({\bigl(\nab{\Yahl}{}\Zhl\bigr)}^{{(1,0)}_{\ld}}_{\scl}\bigr)\Xhl\Bigr)\nonumber\\
&-&2\Bigl(\Xhl g(\Zhl,\Yahl)+2g(\Zhl,\Yahl)\dlsl(\Xhl)\Bigr)\dlsd{\sharp\,{(1,0)}_{\ld}}\nonumber\\
&-&2g(\Zhl,\Yahl)\dlsn\Xhl.
\end{eqnarray}
\begin{eqnarray}\label{e24}
\nab{\nab{\Xhl}{\ld}\Yahl}{\ld}\Zhl&=&{\Bigl(\nab{{\bigl(\nab{\Xhl}{}\Yahl\bigr)}^{{(0,1)}_{\ld}}_{\scl}}
{}\Zhl\Bigr)}^{{(1,0)}_{\ld}}_{\scl}-2g(\Xhl,\Yahl){\bigl(\nab{\dlsd{\sharp\,{(0,1)}_{\ld}}}{}\Zhl\bigr)}^{{(1,0)}_{\ld}}_{\scl}\nonumber\\
&-&2\Bigl(g(\Zhl,{\bigl(\nab{\Xhl}{}\Yahl\bigr)}^{{(0,1)}_{\ld}}_{\scl})-2g(\Xhl,\Yahl)\dlsl(\Zhl)\Bigr)\dlsd{\sharp\,{(1,0)}_{\ld}}.
\end{eqnarray}
\begin{eqnarray}\label{e25}
\nab{\Yahl}{\ld}\nab{\Xhl}{\ld}\Zhl&=&{\Bigl(\nab{\Yahl}{}{\bigl(\nab{\Xhl}{}\Zhl\bigr)}^{{(1,0)}_{\ld}}_{\scl}\Bigr)}^{{(1,0)}_{\ld}}_{\scl}\nonumber\\
&+&2\Bigl(\dlsl(\Xhl){\bigl(\nab{\Yahl}{}\Zhl\bigr)}^{{(1,0)}_{\ld}}_{\scl}+\dlsl(\Zhl){\bigl(\nab{\Yahl}{}\Xhl\bigr)}^{{(1,0)}_{\ld}}_{\scl}\Bigr)\nonumber\\
&-&2\Bigl(g({\bigl(\nab{\Xhl}{}\Zhl\bigr)}^{{(1,0)}_{\ld}}_{\scl},\Yahl)+2g(\Xhl,\Yahl)\dlsl(\Zhl)\nonumber\\
&&+2g(\Zhl,\Yahl)\dlsl(\Xhl)\Bigr)\dlsd{\sharp\,{(1,0)}_{\ld}}\nonumber\\
&+&2\Bigl(\Yahl\dlsl(\Xhl)\Zhl+\Yahl\dlsl(\Zhl)\Xhl\Bigr).
\end{eqnarray}
\begin{eqnarray}\label{e26}
\nab{\nab{\Yahl}{\ld}\Xhl}{\ld}\Zhl&=&{\Bigl(\nab{{\bigl(\nab{\Yahl}{}\Xhl\bigr)}^{{(1,0)}_{\ld}}_{\scl}}
{}\Zhl\Bigr)}^{{(1,0)}_{\ld}}_{\scl}-2g(\Xhl,\Yahl){\bigl(\nab{\dlsd{\sharp\,{(1,0)}_{\ld}}}{}\Zhl\bigr)}^{{(1,0)}_{\ld}}_{\scl}\nonumber\\
&+&2\Bigl(\dlsl\bigl({\bigl(\nab{\Yahl}{}\Xhl\bigr)}^{{(1,0)}_{\ld}}_{\scl}\bigr)\Zhl+\dlsl(\Zhl){\bigl(\nab{\Yahl}{}\Xhl\bigr)}^{{(1,0)}_{\ld}}_{\scl}\Bigr)\nonumber\\
&-&4g(\Xhl,\Yahl)\Bigl(\dlsl(\Zhl)\dlsd{\sharp\,{(1,0)}_{\ld}}+\frac{1}{2}\dlsn\Zhl\Bigr).
\end{eqnarray}
Since 
${\bigl(\nab{X}{}\Zhl\bigr)}^{{(1,0)}_{\ld}}_{\scl}=\nab{X}{}\Zhl+(\nab{X}{}\sig{\ld}{})(\Zhl)\n{\ld}{}$
then we obtain 
\begin{eqnarray}\label{e27}
{\Bigl(\nab{{\bigl(\nab{\Yahl}{}\Xhl\bigr)}^{{(1,0)}_{\ld}}_{\scl}-{\bigl(\nab{\Xhl}{}\Yahl\bigr)}^{{(0,1)}_{\ld}}_{\scl}}{}\Zhl\Bigr)}^{{(1,0)}_{\ld}}_{\scl}&=&{\Bigl(\nab{\nab{\Yahl}{}\Xhl-\nab{\Xhl}{}\Yahl}{}\Zhl\Bigr)}^{{(1,0)}_{\ld}}_{\scl}\\
&+&\Bigl((\nab{\Yahl}{}\sig{\ld}{})(\Xhl)-(\nab{\Xhl}{}\sig{\ld}{})(\Yahl)\Bigr){\bigl(\nab{\n{}{}}{}\Zhl\bigr)}^{{(1,0)}_{\ld}}_{\scl}\nonumber.
\end{eqnarray}
Using (\ref{e23}),(\ref{e24}),(\ref{e25}),(\ref{e26}) and (\ref{e27}), we obtain that :
\begin{eqnarray*}
B^{\ld}(\Xhl,\Yahl)\Zhl&=&{\Bigl(\nab{\Xhl}{}{\bigl(\nab{\Yahl}{}\Zhl\bigr)}^{{(1,0)}_{\ld}}_{\scl}\Bigr)}^{{(1,0)}_{\ld}}_{\scl}-{\Bigl(\nab{\Yahl}{}{\bigl(\nab{\Xhl}{}\Zhl\bigr)}^{{(1,0)}_{\ld}}_{\scl}\Bigr)}^{{(1,0)}_{\ld}}_{\scl}\\
&+&{\Bigl(\nab{\nab{\Yahl}{}\Xhl-\nab{\Xhl}{}\Yahl}{}\Zhl\Bigr)}^{{(1,0)}_{\ld}}_{\scl}
+\Bigl((\nab{\Yahl}{}\sig{\ld}{})(\Xhl)-(\nab{\Xhl}{}\sig{\ld}{})(\Yahl)\Bigr){\bigl(\nab{\n{}{}}{}\Zhl\bigr)}^{{(1,0)}_{\ld}}_{\scl}\\
&+&2\sqrt{-1}g(\Xhl,\Yahl){\bigl(\nab{\Jl\dlsd{\sharp}}{}\Zhl\bigr)}^{{(1,0)}_{\ld}}_{\scl}\\
&-&2\Bigl((\nab{\Yahl}{}\dlsl)(\Xhl)\Zhl+(\nab{\Yahl}{}\dlsl)(\Zhl)\Xhl\\&+&g(\Zhl,\Yahl){\bigl(\nab{\Xhl}{}\dlsd{\sharp\,{(1,0)}_{\ld}}\bigr)}^{{(1,0)}_{\ld}}_{\scl}\Bigr)\\
&-&2\Bigl(g(\Xhl,\Yahl)\dlsn\Zhl+g(\Zhl,\Yahl)\dlsn\Xhl\\&+&2g(\Xhl,\Yahl)\dlsl(\Zhl)\dlsd{\sharp\,{(1,0)}_{\ld}}\Bigr).
\end{eqnarray*}
Since $g\Bigl({\Bigl(\nab{X}{}{\bigl(\nab{Y}{}\Zhl\bigr)}^{{(1,0)}_{\ld}}_{\scl}\Bigr)}^{{(1,0)}_{\ld}}_{\scl},\Wahl\Bigr)=g(\nab{X}{}(\nab{Y}{}\Zhl),\Wahl)$, we obtain that :
\begin{eqnarray}\label{e28}
g(B^{\ld}(\Xhl,\Yahl)\Zhl,\Wahl)&=&g(B(\Xhl,\Yahl)\Zhl,\Wahl)\nonumber\\
&+&g(\nab{2d\sig{}{*}(\Xhl,\Yahl)\Jl\dlsd{\sharp}+((\nab{\Yahl}{}\sig{\ld}{})(\Xhl)-(\nab{\Xhl}{}\sig{\ld}{})(\Yahl))\n{}{}}{}\Zhl,\Wahl)\nonumber\\
&-&2\Bigl((\nab{\Yahl}{}\dlsl)(\Xhl)g(\Zhl,\Wahl)+(\nab{\Yahl}{}\dlsl)(\Zhl)g(\Xhl,\Wahl)\nonumber\\
&+&(\nab{\Xhl}{}\dlsl)(\Wahl)g(\Zhl,\Yahl)\Bigr)\nonumber\\
&-&2\dlsn\Bigl(g(\Xhl,\Yahl)g(\Zhl,\Wahl)+g(\Zhl,\Yahl)g(\Xhl,\Wahl)\Bigr)\nonumber\\
&-&4g(\Xhl,\Yahl)\dlsl(\Zhl)\dlsl(\Wahl),
\end{eqnarray}
with $B(\Xhl,\Yahl)=\nab{\Xhl}{}\nab{\Yahl}{}-\nab{\nab{\Xhl}{}\Yahl}{}-\nab{\Yahl}{}\nab{\Xhl}{}+\nab{\nab{\Yahl}{}\Xhl}{}$.\\
Now since $\om{\ld}=d\sig{\ld}{*}$, we have 
$$T^{\ld}(\Xhl,\Yahl)=d\sig{\ld}{}(\Xhl,\Yahl)\n{\ld}{}+d\sig{\ld}{*}(\Xhl,\Yahl)\n{\ld}{*}.$$
Hence,
$$\nab{T^{\ld}(\Xhl,\Yahl)}{\ld}\Zhl=d\sig{\ld}{}(\Xhl,\Yahl)\nab{\n{\ld}{}}{\ld}\Zhl+d\sig{\ld}{*}(\Xhl,\Yahl)\nab{\n{\ld}{*}}{\ld}\Zhl.$$
Now, we have by (\ref{e21}) :
\begin{eqnarray*}
d\sig{\ld}{*}(\Xhl,\Yahl)\nab{\n{\ld}{*}}{\ld}\Zhl&=&d\sig{\ld}{*}(\Xhl,\Yahl){\bigl(\nab{\n{\ld}{*}}{}\Zhl\bigr)}^{{(1,0)}_{\ld}}_{\scl}\\
&+&d\sig{}{*}(\Xhl,\Yahl)\Bigl(2\nab{\Zhl}{}\Jl\dlsd{\sharp\,{(1,0)}_{\ld}}-4\dlsl(\Zhl)\Jl\dlsd{\sharp\,{(1,0)}_{\ld}}\Bigr).
\end{eqnarray*}
Hence, 
\begin{eqnarray}\label{e29}
g(\nab{T^{\ld}(\Xhl,\Yahl)}{\ld}\Zhl,\Wahl)&=&d\sig{\ld}{}(\Xhl,\Yahl)g(\nab{\n{}{}}{}\Zhl,\Wahl)+d\sig{\ld}{*}(\Xhl,\Yahl)g(\nab{\n{\ld}{*}}{}\Zhl,\Wahl)\nonumber\\
&-&\sqrt{-1}d\sig{}{*}(\Xhl,\Yahl)\Bigl(4\dlsl(\Zhl)\dlsl(\Wahl)-2g(\nab{\Zhl}{}\dlsd{\sharp\,{(1,0)}_{\ld}},\Wahl)\Bigr)\nonumber\\
&=&d\sig{\ld}{}(\Xhl,\Yahl)g(\nab{\n{}{}}{}\Zhl,\Wahl)+d\sig{\ld}{*}(\Xhl,\Yahl)g(\nab{\n{\ld}{*}}{}\Zhl,\Wahl)\nonumber\\
&+&g(\Xhl,\Yahl)\Bigl(4\dlsl(\Zhl)\dlsl(\Wahl)-2(\nab{\Zhl}{}\dlsl)(\Wahl)\Bigr).
\end{eqnarray}
By adding (\ref{e28}) and (\ref{e29}), we obtain that 
\begin{eqnarray}\label{e30}
g(R^{\ld}(\Xhl,\Yahl)\Zhl,\Wahl)&=&g(B(\Xhl,\Yahl)\Zhl,\Wahl)+g(\nab{d\sig{\ld}{*}(\Xhl,\Yahl)\n{\ld}{*}}{}\Zhl,\Wahl)\nonumber\\
&+&g(\nab{2d\sig{}{*}(\Xhl,\Yahl)\Jl\dlsd{\sharp}+(d\sig{\ld}{}(\Xhl,\Yahl)-(\nab{\Xhl}{}\sig{\ld}{})(\Yahl)+(\nab{\Yahl}{}\sig{\ld}{})(\Xhl))\n{}{}}{}\Zhl,\Wahl)\nonumber\\
&-&2\Bigl((\nab{\Yahl}{}\dlsl)(\Xhl)g(\Zhl,\Wahl)+(\nab{\Zhl}{}\dlsl)(\Wahl)g(\Xhl,\Yahl)\nonumber\\
&+&(\nab{\Yahl}{}\dlsl)(\Zhl)g(\Xhl,\Wahl)+(\nab{\Xhl}{}\dlsl)(\Wahl)g(\Zhl,\Yahl)\Bigr)\nonumber\\
&-&2\dlsn\Bigl(g(\Xhl,\Yahl)g(\Zhl,\Wahl)+g(\Zhl,\Yahl)g(\Xhl,\Wahl)\Bigr).
\end{eqnarray}
Now, since $\Xhl=\Xh-2\sqrt{-1}\dls(\Xh)\nu{}{}$ and $\Yahl=\Yah+2\sqrt{-1}\dls(\Yah)\nu{}{}$ with $\Xh\in\Gamma(\So{1}{0}),\;\Yah\in\Gamma(\So{0}{1})$, we obtain using (\ref{e11}) :
\begin{eqnarray*}
T(\Xhl,\Yahl)&=&\Bigl(d\sig{}{}(\Xh,\Yah)-\sqrt{-1}\dls(\Xh)(\Lc_{\n{}{}}\sig{}{})(\Yh)-\sqrt{-1}\dls(\Yh)(\Lc_{\n{}{}}\sig{}{})(\Xh)\Bigr)\n{}{}\nonumber\\
&+&d\sig{}{*}(\Xh,\Yah)\n{}{*}.
\end{eqnarray*}
We deduce that $\sig{\ld}{*}(T(\Xhl,\Yahl))=d\sig{\ld}{*}(\Xhl,\Yahl)$ and
${T(\Xhl,\Yahl)}_{\scl}=2d\sig{}{*}(\Xhl,\Yahl)\Jl\dlsd{\sharp}$.
Hence
\begin{eqnarray*}
T(\Xhl,\Yahl)&=&{T(\Xhl,\Yahl)}_{\scl}+\sig{\ld}{}(T(\Xhl,\Yahl))\n{\ld}{}+\sig{\ld}{*}(T(\Xhl,\Yahl))\n{\ld}{*}\\
&=&2d\sig{}{*}(\Xhl,\Yahl)\Jl\dlsd{\sharp}+(d\sig{\ld}{}(\Xhl,\Yahl)-(\nab{\Xhl}{}\sig{\ld}{})(\Yahl)+(\nab{\Yahl}{}\sig{\ld}{})(\Xhl))\n{}{}\\
&+&d\sig{\ld}{*}(\Xhl,\Yahl)\n{\ld}{*}.
\end{eqnarray*}
By substituting in (\ref{e30}) and noticing that 
$$g(R(\Xhl,\Yahl)\Zhl,\Wahl)=g(B(\Xhl,\Yahl)\Zhl,\Wahl)+g(\nab{T(\Xhl,\Yahl)}{}\Zhl,\Wahl)$$ 
we obtain 
\begin{eqnarray}\label{e31}
g(R^{\ld}(\Xhl,\Yahl)\Zhl,\Wahl)&=&g(R(\Xhl,\Yahl)\Zhl,\Wahl)\\
&-&2\Bigl((\nab{\Yahl}{}\dlsl)(\Xhl)g(\Zhl,\Wahl)+(\nab{\Zhl}{}\dlsl)(\Wahl)g(\Xhl,\Yahl)\nonumber\\
&+&(\nab{\Yahl}{}\dlsl)(\Zhl)g(\Xhl,\Wahl)+(\nab{\Xhl}{}\dlsl)(\Wahl)g(\Zhl,\Yahl)\Bigr)\nonumber\\
&-&2\dlsn\Bigl(g(\Xhl,\Yahl)g(\Zhl,\Wahl)+g(\Zhl,\Yahl)g(\Xhl,\Wahl)\Bigr)\nonumber.
\end{eqnarray}
Now,
\begin{eqnarray*}
g(R(\Xhl,\Yahl)\Zhl,\Wahl)&=&g(R(\Xh,\Yah)\Zh,\Wah)-2\sqrt{-1}\dls(\Xh)g(R(\n{}{},\Yah)\Zh,\Wah)\\
&+&2\sqrt{-1}\dls(\Yah)g(R(\n{}{},\Xh)\Wah,\Zh).
\end{eqnarray*}
Under the assumptions, we have using (\ref{e12}) that $g(R(\n{}{},\Yah)\Zh,\Wah)=g(R(\n{}{},\Xh)\Wah,\Zh)=0$.
We deduce that :
\begin{equation}\label{e32}
g(R(\Xhl,\Yahl)\Zhl,\Wahl)=g(R(\Xh,\Yah)\Zh,\Wah).
\end{equation}
Now, since $d\alpha(X,Y)=(\nab{X}{}\alpha)(Y)-(\nab{Y}{}\alpha)(X)+\alpha(T(X,Y))$, then apply this formula to $\alpha=\dls=d\ld-\sig{}{*}\otimes d\ld(\n{}{*})$ together with $d\ld(\nu{}{})=0$ and ${T(\nu{}{},X_{\Sc})}_{\Sc}=0$, yields to $(\nab{\nu{}{}}{}\dls)(X_{\Sc})=0$. It follows that $(\nab{\Xhl}{}\dlsl)(\Yahl)=(\nab{\Xh}{}\dls)(\Yah)$.
We deduce from (\ref{e31}) and (\ref{e32}) that : 
\begin{eqnarray}\label{e33}
g(R^{\ld}(\Xhl,\Yahl)\Zhl,\Wahl)&=&g(R(\Xh,\Yah)\Zh,\Wah)\\
&-2&\Bigl((\nab{\Yah}{}\dls)(\Xh)g(\Zh,\Wah)+(\nab{\Zh}{}\dls)(\Wah)g(\Xh,\Yah)\nonumber\\
&+&(\nab{\Yah}{}\dls)(\Zh)g(\Xh,\Wah)+(\nab{\Xh}{}\dls)(\Wah)g(\Zh,\Yah)\Bigr)\nonumber\\
&-&2\dlsn\Bigl(g(\Xh,\Yah)g(\Zh,\Wah)+g(\Zh,\Yah)g(\Xh,\Wah)\Bigr)\nonumber.
\end{eqnarray}
Now we have 
\begin{eqnarray}\label{e34}
R^{\Cc}_{\ld}(\Xhl,\Yahl,\Zhl,\Wahl)&=&\frac{1}{4}\Bigl(\gl(R^{\ld}(\Xhl,\Yahl)\Zhl,\Wah)+\gl(R^{\ld}(\Zhl,\Yahl)\Xhl,\Wahl)\nonumber\\
&+&\gl(R^{\ld}(\Xhl,\Wahl)\Zhl,\Yahl)+\gl(R^{\ld}(\Zhl,\Wahl)\Xhl,\Yahl)\Bigr)\nonumber\\
&=&\frac{e^{2\ld}}{4}\Bigl(g(R^{\ld}(\Xhl,\Yahl)\Zhl,\Wah)+g(R^{\ld}(\Zhl,\Yahl)\Xhl,\Wahl)\nonumber\\
&+&g(R^{\ld}(\Xhl,\Wahl)\Zhl,\Yahl)+g(R^{\ld}(\Zhl,\Wahl)\Xhl,\Yahl)\Bigr).
\end{eqnarray}
Substituting (\ref{e33}) and its permutations in (\ref{e34}) and using  
$$(\nab{X}{sym}\dls)(Y)=(\nab{X}{}\dls)(Y)+(\nab{Y}{}\dls)(X)$$ 
yields the formula for $R^{\Cc}_{\ld}(\Xhl,\Yahl,\Zhl,\Wahl)$.\\
Now, setting $Z_{\ld}^i=e^{-\ld}(Z_i-2\sqrt{-1}\dls(Z_i)\n{}{})$ with $g(Z_i,\overline{Z_j})=\delta_{ij}$, we have 
\begin{eqnarray*}
Ric^{\Cc}_{\ld}(\Xhl,\Yahl)&=&\sum_{i=1}^{m}R^{\Cc}_{\ld}(\Xhl,\Yahl,Z_{\ld}^i,\overline{Z_{\ld}^i})\\
&=&e^{-2\ld}\sum_{i=1}^{m}R^{\Cc}_{\ld}(\Xhl,\Yahl,Z_i-2\sqrt{-1}\dls(Z_i)\n{}{},\overline{Z_i}+2\sqrt{-1}\dls(\overline{Z_i})\n{}{})\\
&=&\sum_{i=1}^{m}R^{\Cc}(\Xh,\Yah,Z_i,\overline{Z_i})\\
&-&\sum_{i=1}^{m}\Bigl((\nab{\Xh}{sym}\dls)(\Yah)g(Z_i,\overline{Z_i})+(\nab{Z_i}{sym}\dls)(\overline{Z_i})g(\Xh,\Yah)\\
&+&(\nab{\Yah}{sym}\dls)(Z_i)g(\Xh,\overline{Z_i})+(\nab{\Xh}{sym}\dls)(\overline{Z_i})g(Z_i,\Yah)\Bigr)\\
&-&2\dlsn\sum_{i=1}^{m}\Bigl(g(\Xh,\Yah)g(Z_i,\overline{Z_i})+g(Z_i,\Yah)g(\Xh,\overline{Z_i})\Bigr)\\
&=&Ric^{\Cc}(\Xh,\Yah)-(m+2)(\nab{\Xh}{sym}\dls)(\Yah)\\
&-&g(\Xh,\Yah)\bigl(\tr(\nab{}{sym}\dls)+2(m+1)\dlsn\bigr),
\end{eqnarray*}
with $\ds \tr(\nab{}{sym}\dls)=\sum_{i=1}^{m}(\nab{Z_i}{sym}\dls)(\overline{Z_i})$.
For $s^{\Cc}_{\ld}$, we have 
\begin{eqnarray*}
s^{\Cc}_{\ld}&=&e^{-2\ld}\sum_{i=1}^{m}Ric^{\Cc}_{\ld}(Z_i-2\sqrt{-1}\dls(Z_i)\n{}{},\overline{Z_i}+2\sqrt{-1}\dls(\overline{Z_i})\n{}{})\\
&=&e^{-2\ld}\Bigl(\sum_{i=1}^{m}Ric^{\Cc}(Z_i,\overline{Z_i})-(m+2)\sum_{i=1}^{m}(\nab{Z_i}{sym}\dls)(\overline{Z_i})\\
&-&\sum_{i=1}^{m}g(Z_i,\overline{Z_i})\bigl(\tr(\nab{}{sym}\dls)+2(m+1)\dlsn\bigr)\\
&=&e^{-2\ld}\Bigl(s^{\Cc}-2(m+1)\tr(\nab{}{sym}\dls)-2m(m+1)\dlsn\Bigr).
\end{eqnarray*}
We deduce that 
\begin{eqnarray}\label{e35}
P^{\Cc}_{\ld}(\Xhl,\Yahl)&=&\frac{1}{m+2}\Bigl(Ric^{\Cc}_{\ld}(\Xhl,\Yahl)-\frac{s^{\Cc}_{\ld}}{2(m+1)}\gl(\Xhl,\Yahl)\Bigr)\nonumber\\
&=&\frac{1}{m+2}\Bigl(Ric^{\Cc}(\Xh,\Yah)-\frac{s^{\Cc}}{2(m+1)}g(\Xh,\Yah)\Bigr)-(\nab{\Xh}{sym}\dls)(\Yah)\nonumber\\
&-&\dlsn g(\Xh,\Yah)\nonumber\\
&=&P^{\Cc}(\Xh,\Yah)-(\nab{\Xh}{sym}\dls)(\Yah)-\dlsn g(\Xh,\Yah).
\end{eqnarray}
Now
\begin{eqnarray}\label{e36}
C^{\Mc}_{\ld}(\Xhl,\Yahl,\Zhl,\Wahl)&=&R^{\Cc}_{\ld}(\Xhl,\Yahl,\Zhl,\Wahl)\nonumber\\
&-&P^{\Cc}_{\ld}(\Xhl,\Yahl)\gl(\Zhl,\Wahl)-P^{\Cc}_{\ld}(\Zhl,\Wahl)\gl(\Xhl,\Yahl)\nonumber\\
&-&P^{\Cc}_{\ld}(\Zhl,\Yahl)\gl(\Xhl,\Wahl)-P^{\Cc}_{\ld}(\Xhl,\Wahl)\gl(\Zhl,\Yahl).
\end{eqnarray}
The formula for $R^{\Cc}_{\ld}(\Xhl,\Yahl,\Zhl,\Wahl)$ together with (\ref{e35}) yields the result for $C^{\Mc}_{\ld}$ by substituting in (\ref{e36}). $\Box$\\

By the transformation rules (\ref{e19}), we obtain also :
\begin{Prop} Let $(M^{2m+2},g,\Nc,\Rc,\Rcs)$ be a Fefferman-Robinson manifold and $\ld\in\Cc_{\nu{}{}}^{\infty}(M)$. Then, for any $\Xhl,\Yhl\in\Gamma(\Sl{1}{0})$, we have :
$$A^{\ld}_{\n{\ld}{*}}(\Xhl,\Yhl)=A_{\n{}{*}}(\Xh,\Yh)+\sqrt{-1}\bigl((\nab{\Xh}{sym}\dls)(\Yh)-2(\dls\odot\dls)(\Xh,\Yh)\bigr).
$$
\end{Prop}

\noindent{Proof.} Since $(M^{2m+2},g,\Nc,\Rc,\Rcs)$ be a Fefferman-Robinson manifold then $(M^{2m+2},\gl,\Nc_{\ld},\Rc_{\ld},\Rcs_{\ld})$ is also a Fefferman-Robinson manifold. Hence
$\Lc_{\n{\ld}{*}}\sig{\ld}{*}=0$ and $[\n{}{},\n{\ld}{*}]_{\scl}=0$, also we have $(\Lc_{\n{\ld}{*}}\gl)(\n{}{},X)=(\Lc_{\n{\ld}{*}}\sig{\ld}{*})(X)+\gl([\n{}{},\n{\ld}{*}],X)=0$ for $X\in\Gamma(\mathrm{Ker}\,\sig{\ld}{*})$.
Now, since $\Xhl=\Xh-2\sqrt{-1}\dls(\Xh)\n{}{}$ and $\Yhl=\Yh-2\sqrt{-1}\dls(\Yh)\n{}{}$ with $\Xh,\Yh\in\Gamma(\So{1}{0})$, we have :
\begin{eqnarray*}
A^{\ld}_{\n{\ld}{*}}(\Xhl,\Yhl)&=&\frac{1}{2}(\Lc_{\n{\ld}{*}}\gl)(\Xh-2\sqrt{-1}\dls(\Xh)\n{}{},\Yh-2\sqrt{-1}\dls(\Yh)\n{}{})
=\frac{1}{2}(\Lc_{\n{\ld}{*}}\gl)(\Xh,\Yh)\\
&=&\frac{1}{2}e^{2\ld}(\Lc_{\n{\ld}{*}}g)(\Xh,\Yh)+e^{2\ld}d\ld(\n{\ld}{*})g(\Xh,\Yh)=\frac{1}{2}e^{2\ld}(\Lc_{\n{\ld}{*}}g)(\Xh,\Yh).
\end{eqnarray*}
Now we have
$$\n{\ld}{*}=e^{-2\ld}\bigl(\n{}{*}-2\sqrt{-1}{(d\ld)}^{\sharp\,(1,0)}_{\Sc}+2\sqrt{-1}{(d\ld)}^{\sharp\,(0,1)}_{\Sc}-2\dlsn\n{}{}\bigr).$$
Using $\Lc_{\mu Z}g=d\mu\odot Z^{\flat}+\mu\Lc_{Z}g$ for $\mu\in\Cc^{\infty}(M)$, we deduce that
\begin{eqnarray}\label{e37}
A^{\ld}_{\n{\ld}{*}}(\Xhl,\Yhl)&=&\frac{1}{2}(\Lc_{\n{}{*}}g)(\Xh,\Yh)-\sqrt{-1}(\Lc_{{(d\ld)}^{\sharp\,(1,0)}_{\Sc}}g)(\Xh,\Yh)+\sqrt{-1}(\Lc_{{(d\ld)}^{\sharp\,(0,1)}_{\Sc}}g)(\Xh,\Yh)\nonumber\\
&-&\dlsn(\Lc_{\n{}{}}g)(\Xh,\Yh)-2\sqrt{-1}(\dls\odot\dls)(\Xh,\Yh).
\end{eqnarray}
Using the formula
$$(\Lc_{Z}g)(X,Y)=g(\nab{X}{}Z,Y)+g(\nab{Y}{}Z,X)+g(T(Z,X),Y)+g(T(Z,Y),X)$$
together with $T$ given by formula (\ref{e11}) with $[\Xh,\Yh]^{0,1}=0$, we obtain
\begin{eqnarray*}
(\Lc_{{(d\ld)}^{\sharp\,(1,0)}_{\Sc}}g)(\Xh,\Yh)&=&
g(T({(d\ld)}^{\sharp\,(1,0)}_{\Sc},\Xh),\Yh)+g(T({(d\ld)}^{\sharp\,(1,0)}_{\Sc},\Yh),\Xh)=0
\end{eqnarray*}
and
\begin{eqnarray}\label{e38}
(\Lc_{{(d\ld)}^{\sharp\,(0,1)}_{\Sc}}g)(\Xh,\Yh)
&=&g(\nab{\Xh}{}{(d\ld)}^{\sharp\,(0,1)}_{\Sc},\Yh)+g(\nab{\Yh}{}{(d\ld)}^{\sharp\,(0,1)}_{\Sc},\Xh)\nonumber\\
&=&(\nab{\Xh}{}\dls)(\Yh)+(\nab{\Yh}{}\dls)(\Xh).
\end{eqnarray}
Since $M$ is a strongly geodetic Fefferman-Robinson manifold then $(\Lc_{\n{}{}}g)(\Xh,\Yh)=0$, hence the formula directly follows from (\ref{e37}) and (\ref{e38}). $\Box$

\section{Second Bianchi identity and applications}
In this section we assume that $(M^{2m+2},g,\Nc,\Rc,\Rcs)$ is a Fefferman-Robinson manifold endowed with its Chern-Robinson connection $\nabla$ and we derive the second Bianchi identity ($d^{\nabla}R=0$) for the curvature tensor $R$ of $\nabla$.

\begin{Prop} Let $(M^{2m+2},g,\Nc,\Rc,\Rcs)$ be a Fefferman-Robinson manifold. Then, we have :\\
$(\nab{\Xh}{}R^{\Cc})(\Yh,\Zah,\Uh,\Wah)-(\nab{\Yh}{}R^{\Cc})(\Xh,\Zah,\Uh,\Wah)$
\begin{eqnarray*}
&=&\omega(\Yh,\Wah)(\nab{\Zah}{}A_{\n{}{*}})(\Xh,\Uh)-\omega(\Xh,\Wah)(\nab{\Zah}{}A_{\n{}{*}})(\Yh,\Uh)\\
&+&\omega(\Yh,\Zah)(\nab{\Wah}{}A_{\n{}{*}})(\Xh,\Uh)-\omega(\Xh,\Zah)(\nab{\Wah}{}A_{\n{}{*}})(\Yh,\Uh)
\end{eqnarray*}
\begin{eqnarray*}
(\nab{\n{}{}}{}R^{\Cc})(\Yh,\Zah,\Uh,\Wah)&=&0\\
(\nab{\n{}{*}}{}R^{\Cc})(\Yh,\Zah,\Uh,\Wah)&=&
(\nab{\Yh,\Uh}{2}A_{\n{}{*}})(\Zah,\Wah)+(\nab{\Zah,\Wah}{2}A_{\n{}{*}})(\Yh,\Uh)\\
&+&\omega(\Uh,\Zah)\sum_{i=1}^{m}A_{\n{}{*}}(\Yh,Z_i)A_{\n{}{*}}(\overline{Z_i},\Wah)\\
&-&\omega(\Yh,\Wah)\sum_{i=1}^{m}A_{\n{}{*}}(\Uh,Z_i)A_{\n{}{*}}(\overline{Z_i},\Zah)
\end{eqnarray*}
$(\nab{\Xh,\Wah}{2}A_{\n{}{*}})(\Yh,\Uh)-(\nab{\Yh,\Wah}{2}A_{\n{}{*}})(\Xh,\Uh)$
\begin{eqnarray*}
&=&\omega(\Yh,\Wah)(\nab{\n{}{*}}{}A_{\n{}{*}})(\Xh,\Uh)-\omega(\Xh,\Wah)(\nab{\n{}{*}}{}A_{\n{}{*}})(\Yh,\Uh)\\
&+&\sum_{i=1}^{m}A_{\n{}{*}}(\Xh,Z_i)R(\Yh,\overline{Z_i},\Uh,\Wah)
-\sum_{i=1}^{m}A_{\n{}{*}}(\Yh,Z_i)R(\Xh,\overline{Z_i},\Uh,\Wah).
\end{eqnarray*}
\end{Prop}

\noindent{Proof.}
For a Fefferman-Robinson manifold, we have by (\ref{e12}) together with $A_{\n{}{}}=0$ and $N=0$ that
$$R(\Xh,\Yah,\Zh,\Wah)=R(\Zh,\Yah,\Xh,\Wah)\,\mathrm{and}\, R(\Xh,\Yah,\Zh,\Wah)=R(\Xh,\Wah,\Zh,\Yah).$$
Also $R^{\Cc}(\Xh,\Yah,\Zh,\Wah)=R(\Xh,\Yah,\Zh,\Wah)$. Applying the second Bianchi identity, we have :
\begin{eqnarray*}
(\nab{X}{}R)(Y,Z,U,W)&+&(\nab{Z}{}R)(X,Y,U,W)+(\nab{Y}{}R)(Z,X,U,W)=\\
&-&R(T(X,Y),Z,U,W)-R(T(Z,X),Y,U,W)-R(T(Y,Z),X,U,W).
\end{eqnarray*}
We have
\begin{eqnarray*}
(\nab{\Xh}{}R)(\Yh,\Zah,\Uh,\Wah)&-&(\nab{\Yh}{}R)(\Xh,\Zah,\Uh,\Wah)\\
&=&-(\nab{\Zah}{}R)(\Xh,\Yh,\Uh,\Wah)-R(T(\Xh,\Yh),\Zah,\Uh,\Wah)\\
&+&R(T(\Xh,\Zah),\Yh,\Uh,\Wah)-R(T(\Yh,\Zah),\Xh,\Uh,\Wah).
\end{eqnarray*}
Using (\ref{e11}) and (\ref{e12}), we obtain the first formula.
The second and third formula follow from 
\begin{eqnarray*}
(\nab{\n{}{(*)}}{}R)(\Yh,\Zah,\Uh,\Wah)&=&(\nab{\Yh}{}R)(\n{}{(*)},\Zah,\Uh,\Wah)-(\nab{\Zah}{}R)(\n{}{(*)},\Yh,\Uh,\Wah)\\
&-&R(T(\n{}{(*)},\Yh),\Zah,\Uh,\Wah)+R(T(\n{}{(*)},\Zah),\Yh,\Uh,\Wah)\\
&+&R(\n{}{(*)},T(\Yh,\Zah),\Uh,\Wah),
\end{eqnarray*}
and, the fourth formula follows from
\begin{eqnarray*}
(\nab{\n{}{*}}{}R)(\Yh,\Zh,\Uh,\Wah)
&+&(\nab{\Zh}{}R)(\n{}{*},\Yh,\Uh,\Wah)-(\nab{\Yh}{}R)(\n{}{*},\Zh,\Uh,\Wah)\\
&=&R(T(\n{}{*},\Zh),\Yh,\Uh,\Wah)-R(T(\n{}{*},\Yh),\Zh,\Uh,\Wah)\\
&+&R(\n{}{*},T(\Yh,\Zh),\Uh,\Wah).\quad\Box
\end{eqnarray*}

\noindent\\
Now, consider $E\to M$ complex vector bundle over $M$ and the complex vector bundle $\wedge^{p,q}{(\Sc^{\C})}^*=\wedge^{p}{(\So{1}{0})}^*\otimes\wedge^{q}{(\So{0}{1})}^*\to M$. We denote by $\Omega_{\Sc}^{p,q}(M;E)=\Gamma(\wedge^{p,q}{(\Sc^{\C})}^*\otimes E)$ the bundle of $E$-valued $(p,q)$ forms on $\Sc$.\\

For an adapted basis $(Z_1\,\ldots,Z_m,\overline{Z_1},\ldots,\overline{Z_m})$ of $S^{\C}$, we consider, for any $i\in\{1,\ldots,m\}$,\\
$Z_i^{*}=\overline{Z_i}^{\flat}\in{(\So{1}{0})}^*$ and $\overline{Z}_i^*=Z_i^{\flat}\in{(\So{0}{1})}^*$.\\ We define
$$t:\wedge^{p,q}{(\Sc^{\C})}^*\otimes\wedge^{p',q'}{(\Sc^{\C})}^*\to \wedge^{p-1,q}{(\Sc^{\C})}^*\otimes\wedge^{p'+1,q'}{(\Sc^{\C})}^*$$
and
$$\overline{t}:\wedge^{p,q}{(\Sc^{\C})}^*\otimes\wedge^{p',q'}{(\Sc^{\C})}^*\to \wedge^{p,q-1}{(\Sc^{\C})}^*\otimes\wedge^{p',q'+1}{(\Sc^{\C})}^*$$
by :
\begin{equation}\label{e39}
t(\alpha\otimes\beta)=-\sum_{i=1}^{m}(Z_i\lrcorner\alpha)\otimes(Z_i^{*}\wedge\beta)\quad\mathrm{and}\quad
\overline{t}(\alpha\otimes\beta)=-\sum_{i=1}^{m}(\overline{Z_i}\lrcorner\alpha)\otimes(\overline{Z}_i^*\wedge\beta),
\end{equation}
with $\alpha\in\wedge^{p,q}{(\Sc^{\C})}^*,\beta\in\wedge^{p',q'}{(\Sc^{\C})}^*$.

\newpage

For $\phi\in\Omega_{\Sc}^{1,1}(M;\wedge^{p',q'}{(\Sc^{\C})}^*)$, we denote by $\phi^{1,0}:=\overline{t}(\phi)\in\Omega_{\Sc}^{1,0}(M;\wedge^{p',q'+1}{(\Sc^{\C})}^*)$ and by $\phi^{0,1}:=t(\phi)\in\Omega_{\Sc}^{0,1}(M;\wedge^{p'+1,q'}{(\Sc^{\C})}^*)$.
Note that for $\phi\in\Omega_{\Sc}^{2,0}(M;\wedge^{p',q'}{(\Sc^{\C})}^*)$, $\phi^{1,0}:=t(\phi)\in\Omega_{\Sc}^{1,0}(M;\wedge^{p'+1,q'}{(\Sc^{\C})}^*)$.\\

\noindent
Now, we define
$$c:\wedge^{p,q}{(\Sc^{\C})}^*\otimes\wedge^{p',q'}{(\Sc^{\C})}^*\to(\wedge^{p-1,q}{(\Sc^{\C})}^*\otimes\wedge^{p',q'-1}{(\Sc^{\C})}^*)\oplus(\wedge^{p,q-1}{(\Sc^{\C})}^*\otimes\wedge^{p'-1,q'}{(\Sc^{\C})}^*)$$
by :
$$c(\alpha\otimes\beta)=\sum_{i=1}^{m}\bigl(Z_i\lrcorner\alpha\otimes\overline{Z_i}\lrcorner\beta+\overline{Z_i}\lrcorner\alpha\otimes Z_i\lrcorner\beta\bigr).$$
For $\phi\in\Omega_{\Sc}^{p,q}(M;\wedge^{1,1}{(\Sc^{\C})}^*)$ with $p+q=1$,
we have, for $X,Y\in\Gamma(\Sc^{\C})$ :
$$c(\phi)(X,Y)=\sum_{i}\bigl(\phi(Z_i,X)(\overline{Z_i},Y)+\phi(\overline{Z_i},X)(Z_i,Y)\bigr).$$
Let $k\in\N^*$ and $p,q,p',q'$ such that $p+q=p'+q'=k$.
For $\phi^{p,q}=\alpha^{p,q}\otimes\sigma\in\wedge^{p,q}{(\Sc^{\C})}^* \otimes\wedge^{r,s}{(\Sc^{\C})}^*,\psi^{p',q'}=\beta^{p',q'}\otimes\sigma'\in\wedge^{p',q'}{(\Sc^{\C})}^*\otimes\wedge^{r',s'}{(\Sc^{\C})}^*$, we define
$\phi^{p,q}\underset{k}{\circ}\psi^{p',q'}\in\wedge^{r,s}{(\Sc^{\C})}^*\otimes \wedge^{r',s'}{(\Sc^{\C})}^*$ by :
$$
\phi^{p,q}\underset{k}{\circ}\psi^{p',q'}=\frac{1}{k!}\bigl((\underbrace{c\circ\ldots\circ c}_{k})(\alpha^{p,q}\otimes\beta^{p',q'})\bigr)\otimes(\sigma\otimes\sigma').
$$
In the following
$$
{[\phi,\psi]}_k=\phi\underset{k}{\circ}\psi-\psi\underset{k}{\circ}\phi\quad\mathrm{and}\quad{\{\phi,\psi\}}_k=\phi\underset{k}{\circ}\psi+\psi\underset{k}{\circ}\phi.
$$

\noindent Note that, for $\phi^{1,0}\in\Omega_{\Sc}^{1,0}(M;\wedge^{r,s}{(\Sc^{\C})}^*)$ and $\psi^{0,1}\in\Omega_{\Sc}^{0,1}(M;\wedge^{r',s'}{(\Sc^{\C})}^*)$ with $r+s=r'+s'=1$, then ${[\phi^{1,0},\psi^{0,1}]}_1\in\Omega_{\Sc}^{r+r',s+s'}(M)$ and :
$${[\phi^{1,0},\psi^{0,1}]}_1=\sqrt{-1}\Lambda_{\Sc}(\phi^{1,0}\wedge\psi^{0,1}).$$

Let $\Omega\in\Omega_{\Sc}^{1,1}(M)$ and $\Omega^{1,0}=\overline{t}(\Omega)\in\Omega_{\Sc}^{1,0}(M;{(\So{0}{1})}^*)$, $\Omega^{0,1}=t(\Omega)\in\Omega_{\Sc}^{0,1}(M;{(\So{1}{0})}^*)$, then
\begin{equation}\label{e40}
{[\omega^{1,0},\Omega^{0,1}]}_{1}=\sqrt{-1}\Omega \quad\mathrm{and}\quad{[\omega^{0,1},\Omega^{1,0}]}_{1}=-\sqrt{-1}\Omega.
\end{equation}

Note also that for $\phi\in\Omega_{\Sc}^{1,1}(M;\wedge^{r,s}{(\Sc^{\C})}^*)$,
$$\omega\underset{2}{\circ}\phi=\Lambda_{\Sc}(\phi).$$

\noindent{\bf Differential operators $\del,\delb,\sdel,\sdelb$}\\

Assume that $E\to M$ is a complex vector bundle over $(M,g,\Nc,\Rc,\Rcs)$, then we define on $\Omega_{\Sc}^{p,q}(M;E)$ the following differential operators :
\begin{eqnarray*}
L_{\Sc}&:&\Omega_{\Sc}^{p,q}(M;E)\to\Omega_{\Sc}^{p+1,q+1}(M;E);\quad\alpha\mapsto L_{\Sc}\alpha=\omega\wedge\alpha\\
\Lambda_{\Sc}&:&\Omega_{\Sc}^{p,q}(M;E)\to\Omega_{\Sc}^{p-1,q-1}(M;E);\quad\alpha\mapsto\Lambda_{\Sc}\alpha=-\sqrt{-1}\sum_{i=1}^{m}\overline{Z_i}\lrcorner\,(Z_i\lrcorner\,\alpha)\\
\del&:&\Omega_{\Sc}^{p,q}(M;E)\to\Omega_{\Sc}^{p+1,q}(M;E);
\quad\alpha\mapsto\del\alpha=\sum_{i=1}^{m}Z_i^*\wedge\nab{Z_i}{}\alpha\\
\delb&:&\Omega_{\Sc}^{p,q}(M;E)\to\Omega_{\Sc}^{p,q+1}(M;E);
\quad\alpha\mapsto\delb\alpha=\sum_{i=1}^{m}\overline{Z}_i^*\wedge\nab{\overline{Z_i}}{}\alpha\\
\sdel&:&\Omega_{\Sc}^{p,q}(M;E)\to\Omega_{\Sc}^{p-1,q}(M;E);
\quad\alpha\mapsto\sdel\alpha=-\sum_{i=1}^{m}Z_i\lrcorner\,(\nab{\overline{Z_i}}{}\alpha)\\
\sdelb&:&\Omega_{\Sc}^{p,q}(M;E)\to\Omega_{\Sc}^{p,q-1}(M;E);
\quad\alpha\mapsto\sdelb\alpha=-\sum_{i=1}^{m}\overline{Z_i}\lrcorner\,(\nab{Z_i}{}\alpha)\\
\Box_{\Sc}&:&\Omega_{\Sc}^{p,q}(M;E)\to\Omega_{\Sc}^{p,q}(M;E);
\quad\alpha\mapsto\Box_{\Sc}\alpha=(\del\sdel+\sdel\del)\alpha\\\\
\overline{\Box}_{\Sc}&:&\Omega_{\Sc}^{p,q}(M;E)\to\Omega_{\Sc}^{p,q}(M;E);
\quad\alpha\mapsto\overline{\Box}_{\Sc}\alpha=(\delb\sdelb+\sdelb\delb)\alpha.
\end{eqnarray*}

Now, we have type Kahler identities for the previous operators.

\begin{Le} Let $(M^{2m+2},g,\Nc,\Rc,\Rcs)$ be an almost Robinson manifold endowed with its Chern-Robinson connection $\nabla$. Then the following identities hold on $\Omega_{\Sc}^{p,q}(M;E)$ :
\begin{enumerate}
\item $[\Lambda_{\Sc},L_{\Sc}]=(m-p-q)I_{\Sc}$.
\item $[\del,L_{\Sc}]=0,\,[\delb,L_{\Sc}]=0,\,[\sdel,\Lambda_{\Sc}]=0,\,[\sdelb,\Lambda_{\Sc}]=0$
\item $[\del,\Lambda_{\Sc}]=-\sqrt{-1}\,\sdelb,\,[\delb,\Lambda_{\Sc}]=\sqrt{-1}\,\sdel,\,[\sdel,L_{\Sc}]=-\sqrt{-1}\,\delb,\,[\sdelb,L_{\Sc}]=\sqrt{-1}\,\del.$
\item $[\del\delb+\delb\del,\Lambda_{\Sc}]=\sqrt{-1}(\Box_{\Sc}-\overline{\Box}_{\Sc}),\,[\Box_{\Sc},\Lambda_{\Sc}]=-\sqrt{-1}(\sdel\sdelb+\sdelb\sdel),\,\del\sdelb+\sdelb\del=\sqrt{-1}[\del^2,\Lambda_{\Sc}]$.\\
\end{enumerate}
\end{Le}

For $A\in\Gamma(\otimes^2{(\So{1}{0})}^*)$ and $B\in\Gamma(\otimes^2{(\So{0}{1})}^*)$,
we define $A^{1,0}\in\Omega_{\Sc}^{1,0}(M;{(\So{1}{0})}^*)$ and $B^{0,1}\in\Omega_{\Sc}^{0,1}(M;{(\So{0}{1})}^*)$ respectively by:
$$A^{1,0}(\Xh)(\Yh)=A(\Xh,\Yh)\quad \mathrm{and}\quad B^{0,1}(\Xah)(\Yah)=B(\Xah,\Yah).$$

Now, we reformulate the Bianchi identities in terms of operators $\del$ and $\delb$.

\newpage

\begin{Prop} Let $(M^{2m+2},g,\Nc,\Rc,\Rcs)$ be a Fefferman-Robinson manifold endowed with its Chern-Robinson connection $\nabla$ and let $R_{\Cc}^{1,1}=R^{\Cc}$ viewed as an element of $\Omega_{\Sc}^{1,1}(M;\wedge^{1,1}{(\Sc^{\C})}^*)$. Then, we have :
\begin{eqnarray}
\del A^{1,0}_{\n{}{*}}&=&0,\quad \nab{\n{}{}}{}A^{1,0}_{\n{}{*}}=0,\quad \nab{\n{}{}}{}R_{\Cc}^{1,1}=0,\quad R^{2,0}=-\omega^{1,0}\wedge A^{1,0}_{\n{}{*}}\\
\del R_{\Cc}^{1,1}&=&L_{\Sc}({(\delb A^{1,0}_{\n{}{*}})}^{1,0})-\omega^{1,0}\wedge\delb A^{1,0}_{\n{}{*}}\\
\nab{\n{}{*}}{}R_{\Cc}^{1,1}&=&-(\delb{(\delb A^{1,0}_{\n{}{*}})}^{1,0}+\del{(\del A^{0,1}_{\n{}{*}})}^{0,1})+\omega^{0,1}\wedge {[A^{1,0}_{\n{}{*}},A^{0,1}_{\n{}{*}}]}_{1}^{1,0}-\omega^{1,0}\wedge {[A^{1,0}_{\n{}{*}},A^{0,1}_{\n{}{*}}]}_{1}^{0,1}\\
\del{(\delb A^{1,0}_{\n{}{*}})}^{1,0}&=&\omega^{1,0}\wedge\nab{\n{}{*}}{}A^{1,0}_{\n{}{*}}-\sqrt{-1}(A^{1,0}_{\n{}{*}}\wedge\omega^{0,1})\underset{2}{\circ}R_{\Cc}^{1,1}\\
\sdelb R_{\Cc}^{1,1}&=&{(\del Ric_{\Cc}^{0,1})}^{1,0}-\sqrt{-1}m{(\delb A^{1,0}_{\n{}{*}})}^{1,0}+\omega^{1,0}\wedge\sdel A^{1,0}_{\n{}{*}}\nonumber\\
\del Ric_{\Cc}^{1,0}&=&\omega^{1,0}\wedge{(\sdel A^{1,0}_{\n{}{*}})}^{1,0}\nonumber\\
\nab{\n{}{*}}{}Ric^{\Cc}&=&\delb{(\sdel A^{1,0}_{\n{}{*}})}^{1,0}-\del{(\sdelb A^{0,1}_{\n{}{*}})}^{0,1}\\
\Box_{\Sc}A^{1,0}_{\n{}{*}}&=&{[A^{1,0}_{\n{}{*}},Ric_{\Cc}^{0,1}]}_{1}^{1,0}\\
\overline{\Box}_{\Sc}A^{1,0}_{\n{}{*}}&=&\sqrt{-1}(m-1)\nab{\n{}{*}}{}A^{1,0}_{\n{}{*}}-{(Ric_{\Cc}^{0,1}\underset{1}{\circ}A^{1,0}_{\n{}{*}})}^{1,0}-\bigl(\overset{\circ}{R_{\Cc}^{1,1}}(A_{\n{}{*}}){\bigr)}^{1,0}\nonumber\\
\sdelb Ric_{\Cc}^{0,1}&=&\partial_{\Sc} s^{\Cc}+\sqrt{-1}(m-1){(\sdel A^{1,0}_{\n{}{*}})}^{1,0}\\
ds^{\Cc}(\n{}{*})&=&\sdel{(\sdel A^{1,0}_{\n{}{*}})}^{1,0}+\sdelb{(\sdelb A^{0,1}_{\n{}{*}})}^{0,1},
\end{eqnarray}
with
$${(\sdel A^{1,0}_{\n{}{*}})}^{1,0}(\Xh)=(\sdel A^{1,0}_{\n{}{*}})(\Xh),\quad{(\sdelb A^{0,1}_{\n{}{*}})}^{0,1}(\Xah)=(\sdelb A^{0,1}_{\n{}{*}})(\Xah)$$
and
$$\overset{\circ}{R_{\Cc}^{1,1}}(A_{\n{}{*}})(\,.\,,\,.\,)=\sum_{i,j=1}^{m}R_{\Cc}^{1,1}(\,.\,,\overline{Z_i})(\,.\,,\overline{Z_j})A_{\n{}{*}}(Z_i,Z_j).$$
\end{Prop}

\noindent{Proof.}
For $\gamma\in\Omega_{\Sc}^{1,1}(M;{(\So{1}{0})}^*)$ and $\mu\in\Omega_{\Sc}^{1,1}(M;{(\So{0}{1})}^*)$, then $\gamma^{1,0}=\overline{t}(\gamma)\in\Omega_{\Sc}^{1,0}(M;\wedge^{1,1}{(\Sc^{\C})}^*)$ and $\mu^{0,1}=t(\mu)\in\Omega_{\Sc}^{0,1}(M;\wedge^{1,1}{(\Sc^{\C})}^*)$ are defined using (\ref{e39}) by :
\begin{eqnarray*}
\gamma^{1,0}(\Xh)(\Zh,\Wah)&=&\gamma(\Xh,\Zh)(\Wah)-\gamma(\Xh,\Wah)(\Zh)=-\gamma(\Xh,\Wah)(\Zh)\\
\mu^{0,1}(\Yah)(\Zh,\Wah)&=&\mu(\Yah,\Zh)(\Wah)-\mu(\Yah,\Wah)(\Zh)=-\mu(\Zh,\Yah)(\Wah).
\end{eqnarray*}
We deduce, using the formulas for $\del$ and $\delb$ that :
\begin{eqnarray*}
(\del\gamma^{1,0})(\Xh,\Yh)(\Zh,\Wah)&=&-(\del\gamma)(\Xh,\Yh,\Wah)(\Zh)\\
(\delb\gamma^{1,0})(\Xh,\Yah)(\Zh,\Wah)&=&(\nab{\Yah}{}\gamma)(\Xh,\Wah)(\Zh)\\
(\del\mu^{0,1})(\Xh,\Yah)(\Zh,\Wah)&=&-(\nab{\Xh}{}\mu)(\Zh,\Yah)(\Wah)\\
(\delb\mu^{0,1})(\Xah,\Yah)(\Zh,\Wah)&=&(\delb\mu)(\Xah,\Yah,\Zh)(\Wah).
\end{eqnarray*}
Propositions 3.2 and 5.1 together with the previous formulas yields formulas (41),(42),(43) and (44). Now, we can verifie that
\begin{equation}\label{e49}
\Lambda_{\Sc}(\omega^{1,0}\wedge\gamma)=-\gamma^{1,0}+\omega^{1,0}\wedge(\Lambda_{\Sc}\gamma)\quad\mathrm{and}\quad\Lambda_{\Sc}(\omega^{0,1}\wedge\mu)=-\mu^{0,1}+\omega^{0,1}\wedge(\Lambda_{\Sc}\mu).
\end{equation}
Using the relations $[\del,\Lambda_{\Sc}]=-\sqrt{-1}\,\sdelb$ and $[\delb,\Lambda_{\Sc}]=\sqrt{-1}\,\sdel$, we obtain by (42) and (\ref{e49}) :
\begin{eqnarray*}
\sdelb R_{\Cc}^{1,1}&=&\sqrt{-1}(\del\Lambda_{\Sc} R_{\Cc}^{1,1}-\Lambda_{\Sc}\del R_{\Cc}^{1,1})\\
&=&\sqrt{-1}(\del\Lambda_{\Sc} R_{\Cc}^{1,1})-\sqrt{-1}\Lambda_{\Sc}L_{\Sc}({(\delb A^{1,0}_{\n{}{*}})}^{1,0})+\sqrt{-1}\Lambda_{\Sc}(\omega^{1,0}\wedge\delb A^{1,0}_{\n{}{*}})\\
&=&{(\del Ric_{\Cc}^{0,1})}^{1,0}-\sqrt{-1}(m-1){(\delb A^{1,0}_{\n{}{*}})}^{1,0}-\sqrt{-1}{(\delb A^{1,0}_{\n{}{*}})}^{1,0}+\sqrt{-1}\omega^{1,0}\wedge\Lambda_{\Sc}\delb A^{1,0}_{\n{}{*}}\\
&=&{(\del Ric_{\Cc}^{0,1})}^{1,0}-\sqrt{-1}m{(\delb A^{1,0}_{\n{}{*}})}^{1,0}+\omega^{1,0}\wedge\sdel A^{1,0}_{\n{}{*}}.
\end{eqnarray*}
For $\phi\in\Omega_{\Sc}^{p,q}(M;\wedge^{1,1}{(\Sc^{\C})}^*)$, then
$\tr\,\phi\in\Omega_{\Sc}^{p,q}(M)$ is given by :
$$(\tr\,\phi)(X_1^{1,0},\ldots,X_p^{1,0},Y_1^{0,1},\ldots,Y_q^{0,1})=\sum_{i}\phi(X_1^{1,0},\ldots,X_p^{1,0},Y_1^{0,1},\ldots,Y_q^{0,1})(Z_i,\overline{Z_i}).$$
Now, we have
$$
(\del Ric_{\Cc}^{1,0})(\Xh,\Yh)(\Zah)=\tr\,(\del R_{\Cc}^{1,1})(\Xh,\Yh,\Zah)=\sum_{i}(\del R_{\Cc}^{1,1})(\Xh,\Yh,\Zah)(Z_i,\overline{Z_i}).
$$
Since
\begin{eqnarray*}
(\del R_{\Cc}^{1,1})(\Xh,\Yh,\Zah)(\Uh,\Wah)&=&(\del R_{\Cc}^{1,1})(\Xh,\Uh,\Wah)(\Yh,\Zah)\\
&-&(\del R_{\Cc}^{1,1})(\Yh,\Uh,\Wah)(\Xh,\Zah),
\end{eqnarray*}
we deduce that
\begin{eqnarray*}
(\del Ric_{\Cc}^{1,0})(\Xh,\Yh)(\Zah)&=&\sqrt{-1}\bigl((\Lambda_{\Sc}\del R_{\Cc}^{1,1})(\Xh)(\Yh,\Zah)-(\Lambda_{\Sc}\del R_{\Cc}^{1,1})(\Yh)(\Xh,\Zah)\bigr)\\
&=&\sqrt{-1}m\bigl({(\delb A^{1,0}_{\n{}{*}})}^{1,0}(\Xh)(\Yh,\Zah)-{(\delb A^{1,0}_{\n{}{*}})}^{1,0}
(\Yh)(\Xh,\Zah)\bigr)\\
&+&(\omega^{1,0}\wedge\sdel A^{1,0}_{\n{}{*}})(\Yh)(\Xh,\Zah)-(\omega^{1,0}\wedge\sdel A^{1,0}_{\n{}{*}})(\Xh)(\Yh,\Zah).
\end{eqnarray*}

Using that ${(\delb A^{1,0}_{\n{}{*}})}^{1,0}(\Xh)(\Yh,\Zah)={(\delb A^{1,0}_{\n{}{*}})}^{1,0}(\Yh)(\Xh,\Zah)$, we deduce that
$$(\del Ric_{\Cc}^{1,0})(\Xh,\Yh)(\Zah)=(\omega^{1,0}\wedge{(\sdel A^{1,0}_{\n{}{*}})}^{1,0})(\Xh,\Yh)(\Zah).$$
Now
\begin{eqnarray*}
\nab{\n{}{*}}{}Ric^{\Cc}&=&\tr\,(\nab{\n{}{*}}{}R_{\Cc}^{1,1})(\Xh,\Yah)=\sum_{i}(\nab{\n{}{*}}{}R_{\Cc}^{1,1})(\Xh,\Yah)(Z_i,\overline{Z_i})\\
&=&-\sum_{i}\bigl(\delb{(\delb A^{1,0}_{\n{}{*}})}^{1,0}+\del{(\del A^{0,1}_{\n{}{*}})}^{0,1}\bigr)(\Xh,\Yah)(Z_i,\overline{Z_i})\\
&+&\sum_{i}\bigl(\omega^{0,1}\wedge {[A^{1,0}_{\n{}{*}},A^{0,1}_{\n{}{*}}]}_{1}^{1,0}-\omega^{1,0}\wedge {[A^{1,0}_{\n{}{*}},A^{0,1}_{\n{}{*}}]}_{1}^{0,1})\bigr)(\Xh,\Yah)(Z_i,\overline{Z_i})\\
&=&(\delb{(\sdel A^{1,0}_{\n{}{*}})}^{1,0}-\del{(\sdelb A^{0,1}_{\n{}{*}})}^{0,1})(\Xh,\Yah).
\end{eqnarray*}
Using the relation $[\del\delb+\delb\del,\Lambda_{\Sc}]=\sqrt{-1}(\Box_{\Sc}-\overline{\Box}_{\Sc})$, we obtain
\begin{eqnarray*}
(\overline{\Box}_{\Sc}-\Box_{\Sc})A^{1,0}_{\n{}{*}}&=&-\sqrt{-1}\Lambda_{\Sc}(\del\delb+\delb\del)A^{1,0}_{\n{}{*}}=-\sqrt{-1}\Lambda_{\Sc}\del\delb A^{1,0}_{\n{}{*}}=c\bigl(\del{(\delb A^{1,0}_{\n{}{*}})}^{1,0}\bigr)\\
&=&c\bigl(\omega^{1,0}\wedge\nab{\n{}{*}}{}A^{1,0}_{\n{}{*}}\bigr)-\sqrt{-1}c\bigl((A^{1,0}_{\n{}{*}}\wedge\omega^{0,1})\underset{2}{\circ}R_{\Cc}^{1,1}
\bigr).
\end{eqnarray*}
Now, for $A\in\Gamma(\odot^2{(\So{1}{0})}^*)$ and $\phi\in\Omega_{\Sc}^{1,1}(M;\wedge^{1,1}{(\Sc^{\C})}^*)$, we have
\begin{eqnarray}\label{e50}
c(\omega^{1,0}\wedge A^{1,0})(\Xh,\Yh)&=&\sqrt{-1}(m-1) A(\Xh,\Yh)\nonumber\\
((A^{1,0}\wedge\omega^{0,1})\underset{2}{\circ}\phi)(\Xh,\Yh)(\Uh,\Wah)&=&-\sqrt{-1}\sum_{i}\bigl(\phi(\Xh,\overline{Z_i})(\Uh,\Wah)A(\Yh,Z_i)\nonumber\\
&&-\phi(\Yh,\overline{Z_i})(\Uh,\Wah)A(\Xh,Z_i)\bigr).
\end{eqnarray}
We deduce from (\ref{e50}) that
$$
(\overline{\Box}_{\Sc}-\Box_{\Sc})A^{1,0}_{\n{}{*}}
=\sqrt{-1}(m-1)\nab{\n{}{*}}{}A^{1,0}_{\n{}{*}}-{(A^{1,0}_{\n{}{*}}\underset{1}{\circ} Ric_{\Cc}^{0,1})}^{1,0}-{\bigl(\overset{\circ}{R_{\Cc}^{1,1}}(A_{\n{}{*}})\bigr)}^{1,0}.
$$
Now,
\begin{eqnarray*}
\del{(\delb A^{1,0}_{\n{}{*}})}^{1,0}(\Xh,\Yh)(\Uh,\Wah)&=&\del{(\delb A^{1,0}_{\n{}{*}})}^{1,0}(\Uh,\Yh)(\Xh,\Wah)\\&-&\del{(\delb A^{1,0}_{\n{}{*}})}^{1,0}(\Uh,\Xh)(\Yh,\Wah).
\end{eqnarray*}
We deduce that
\begin{eqnarray*}
(\Box_{\Sc}A^{1,0}_{\n{}{*}})(\Xh)(\Yh)&=&(\del{(\sdel A^{1,0}_{\n{}{*}})}^{1,0})(\Xh,\Yh)=\tr\,(\del{(\delb A^{1,0}_{\n{}{*}})}^{1,0})(\Xh,\Yh)\\
&=&\sum_{i}\del{(\delb A^{1,0}_{\n{}{*}})}^{1,0}(\Xh,\Yh)(Z_i,\overline{Z_i})\\
&=&c\bigl(\del{(\delb A^{1,0}_{\n{}{*}})}^{1,0}\bigr)(\Xh,\Yh)-c\bigl(\del{(\delb A^{1,0}_{\n{}{*}})}^{1,0}\bigr)(\Yh,\Xh)\\
&=&{[A^{1,0}_{\n{}{*}},Ric_{\Cc}^{0,1}]}_{1}(\Xh,\Yh)={[A^{1,0}_{\n{}{*}},Ric_{\Cc}^{0,1}]}_{1}^{1,0}(\Xh)(\Yh).
\end{eqnarray*}
Since $\Box_{\Sc}A^{1,0}_{\n{}{*}}={[A^{1,0}_{\n{}{*}},Ric_{\Cc}^{0,1}]}_{1}^{1,0}$, we obtain
\begin{eqnarray*}
\overline{\Box}_{\Sc}A^{1,0}_{\n{}{*}}&=&{[A^{1,0}_{\n{}{*}},Ric_{\Cc}^{0,1}]}_{1}^{1,0}+\sqrt{-1}(m-1)\nab{\n{}{*}}{}A^{1,0}_{\n{}{*}}-{(A^{1,0}_{\n{}{*}}\underset{1}{\circ} Ric_{\Cc}^{0,1})}^{1,0}-\overset{\circ}{R_{\Cc}^{1,1}}(A_{\n{}{*}})\\
&=&\sqrt{-1}(m-1)\nab{\n{}{*}}{}A^{1,0}_{\n{}{*}}-{(Ric_{\Cc}^{0,1}\underset{1}{\circ}A^{1,0}_{\n{}{*}})}^{1,0}-\overset{\circ}{R_{\Cc}^{1,1}}(A_{\n{}{*}}).
\end{eqnarray*}
We have
\begin{eqnarray*}
(\sdelb Ric_{\Cc}^{0,1})(\Xh)&=&\tr\,(\sdelb R_{\Cc}^{1,1})(\Xh)\\
&=&\tr\,{(\del Ric_{\Cc}^{0,1})}^{1,0}(\Xh)-\sqrt{-1}m\tr\,{(\delb A^{1,0}_{\n{}{*}})}^{1,0}(\Xh)+\tr\,(\omega^{1,0}\wedge\sdel A^{1,0}_{\n{}{*}})(\Xh)\\
&=&\sum_{i}({(\del Ric_{\Cc}^{0,1})}^{1,0}-\sqrt{-1}m{(\delb A^{1,0}_{\n{}{*}})}^{1,0}+\omega^{1,0}\wedge\sdel A^{1,0}_{\n{}{*}})(\Xh)(Z_i,\overline{Z_i})\\
&=&(\partial_{\Sc} s^{\Cc}+\sqrt{-1}(m-1){(\sdel A^{1,0}_{\n{}{*}})}^{1,0})(\Xh).
\end{eqnarray*}
To conclude, we have by (45)
$$ds^{\Cc}(\n{}{*})=\sqrt{-1}\Lambda_{\Sc}(\nab{\n{}{*}}{}Ric^{\Cc})
=\sqrt{-1}\Lambda_{\Sc}(\delb{(\sdel A^{1,0}_{\n{}{*}})}^{1,0}-\del{(\sdelb A^{0,1}_{\n{}{*}})}^{0,1})=\sdel{(\sdel A^{1,0}_{\n{}{*}})}^{1,0}+\sdelb{(\sdelb A^{0,1}_{\n{}{*}})}^{0,1}.\quad\Box$$
\newpage
Let $\rho\in\Omega_{\Sc}^{1}(M)$ given by
$$\rho=\frac{1}{2(m+1)}d_{\Sc} s^{\Cc}+\sqrt{-1}({(\sdel A^{1,0}_{\n{}{*}})}^{1,0}-{(\sdelb A^{0,1}_{\n{}{*}})}^{0,1}).$$
We define $C_{\Yc}\in\Omega_{\Sc}^{1,1}(M;{(\So{1}{0})}^*)$, $\overline{C}_{\Yc}\in\Omega_{\Sc}^{1,1}(M;{(\So{0}{1})}^*)$, $Q_{\Cc}\in\Gamma(\bigodot^2{(\So{1}{0})}^*\oplus\bigodot^2{(\So{0}{1})}^*)$ and $B_{\Cc}\in\Omega_{\Sc}^{1,1}(M)$ respectively by :
\begin{eqnarray*}
C_{\Yc}&=&\delb A^{1,0}_{\n{}{*}}+\sqrt{-1}\del P_{\Cc}^{0,1}+\frac{1}{m+2}\omega^{0,1}\wedge\rho^{1,0},\quad C_{\Yc_{0}}=C_{\Yc}+\frac{2}{m+2}\omega\otimes\rho^{1,0}\\
\overline{C}_{\Yc}&=&\del A^{0,1}_{\n{}{*}}+\sqrt{-1}\,\delb P_{\Cc}^{1,0}+\frac{1}{m+2}\omega^{1,0}\wedge\rho^{0,1},\quad\overline{C}_{\Yc_{0}}=\overline{C}_{\Yc}+\frac{2}{m+2}\omega\otimes\rho^{0,1}\\
Q_{\Cc}&=&-\nab{\n{}{*}}{}A^{2,0}_{\n{}{*}}-\nab{\n{}{*}}{}A^{0,2}_{\n{}{*}}-\sqrt{-1}\bigl({\{A^{1,0}_{\n{}{*}},P_{\Cc}^{0,1}\}}_1+{\{A^{0,1}_{\n{}{*}},P_{\Cc}^{1,0}\}}_1\bigr)+\frac{1}{m+2}(\nab{(1,0)}{sym}\rho^{1,0}+\nab{(0,1)}{sym}\rho^{0,1})\\
B_{\Cc}&=&\frac{1}{m+2}(\del\rho^{0,1}-\delb\rho^{1,0})-{[A^{1,0}_{\n{}{*}},A^{0,1}_{\n{}{*}}]}_1+{[P_{\Cc}^{1,0},P_{\Cc}^{0,1}]}_1.
\end{eqnarray*}

We call $C_{\Yc}$ (and also $\overline{C}_{\Yc}$) the Cotton-Robinson tensor, $Q_{\Cc}$ the Cartan-Robinson tensor and $B_{\Cc}$ the Bach-Robinson tensor.

\begin{Prop} Let $(M^{2m+2},g,\Nc,\Rc,\Rcs)$ be a Fefferman-Robinson manifold. Then, we have :
\begin{eqnarray}
\del C_{\Mc}^{1,1}&=&L_{\Sc}(C_{\Yc}^{1,0})-\omega^{1,0}\wedge C_{\Yc}=L_{\Sc}(C_{\Yc_{0}}^{1,0})-\omega^{1,0}\wedge C_{\Yc_{0}}\nonumber\\
\delb C_{\Mc}^{1,1}&=&L_{\Sc}(\overline{C}_{\Yc}^{0,1})-\omega^{0,1}\wedge\overline{C}_{\Yc}=L_{\Sc}(\overline{C}_{\Yc_{0}}^{0,1})-\omega^{0,1}\wedge \overline{C}_{\Yc_{0}}\nonumber\\
\del C_{\Yc_{0}}^{1,0}&=&-\bigl(\omega^{1,0}\wedge Q_{\Cc}^{1,0}+\sqrt{-1}(A^{1,0}_{\n{}{*}}\wedge\omega^{0,1})\underset{2}{\circ} C_{\Mc}^{1,1}\bigr)\\
\delb C_{\Yc_{0}}^{1,0}+\del\overline{C}_{\Yc_{0}}^{0,1}&=&-\bigl(B_{\Cc}^{1,0}\wedge\omega^{0,1}-\omega^{1,0}\wedge B_{\Cc}^{0,1}+\nab{\n{}{*}}{}C_{\Mc}^{1,1}+C_{\Mc}^{1,1}\underset{2}{\circ}(\omega^{1,0}\wedge P_{\Cc}^{0,1}-P_{\Cc}^{1,0}\wedge\omega^{0,1})\bigr)\\
\sdelb C_{\Mc}^{1,1}&=&-\sqrt{-1}mC_{\Yc_{0}}^{1,0},\quad\sdel C_{\Mc}^{1,1}=\sqrt{-1}m\overline{C}_{\Yc_{0}}^{0,1}\\
\sdelb C_{\Yc_{0}}&=&-\Bigl(\sqrt{-1}(m-1)Q_{\Cc}^{1,0}+\bigl(\overset{\circ}{C_{\Mc}^{1,1}}(A_{\n{}{*}}){\bigr)}^{1,0}\Bigr)\nonumber\\
\sdelb \overline{C}_{\Yc_{0}}&=&\sdel C_{\Yc_{0}}^{1,0}=-\sqrt{-1}\Bigl(mB_{\Cc_{0}}^{1,0}+{\bigl(C_{\Mc}^{1,1}\underset{2}{\circ} P^{\Cc}\bigr)}^{1,0}\Bigr),
\end{eqnarray}
with $B_{\Cc_{0}}=B_{\Cc}-\dfrac{\sqrt{-1}}{m}\omega\otimes\Bigl(\dfrac{1}{m+2}(\sdelb\rho^{0,1}+\sdel\rho^{1,0})+\tr\,{[A^{1,0}_{\n{}{*}},A^{0,1}_{\n{}{*}}]}_1-\tr\,{[P_{\Cc}^{1,0},P_{\Cc}^{0,1}]}_1\Bigr)$.
\end{Prop}

The proof of this proposition used the following lemmas.

\begin{Le} For $\rho^{1,0}\in\Omega_{\Sc}^{1,0}(M)$, $\rho^{0,1}\in\Omega_{\Sc}^{0,1}(M)$ and $\Omega\in\Omega_{\Sc}^{1,1}(M)$, we have the relations :
\begin{eqnarray*}
\del{(\omega\odot\Omega)}^{1,1}&=&\omega\wedge{(\del\Omega^{0,1})}^{1,0}+\del\Omega\otimes\omega\\
\del{(\omega^{1,0}\wedge\rho^{0,1})}^{0,1}&=&-\del\rho^{0,1}\otimes\omega\\
\del{(\omega^{0,1}\wedge\rho^{1,0})}^{1,0}&=&-\del\rho^{1,0}\otimes\omega,
\quad\delb{(\omega^{0,1}\wedge\rho^{1,0})}^{1,0}=-\delb \rho^{1,0}\otimes\omega\\
\del{(\omega\otimes\rho^{0,1})}^{0,1}&=&-\omega^{0,1}\wedge{(\del\rho^{0,1})}^{1,0},\quad \delb{(\omega\otimes\rho^{1,0})}^{1,0}=-\omega^{1,0}\wedge{(\delb\rho^{1,0})}^{0,1}\\
\del{(\omega\otimes\rho^{1,0})}^{1,0}&=&-\omega^{1,0}\wedge{(\nab{}{}\rho^{1,0})}^{1,0}.
\end{eqnarray*}
\end{Le}

\begin{Le} Let $A^{1,0}\in\Omega_{\Sc}^{1,0}(M;{(\So{1}{0})}^*),B^{0,1}\in\Omega_{\Sc}^{0,1}(M;{(\So{0}{1})}^*)$, $P^{1,0},R^{1,0}\in\Omega_{\Sc}^{1,0}(M;{(\So{0}{1})}^*)$, $Q^{0,1},S^{0,1}\in\Omega_{\Sc}^{0,1}(M;{(\So{1}{0})}^*)$ and $\Omega\in\Omega_{\Sc}^{1,1}(M)$, we have the relations :
\begin{eqnarray*}
(P^{1,0}\wedge Q^{0,1})\underset{2}{\circ}(R^{1,0}\wedge S^{0,1})&=&-{[R^{1,0},Q^{0,1}]}_{1}^{1,0}\wedge{[P^{1,0},S^{0,1}]}_{1}^{0,1}\\
(A^{1,0}\wedge Q^{0,1})\underset{2}{\circ}(R^{1,0}\wedge S^{0,1})&=&-{(A^{1,0}\underset{1}{\circ}S^{0,1})}^{1,0}\wedge{[R^{1,0},Q^{0,1}]}_{1}^{1,0}\\
(\omega^{1,0}\wedge A^{1,0})\underset{2}{\circ}(\omega^{0,1}\wedge B^{0,1})&=& A^{1,0}\wedge B^{0,1}+\sqrt{-1}{[A^{1,0},B^{0,1}]}_{1}^{1,0}\wedge\omega^{0,1}\\
(A^{1,0}\wedge Q^{0,1})\underset{2}{\circ}(\omega\odot\Omega)&=&-\sqrt{-1}{[A^{1,0},Q^{0,1}]}_{1}\otimes\Omega-{[A^{1,0},{[\Omega^{1,0},Q^{0,1}]}_{1}^{0,1}]}_{1}\otimes\omega\\
(\omega\odot\Omega)\underset{2}{\circ}(R^{1,0}\wedge S^{0,1})&=&-\sqrt{-1}\Omega\otimes{[R^{1,0},S^{0,1}]}_{1}-\omega\otimes{[R^{1,0},{[\Omega^{1,0},S^{0,1}]}_{1}^{0,1}]}_{1}.
\end{eqnarray*}
\end{Le}
(We omit the proofs which are respectively direct consequences of definitions of $\del$ and $\delb$ for Lemma 5.2 and definitions of products $\underset{1}{\circ}$ and $\underset{2}{\circ}$ for Lemma 5.3).\\

\noindent{Proof of Proposition 5.3}. Since $C_{\Mc}^{1,1}\in\Omega_{\Sc}^{1,1}(M;\wedge^{1,1}{(\Sc^{\C})}^*)$, then we have by (\ref{e16}) :
$$C_{\Mc}^{1,1}=R_{\Cc}^{1,1}+\sqrt{-1}\bigl(P_{\Cc}^{1,0}\wedge\omega^{0,1}+\omega^{1,0}\wedge P_{\Cc}^{0,1}+{(\omega\odot P^{\Cc})}^{1,1}\bigr).$$
Now, we have for $P^{1,0}\in\Omega_{\Sc}^{1,0}(M;{(\So{0}{1})}^*)$ and $Q^{0,1}\in\Omega_{\Sc}^{0,1}(M;{(\So{1}{0})}^*)$ the formula $\del(P^{1,0}\wedge Q^{0,1})=\del P^{1,0}\wedge Q^{0,1}-P^{1,0}\wedge\del Q^{0,1}$. Lemma 5.2 together with the previous formula, the relations $\del P_{\Cc}=\dfrac{\sqrt{-1}}{m+2}L_{\Sc}(\rho^{1,0})$, $\del P_{\Cc}^{1,0}=-\dfrac{\sqrt{-1}}{m+2}\omega^{1,0}\wedge\rho^{1,0}$ and (42) yields
\begin{eqnarray*}
\del C_{\Mc}^{1,1}&=&\del R_{\Cc}^{1,1}+\sqrt{-1}\bigl(\del P_{\Cc}^{1,0}\wedge \omega^{0,1}-\omega^{1,0}\wedge\del P_{\Cc}^{0,1}+\del P^{\Cc}\otimes\omega+\omega\wedge{(\del P_{\Cc}^{0,1})}^{1,0}\bigr)\\
&=&L_{\Sc}({(\delb A^{1,0}_{\n{}{*}})}^{1,0})-\omega^{1,0}\wedge\delb A^{1,0}_{\n{}{*}}
-\frac{1}{m+2}\omega^{1,0}\wedge\omega^{0,1}\wedge\rho^{1,0}-\sqrt{-1}\omega^{1,0}\wedge\del P_{\Cc}^{0,1}\\
&-&\frac{1}{m+2}L_{\Sc}(\rho^{1,0})\otimes\omega+
\sqrt{-1}L_{\Sc}({(\del P_{\Cc}^{0,1})}^{1,0}).
\end{eqnarray*}
Since $L_{\Sc}\bigl({(\omega^{0,1}\wedge\rho^{1,0})}^{1,0}\bigr)=-L_{\Sc}(\rho^{1,0})\otimes\omega$, we deduce that
\begin{eqnarray*}
\del C_{\Mc}^{1,1}&=&L_{\Sc}\bigl({(\delb A^{1,0}_{\n{}{*}})}^{1,0}+\sqrt{-1}{(\del P_{\Cc}^{0,1})}^{1,0}+\frac{1}{m+2}{(\omega^{0,1}\wedge\rho^{1,0})}^{1,0}\bigr)\\
&-&\omega^{1,0}\wedge\bigl(\delb A^{1,0}_{\n{}{*}}+\sqrt{-1}\del P_{\Cc}^{0,1}+\frac{1}{m+2}\omega^{0,1}\wedge\rho^{1,0}\bigr)\\
&=&L_{\Sc}(C_{\Yc}^{1,0})-\omega^{1,0}\wedge C_{\Yc}.
\end{eqnarray*}
Since $L_{\Sc}({(\omega\otimes\rho^{1,0})}^{1,0})=\omega^{1,0}\wedge(\omega\otimes\rho^{1,0})$, we deduce that
$$\del C_{\Mc}^{1,1}=L_{\Sc}(C_{\Yc_{0}}^{1,0})-\omega^{1,0}\wedge C_{\Yc_{0}}.$$
Using the fact that $\delb R_{\Cc}^{1,1}=L_{\Sc}({(\delb A^{0,1}_{\n{}{*}})}^{0,1})-\omega^{0,1}\wedge\del A^{0,1}_{\n{}{*}}$, we obtain by the same way the formula for $\delb C_{\Mc}^{1,1}$. Now, using Lemma 5.2, we have
\begin{eqnarray*}
\del C_{\Yc_{0}}^{1,0}&=&\del C_{\Yc}^{1,0}+\frac{2}{m+2}\del{(\omega\otimes\rho^{1,0})}^{1,0}\\
&=&\del{(\delb A^{1,0}_{\n{}{*}})}^{1,0}+\sqrt{-1}\del{(\del P_{\Cc}^{0,1})}^{1,0}+\frac{1}{m+2}\del{(\omega^{0,1}\wedge\rho^{1,0})}^{1,0}+\frac{2}{m+2}\del{(\omega\otimes\rho^{1,0})}^{1,0}\\
&=&\del{(\delb A^{1,0}_{\n{}{*}})}^{1,0}+\sqrt{-1}\del{(\del P_{\Cc}^{0,1})}^{1,0}-\frac{1}{m+2}\del\rho^{1,0}\otimes\omega-\frac{2}{m+2}\omega^{1,0}\wedge{(\nab{}{}\rho^{1,0})}^{1,0}.
\end{eqnarray*}
Now, the following formula holds
\begin{eqnarray}\label{e55}
\del{(\del P_{\Cc}^{0,1})}^{1,0}&=&\sqrt{-1} A^{1,0}_{\n{}{*}}\wedge P_{\Cc}^{0,1}+\omega^{1,0}\wedge {(A^{1,0}_{\n{}{*}}\underset{1}{\circ}P_{\Cc}^{0,1})}^{1,0}.
\end{eqnarray}
Moreover, by (46) we have
$$\del\rho^{1,0}=\sqrt{-1}\del{(\sdel A^{1,0}_{\n{}{*}})}^{1,0}=\sqrt{-1}\Box_{\Sc}A^{1,0}_{\n{}{*}}=\sqrt{-1}(m+2){[A^{1,0}_{\n{}{*}},P_{\Cc}^{0,1}]}_{1}^{1,0}.$$
Using
${(\nab{}{}\rho^{1,0})}^{1,0}=\dfrac{1}{2}\nab{(1,0)}{sym}\rho^{1,0}+\dfrac{\sqrt{-1}(m+2)}{2}{[A^{1,0}_{\n{}{*}},P_{\Cc}^{0,1}]}_{1}^{1,0}$ formula (44) and (\ref{e55}), we deduce that
\begin{eqnarray*}
\del C_{\Yc_{0}}^{1,0}&=&\omega^{1,0}\wedge\nab{\n{}{*}}{}A^{1,0}_{\n{}{*}}-\sqrt{-1}(A^{1,0}_{\n{}{*}}\wedge\omega^{0,1})\underset{2}{\circ}R_{\Cc}^{1,1}-A^{1,0}_{\n{}{*}}\wedge P_{\Cc}^{0,1}+\sqrt{-1}\omega^{1,0}\wedge {(A^{1,0}_{\n{}{*}}\underset{1}{\circ}P_{\Cc}^{0,1})}^{1,0}\\
&-&\sqrt{-1}{[A^{1,0}_{\n{}{*}},P_{\Cc}^{0,1}]}_{1}^{1,0}\otimes\omega
-\sqrt{-1}\omega^{1,0}\wedge {[A^{1,0}_{\n{}{*}},P_{\Cc}^{0,1}]}_{1}^{1,0}-\frac{1}{m+2}\omega^{1,0}\wedge \nab{(1,0)}{sym}\rho^{1,0}\\
&=&\omega^{1,0}\wedge\Bigl(\nab{\n{}{*}}{}A^{1,0}_{\n{}{*}}+\sqrt{-1}{\{A^{1,0}_{\n{}{*}},P_{\Cc}^{0,1}\}}_{1}^{1,0}-\frac{1}{m+2}\nab{(1,0)}{sym}\rho^{1,0}\Bigr)\\
&-&\sqrt{-1}\Bigl((A^{1,0}_{\n{}{*}}\wedge\omega^{0,1})\underset{2}{\circ}R_{\Cc}^{1,1}-\sqrt{-1}A^{1,0}_{\n{}{*}}\wedge P_{\Cc}^{0,1}+{(A^{1,0}_{\n{}{*}}\underset{1}{\circ}P_{\Cc}^{0,1})}^{1,0}\wedge\omega^{1,0}+{[A^{1,0}_{\n{}{*}},P_{\Cc}^{0,1}]}_{1}^{1,0}\otimes\omega\Bigr).
\end{eqnarray*}
Now, we have
\begin{eqnarray*}
(A^{1,0}_{\n{}{*}}\wedge\omega^{0,1})\underset{2}{\circ} C_{\Mc}^{1,1}&=&(A^{1,0}_{\n{}{*}}\wedge\omega^{0,1})\underset{2}{\circ}R_{\Cc}^{1,1}
+\sqrt{-1}(A^{1,0}_{\n{}{*}}\wedge\omega^{0,1})\underset{2}{\circ}(P_{\Cc}^{1,0}\wedge\omega^{0,1})\\
&+&\sqrt{-1}(A^{1,0}_{\n{}{*}}\wedge\omega^{0,1})\underset{2}{\circ}(\omega^{1,0}\wedge P_{\Cc}^{0,1})
+\sqrt{-1}(A^{1,0}_{\n{}{*}}\wedge\omega^{0,1})\underset{2}{\circ}{(\omega\odot P^{\Cc})}^{1,1}.
\end{eqnarray*}
From formula (\ref{e40}) together with Lemma 5.3 and from the relation $A^{1,0}_{\n{}{*}}\underset{1}{\circ}\omega^{0,1}=-\sqrt{-1}A^{2,0}_{\n{}{*}}$ follows :
$$
(A^{1,0}_{\n{}{*}}\wedge\omega^{0,1})\underset{2}{\circ} C_{\Mc}^{1,1}
=(A^{1,0}_{\n{}{*}}\wedge\omega^{0,1})\underset{2}{\circ}R_{\Cc}^{1,1}-\sqrt{-1}A^{1,0}_{\n{}{*}}\wedge P_{\Cc}^{0,1}+{(A^{1,0}_{\n{}{*}}\underset{1}{\circ}P_{\Cc}^{0,1})}^{1,0}\wedge\omega^{1,0}+{[A^{1,0}_{\n{}{*}},P_{\Cc}^{0,1}]}_{1}^{1,0}\otimes\omega.
$$
Hence
$$\del C_{\Yc_{0}}^{1,0}=-\bigl(\omega^{1,0}\wedge Q_{\Cc}^{1,0}+\sqrt{-1}(A^{1,0}_{\n{}{*}}\wedge\omega^{0,1})\underset{2}{\circ} C_{\Mc}^{1,1}\bigr).$$
Using Lemma 5.2 we have
\begin{eqnarray*}
\delb C_{\Yc_{0}}^{1,0}+\del\overline{C}_{\Yc_{0}}^{0,1}&=&\del{(\del A^{0,1}_{\n{}{*}})}^{0,1}+\delb{(\delb A^{1,0}_{\n{}{*}})}^{1,0}+\sqrt{-1}\bigl(\del{(\delb P_{\Cc}^{1,0})}^{0,1}+\delb{(\del P_{\Cc}^{0,1})}^{1,0}\bigr)\\
&-&\frac{1}{m+2}(\del\rho^{0,1}+\delb\rho^{1,0})\otimes\omega-\frac{2}{m+2}\bigl(\omega^{1,0}\wedge{(\delb\rho^{1,0})}^{0,1}+\omega^{0,1}\wedge{(\del\rho^{0,1})}^{1,0}\bigr).
\end{eqnarray*}
Now we have the formula :
\begin{eqnarray}\label{e56}
\del{(\delb P_{\Cc}^{1,0})}^{0,1}+\delb{(\del P_{\Cc}^{0,1})}^{1,0}&=&\sqrt{-1}R_{\Cc}^{1,1}\underset{2}{\circ}(\omega^{1,0}\wedge P_{\Cc}^{0,1}-P_{\Cc}^{1,0}\wedge\omega^{0,1})
-\omega\otimes\nab{\n{}{*}}{}P^{\Cc}.
\end{eqnarray}
By (43) and (48) we have the relation $\nab{\n{}{*}}{}P^{\Cc}=-\dfrac{\sqrt{-1}}{m+2}(\del\rho^{0,1}+\delb\rho^{1,0})$. Now (43) and (\ref{e56}) yields
\begin{eqnarray*}
\delb C_{\Yc_{0}}^{1,0}+\del\overline{C}_{\Yc_{0}}^{0,1}&=&-\nab{\n{}{*}}{}R_{\Cc}^{1,1}+\omega^{0,1}\wedge {[A^{1,0}_{\n{}{*}},A^{0,1}_{\n{}{*}}]}_{1}^{1,0}-\omega^{1,0}\wedge {[A^{1,0}_{\n{}{*}},A^{0,1}_{\n{}{*}}]}_{1}^{0,1}\\
&-&R_{\Cc}^{1,1}\underset{2}{\circ}(\omega^{1,0}\wedge P_{\Cc}^{0,1}-P_{\Cc}^{1,0}\wedge\omega^{0,1})
-\sqrt{-1}\omega\otimes\nab{\n{}{*}}{}P^{\Cc}-\frac{1}{m+2}(\del\rho^{0,1}+\delb\rho^{1,0})\otimes\omega\\
&+&\frac{1}{m+2}\Bigl(\omega^{1,0}\wedge{(\del\rho^{0,1}-\delb\rho^{1,0})}^{0,1}-\omega^{0,1}\wedge{(\del\rho^{0,1}-\delb\rho^{1,0})}^{1,0}\Bigr)\\
&-&\frac{1}{m+2}\Bigl(\omega^{1,0}\wedge{(\del\rho^{0,1}+\delb\rho^{1,0})}^{0,1}+\omega^{0,1}\wedge{(\del\rho^{0,1}+\delb\rho^{1,0})}^{1,0}\Bigr)\\
&=&\omega^{1,0}\wedge\Bigl(\frac{1}{m+2}({(\del\rho^{0,1}-\delb\rho^{1,0})}^{0,1}-{[A^{1,0}_{\n{}{*}},A^{0,1}_{\n{}{*}}]}^{0,1}_1\Bigr)\\
&-&\omega^{0,1}\wedge\Bigl(\frac{1}{m+2}({(\del\rho^{0,1}-\delb\rho^{1,0})}^{1,0}-{[A^{1,0}_{\n{}{*}},A^{0,1}_{\n{}{*}}]}^{1,0}_1\Bigr)\\
&-&\Bigl(\nab{\n{}{*}}{}R_{\Cc}^{1,1}
+\sqrt{-1}\omega\otimes\nab{\n{}{*}}{}P^{\Cc}+\sqrt{-1}\nab{\n{}{*}}{}P^{\Cc}\otimes\omega\\
&&+\sqrt{-1}\omega^{1,0}\wedge\nab{\n{}{*}}{}P_{\Cc}^{0,1}+\sqrt{-1}\omega^{0,1}\wedge\nab{\n{}{*}}{}P_{\Cc}^{1,0})+R_{\Cc}^{1,1}\underset{2}{\circ}(\omega^{1,0}\wedge P_{\Cc}^{0,1}-P_{\Cc}^{1,0}\wedge\omega^{0,1})\Bigr).
\end{eqnarray*}
Now
$$
\nab{\n{}{*}}{}C_{\Mc}^{1,1}=\nab{\n{}{*}}{}R_{\Cc}^{1,1}
+\sqrt{-1}\omega^{1,0}\wedge\nab{\n{}{*}}{}P_{\Cc}^{0,1}+\sqrt{-1}\omega^{0,1}\wedge\nab{\n{}{*}}{}P_{\Cc}^{1,0}+\sqrt{-1}\omega\odot\nab{\n{}{*}}{}P^{\Cc}
$$
and using formula (\ref{e40}) together with Lemma 5.3
\begin{eqnarray*}
C_{\Mc}^{1,1}\underset{2}{\circ}(\omega^{1,0}\wedge P_{\Cc}^{0,1}-P_{\Cc}^{1,0}\wedge\omega^{0,1})\bigr)&=&R_{\Cc}^{1,1}\underset{2}{\circ}(\omega^{1,0}\wedge P_{\Cc}^{0,1}-P_{\Cc}^{1,0}\wedge\omega^{0,1})\\
&+&\sqrt{-1}(P_{\Cc}^{1,0}\wedge\omega^{0,1})\underset{2}{\circ}(\omega^{1,0}\wedge P_{\Cc}^{0,1}-P_{\Cc}^{1,0}\wedge\omega^{0,1})\\
&+&\sqrt{-1}(\omega^{1,0}\wedge P_{\Cc}^{0,1})\underset{2}{\circ}(\omega^{1,0}\wedge P_{\Cc}^{0,1}-P_{\Cc}^{1,0}\wedge\omega^{0,1})\\
&+&\sqrt{-1}{(\omega\odot P^{\Cc})}^{1,1}\underset{2}{\circ}(\omega^{1,0}\wedge P_{\Cc}^{0,1}-P_{\Cc}^{1,0}\wedge\omega^{0,1})\\
&=&R_{\Cc}^{1,1}\underset{2}{\circ}(\omega^{1,0}\wedge P_{\Cc}^{0,1}-P_{\Cc}^{1,0}\wedge\omega^{0,1})\\
&+&\omega^{1,0}\wedge{[P_{\Cc}^{1,0},P_{\Cc}^{0,1}]}^{0,1}_1-{[P_{\Cc}^{1,0},P_{\Cc}^{0,1}]}^{1,0}_1\wedge\omega^{0,1}.
\end{eqnarray*}
Finally, we obtain that
$$
\delb C_{\Yc_{0}}^{1,0}+\del\overline{C}_{\Yc_{0}}^{0,1}=
\omega^{1,0}\wedge B_{\Cc}^{0,1}-\omega^{0,1}\wedge B_{\Cc}^{1,0}-\nab{\n{}{*}}{}C_{\Mc}^{1,1}
-C_{\Mc}^{1,1}\underset{2}{\circ}(\omega^{1,0}\wedge P_{\Cc}^{0,1}-P_{\Cc}^{1,0}\wedge\omega^{0,1}).
$$
Now, since $\Lambda_{\Sc}C_{\Mc}^{1,1}=0$, $\Lambda_{\Sc}C_{\Yc_{0}}=0$ and $\Lambda_{\Sc}\overline{C}_{\Yc_{0}}=0$, we have
\begin{eqnarray*}
\sdelb C_{\Mc}^{1,1}&=&-\sqrt{-1}\Lambda_{\Sc}\del C_{\Mc}^{1,1}=-\sqrt{-1}\Lambda_{\Sc}L_{\Sc}(C_{\Yc_{0}}^{1,0})+\sqrt{-1}\Lambda_{\Sc}(\omega^{1,0}\wedge C_{\Yc_{0}})\\
&=&-\sqrt{-1}(m-1)C_{\Yc_{0}}^{1,0}-\sqrt{-1}C_{\Yc_{0}}^{1,0}=-\sqrt{-1}mC_{\Yc_{0}}^{1,0}
\end{eqnarray*}
and
\begin{eqnarray*}
\sdel C_{\Mc}^{1,1}&=&\sqrt{-1}\Lambda_{\Sc}\delb C_{\Mc}^{1,1}=\sqrt{-1}\Lambda_{\Sc}L_{\Sc}(\overline{C}_{\Yc_{0}}^{0,1})-\sqrt{-1}\Lambda_{\Sc}(\omega^{0,1}\wedge\overline{C}_{\Yc_{0}})\\
&=&\sqrt{-1}(m-1)\overline{C}_{\Yc_{0}}^{0,1}+\sqrt{-1}\overline{C}_{\Yc_{0}}^{0,1}=\sqrt{-1}m\overline{C}_{\Yc_{0}}^{0,1}.
\end{eqnarray*}
Using $\Lambda_{\Sc}C_{\Yc_{0}}=0$, we obtain
\begin{eqnarray*}
(\sdelb C_{\Yc_{0}})(\Xh)(\Yh)&=&\sum_{i}(\nab{Z_i}{}C_{\Yc_{0}})(\Xh,\overline{Z_i})(\Yh)\\
&=&\sum_{i}\bigl((\nab{Z_i}{}C_{\Yc_{0}})(\Xh,\overline{Z_i})(\Yh)-(\nab{\Xh}{}C_{\Yc_{0}})(Z_i,\overline{Z_i})(\Yh)\bigr)\\
&=&\sum_{i}(\del C_{\Yc_{0}}^{1,0})(\Xh,Z_i)(\Yh,\overline{Z_i})=c\bigl(\del C_{\Yc_{0}}^{1,0}\bigr)(\Xh,\Yh).
\end{eqnarray*}
It follows from (\ref{e50}) that
$$(\sdelb C_{\Yc_{0}})(\Xh)(\Yh)=\Bigl(\sqrt{-1}(m-1)Q_{\Cc}^{1,0}+\bigl(\overset{\circ}{C_{\Mc}^{1,1}}(A_{\n{}{*}}){\bigr)}^{1,0}\Bigr)(\Xh)(\Yh).$$
Now, we can verifie that $\Lambda_{\Sc}(\delb C_{\Yc_{0}}^{1,0}+\del\overline{C}_{\Yc_{0}}^{0,1})=0$. It follows that $\sdel C_{\Yc_{0}}^{1,0}=\sdelb\overline{C}_{\Yc_{0}}^{0,1}$.\\
Using this relation together with $C_{\Yc_{0}}(\Xh,\Wah)(\Yh)=C_{\Yc_{0}}(\Yh,\Wah)(\Xh)$, we obtain that
\begin{eqnarray*}
(\sdelb \overline{C}_{\Yc_{0}})(\Xh)(\Wah)&=&-(\sdel C_{\Yc_{0}})(\Wah)(\Xh)=
\sum_{i}(\nab{\overline{Z_i}}{}C_{\Yc_{0}})(Z_i,\Wah)(\Xh)\\
&=&\sum_{i}(\nab{\overline{Z_i}}{}C_{\Yc_{0}})(\Xh,\Wah)(Z_i)\\
&=&\sum_{i}\bigl((\nab{\overline{Z_i}}{}C_{\Yc_{0}})(\Xh,\Wah)(Z_i)-(\nab{\Xh}{}\overline{C}_{\Yc_{0}})(Z_i,\overline{Z_i})(\Wah)\bigr)\\
&=&\sum_{i}(\delb C_{\Yc_{0}}^{1,0}+\del\overline{C}_{\Yc_{0}}^{0,1})(\Xh,\overline{Z_i})(Z_i,\Wah)\\
&=&-c\bigl(\delb C_{\Yc_{0}}^{1,0}+\del\overline{C}_{\Yc_{0}}^{0,1}\bigr)(\Xh,\Wah).
\end{eqnarray*}
Now, for $P\in\Omega_{\Sc}^{1,1}(M)$ and $\phi\in\Omega_{\Sc}^{1,1}(M;\wedge^{1,1}{(\Sc^{\C})}^*)$, we have
\begin{eqnarray}\label{e57}
c(\omega^{1,0}\wedge P^{0,1})(\Xh,\Yah)&=&-(\tr\,P)\omega(\Xh,\Yah)\nonumber\\
c(P^{1,0}\wedge\omega^{0,1})(\Xh,\Yah)&=&-\sqrt{-1}m P(\Xh,\Yh)\nonumber\\
((\phi\underset{2}{\circ}(\omega^{1,0}\wedge P^{0,1}-P^{1,0}\wedge\omega^{0,1}))(\Xh,\Yah)(\Uh,\Wah)&=&-\sqrt{-1}\sum_{i}\bigl(\phi(\Xh,\Yah)(Z_i,\Wah)P(\Uh,\overline{Z_i})\nonumber\\
&&-\phi(\Xh,\Yah)(\Uh,\overline{Z_i})P(Z_i,\Wah)\bigr).
\end{eqnarray}
Using (\ref{e57}), we obtain
$$(\sdelb \overline{C}_{\Yc_{0}})(\Xh)(\Wah)
=-\sqrt{-1}\Bigl(mB_{\Cc_{0}}^{1,0}+{\bigl(C_{\Mc}^{1,1}\underset{2}{\circ} P^{\Cc}\bigr)}^{1,0}\Bigr)(\Xh)(\Wah).\quad\Box$$

\begin{Rem}
\begin{enumerate}
\item Since we made the assumption that $M$ is a Fefferman-Robinson manifold, then the tensor $A_{\n{}{}}$ vanishes and consequently most of Bianchi identities derived in Propositions 5.2 and 5.3 are similar to those obtained by \cite{JL} and \cite{CG} for the curvature of the Tanaka-Webster connection (\cite{Ta},\cite{We}) of a strictly pseudoconvex $CR$-manifold.
\item Note that in \cite{TC2} Bianchi identities are derived for the curvature of the Levi-Civita connection associated to almost Robinson manifolds with twist-induced almost Robinson structure (cf. \cite{FLTC2},\cite{TC2}).
\end{enumerate}
\end{Rem}

\begin{Prop} Let $(M^{2m+2},g,\Nc,\Rc,\Rcs)$ be a Einstein-Fefferman-Robinson manifold such that $\nab{\Sc}{}A_{\n{}{*}}=0$. Then
\begin{equation}\label{e58}
 C^{\Mc}=-\frac{s^{\Cc}}{m}\Bigl(\frac{1}{m(m+1)}(\omega^{1,0}\wedge \omega^{0,1}+\omega\otimes\omega)+\frac{1}{\tr\,{[A^{1,0}_{\n{}{*}},A^{0,1}_{\n{}{*}}]}_1}A^{1,0}_{\n{}{*}}\wedge A^{0,1}_{\n{}{*}}\Bigr).
\end{equation}
\end{Prop}

\noindent{Proof.} First, since $M$ is Einstein-Fefferman-Robinson manifold then $Ric^{\Cc}=-\dfrac{\sqrt{-1}s^{\Cc}}{m}\omega$. It follows from (47) that $\partial_{\Sc} s^{\Cc}=-\sqrt{-1}m{(\sdel A^{1,0}_{\n{}{*}})}^{1,0}$. Since $\nab{\Sc}{}A_{\n{}{*}}=0$, we have $\del s^{\Cc}=\delb s^{\Cc}=\rho^{1,0}=\rho^{0,1}=0$ and by (48) $\nab{\n{}{*}}{}s^{\Cc}=0$. So $s^{\Cc}$ is constant. Now, since $\ds P^{\Cc}=-\frac{\sqrt{-1}}{2m(m+1)}s^{\Cc}\omega$, it follows that
\begin{eqnarray*}
 {\{A^{1,0}_{\n{}{*}},P_{\Cc}^{0,1}\}}_1+{\{A^{0,1}_{\n{}{*}},P_{\Cc}^{1,0}\}}_1=-\frac{\sqrt{-1}s^{\Cc}}{2m(m+1)}({\{A^{1,0}_{\n{}{*}},\omega^{0,1}\}}_1+{\{A^{0,1}_{\n{}{*}},\omega^{1,0}\}}_1)=-\frac{s^{\Cc}}{m(m+1)}(A^{2,0}_{\n{}{*}}-A^{0,2}_{\n{}{*}})
\end{eqnarray*}
and
$${[P_{\Cc}^{1,0},P_{\Cc}^{0,1}]}_1=-\frac{{s^{\Cc}}^2}{4m^2{(m+1)}^2}{[\omega^{1,0},\omega^{0,1}]}_1=-\frac{\sqrt{-1}{s^{\Cc}}^2}{4m^2{(m+1)}^2}\omega.$$
We deduce that
\begin{eqnarray*}
Q_{\Cc}&=&-\nab{\n{}{*}}{}A^{2,0}_{\n{}{*}}-\nab{\n{}{*}}{}A^{0,2}_{\n{}{*}}+\sqrt{-1}\dfrac{s^{\Cc}}{m(m+1)}(A^{2,0}_{\n{}{*}}-A^{0,2}_{\n{}{*}})\\
B_{\Cc_{0}}&=&-{[A^{1,0}_{\n{}{*}},A^{0,1}_{\n{}{*}}]}_1-\dfrac{\sqrt{-1}}{m}(\tr\,{[A^{1,0}_{\n{}{*}},A^{0,1}_{\n{}{*}}]}_1)\omega.
\end{eqnarray*}
Now, for $\phi\in\Omega_{\Sc}^{1,1}(M;\wedge^{1,1}{(\Sc^{\C})}^*)$, we have
\begin{eqnarray}\label{e59}
(\phi(\Zh,\Wah)A^{1,0}_{\n{}{*}})(\Xh)(\Yh)
&=&-\sum_{i}\phi(\Zh,\Wah)(\Xh,\overline{Z_i})A_{\n{}{*}}(Z_i,\Yh)\nonumber\\
&&-\sum_{i}\phi(\Zh,\Wah)(\Yh,\overline{Z_i})A_{\n{}{*}}(Z_i,\Xh).
\end{eqnarray}
Using (\ref{e59}), we obtain :
\begin{eqnarray}\label{e60}
{[C_{\Mc}^{1,1}(\Zh,\Wah)A^{1,0}_{\n{}{*}},A^{0,1}_{\n{}{*}}]}_1(\Xh,\Yah)&=&\bigl(C_{\Mc}^{1,1}\underset{2}{\circ} (A^{1,0}_{\n{}{*}}\wedge A^{0,1}_{\n{}{*}})\bigr)(\Zh,\Wah)(\Xh,\Yah)\\&-&\sqrt{-1}\bigl(C_{\Mc}^{1,1}\underset{2}{\circ} ({[A^{1,0}_{\n{}{*}},A^{0,1}_{\n{}{*}}]}_{1}^{1,0}\wedge\omega^{0,1})\bigr)(\Zh,\Wah)(\Xh,\Yah)\nonumber.
\end{eqnarray}
Now, using (\ref{e50}), we also obtain that
\begin{eqnarray}\label{e61}
\bigl(((A^{1,0}_{\n{}{*}}\wedge\omega^{0,1})\underset{2}{\circ} C_{\Mc}^{1,1})\underset{2}{\circ}(\omega^{0,1}\wedge A^{0,1}_{\n{}{*}})\bigr)(\Zh,\Wah)(\Xh,\Yah)=-\bigl(C_{\Mc}^{1,1}\underset{2}{\circ} (A^{1,0}_{\n{}{*}}\wedge A^{0,1}_{\n{}{*}})\bigr)(\Zh,\Wah)(\Xh,\Yah)\nonumber\\
-\sqrt{-1}\bigl(C_{\Mc}^{1,1}\underset{2}{\circ} ({[A^{1,0}_{\n{}{*}},A^{0,1}_{\n{}{*}}]}_{1}^{1,0}\wedge\omega^{0,1})\bigr)(\Zh,\Wah)(\Xh,\Yah).
\end{eqnarray}
Adding (\ref{e60}) and (\ref{e61}), we obtain that
\begin{eqnarray*}
{[C_{\Mc}^{1,1}(\Zh,\Wah)A^{1,0}_{\n{}{*}},A^{0,1}_{\n{}{*}}]}_1(\Xh,\Yah)&=&
-\bigl(((A^{1,0}_{\n{}{*}}\wedge\omega^{0,1})\underset{2}{\circ} C_{\Mc}^{1,1})\underset{2}{\circ}(\omega^{0,1}\wedge A^{0,1}_{\n{}{*}})\bigr)(\Zh,\Wah)(\Xh,\Yah)\\
&&-2
\sqrt{-1}\bigl(C_{\Mc}^{1,1}\underset{2}{\circ} ({[A^{1,0}_{\n{}{*}},A^{0,1}_{\n{}{*}}]}_{1}^{1,0}\wedge\omega^{0,1})\bigr)(\Zh,\Wah)(\Xh,\Yah).
\end{eqnarray*}
Since $\nabla_{\Sc}A_{\n{}{*}}=0$ and $M$ is Einstein-Fefferman then $C_{\Yc_{0}}=0$. Also it follows from (51) that $$((A^{1,0}_{\n{}{*}}\wedge\omega^{0,1})\underset{2}{\circ} C_{\Mc}^{1,1})\underset{2}{\circ}(\omega^{0,1}\wedge A^{0,1}_{\n{}{*}})=\sqrt{-1}(\omega^{1,0}\wedge Q_{\Cc}^{1,0})\underset{2}{\circ}(\omega^{0,1}\wedge A^{0,1}_{\n{}{*}}).$$
By Lemma 5.3, we have $(\omega^{1,0}\wedge Q_{\Cc}^{1,0})\underset{2}{\circ}(\omega^{0,1}\wedge A^{0,1}_{\n{}{*}})=Q_{\Cc}^{1,0}\wedge A^{0,1}_{\n{}{*}}+\sqrt{-1}{[Q_{\Cc}^{1,0},A^{0,1}_{\n{}{*}}]}_{1}^{1,0}\wedge\omega^{0,1}$ and so
\begin{eqnarray}\label{e62}
{[C_{\Mc}^{1,1}(\Zh,\Wah)A^{1,0}_{\n{}{*}},A^{0,1}_{\n{}{*}}]}_1(\Xh,\Yah)
&=&-\sqrt{-1}\Bigl(Q_{\Cc}^{1,0}\wedge A^{0,1}_{\n{}{*}}
+\sqrt{-1}{[Q_{\Cc}^{1,0},A^{0,1}_{\n{}{*}}]}_{1}^{1,0}\wedge\omega^{0,1}\\
&+&2\bigl(C_{\Mc}^{1,1}\underset{2}{\circ} ({[A^{1,0}_{\n{}{*}},A^{0,1}_{\n{}{*}}]}_{1}^{1,0}\wedge\omega^{0,1})\bigr)\Bigr)(\Zh,\Wah)(\Xh,\Yah)\nonumber.
\end{eqnarray}
Since $\ds P^{\Cc}=-\frac{\sqrt{-1}}{2m(m+1)}s^{\Cc}\omega$, then (\ref{e16}) yields
$$C_{\Mc}^{1,1}=R_{\Cc}^{1,1}+\frac{s^{\Cc}}{m(m+1)}\omega^{1,0}\wedge\omega^{0,1}+\frac{s^{\Cc}}{2m(m+1)}{(\omega\odot\omega)}^{1,1}.$$
Also we have
$${[(\omega^{1,0}\wedge\omega^{0,1})(\Zh,\Wah)A^{1,0}_{\n{}{*}},A^{0,1}_{\n{}{*}}]}_1=(A^{1,0}_{\n{}{*}}\wedge A^{0,1}_{\n{}{*}})(\Zh,\Wah)-\sqrt{-1}({[A^{1,0}_{\n{}{*}},A^{0,1}_{\n{}{*}}]}_{1}^{1,0}\wedge\omega^{0,1})(\Zh,\Wah)$$
and
$${[{(\omega\odot\omega)}^{1,1}(\Zh,\Wah)A^{1,0}_{\n{}{*}},A^{0,1}_{\n{}{*}}]}_1
=-4\sqrt{-1}\omega(\Zh,\Wah){[A^{1,0}_{\n{}{*}},A^{0,1}_{\n{}{*}}]}_{1}.$$
We deduce that
\begin{eqnarray}\label{e63}
{[C_{\Mc}^{1,1}(\Zh,\Wah)A^{1,0}_{\n{}{*}},A^{0,1}_{\n{}{*}}]}_1&=&
{[R_{\Cc}^{1,1}(\Zh,\Wah)A^{1,0}_{\n{}{*}},A^{0,1}_{\n{}{*}}]}_1
+\frac{s^{\Cc}}{m(m+1)}\Bigl((A^{1,0}_{\n{}{*}}\wedge A^{0,1}_{\n{}{*}})(\Zh,\Wah)\\
&-&\sqrt{-1}({[A^{1,0}_{\n{}{*}},A^{0,1}_{\n{}{*}}]}_{1}^{1,0}\wedge\omega^{0,1})(\Zh,\Wah)\
-2\sqrt{-1}\omega(\Zh,\Wah){[A^{1,0}_{\n{}{*}},A^{0,1}_{\n{}{*}}]}_{1}\Bigr)\nonumber.
\end{eqnarray}
Now, since $\nab{\Sc}{}A^{1,0}_{\n{}{*}}=\nab{\n{}{}}{}A^{1,0}_{\n{}{*}}=0$, we have
$$R^{\Cc}(\Zh,\Wah)A^{1,0}_{\n{}{*}}=(\nab{\Zh,\Wah}{2}-\nab{\Wah,\Zh}{2})A^{1,0}_{\n{}{*}}+\omega(\Zh,\Wah)\nab{\n{}{*}}{}A^{1,0}_{\n{}{*}}=\omega(\Zh,\Wah)\nab{\n{}{*}}{}A^{1,0}_{\n{}{*}}.$$
Also, using (\ref{e59}) together with $M$ Einstein-Fefferman, we obtain that
$$(\wedge_{\Sc}R_{\Cc}^{1,1})A^{1,0}_{\n{}{*}}=-\sqrt{-1}{\{A^{1,0}_{\n{}{*}},Ric_{\Cc}^{0,1}\}}^{1,0}_{1}=-\frac{s^{\Cc}}{m}{\{A^{1,0}_{\n{}{*}},\omega^{0,1}\}}^{1,0}_{1}=\sqrt{-1}\frac{2s^{\Cc}}{m}A^{1,0}_{\n{}{*}}
=m\nab{\n{}{*}}{}A^{1,0}_{\n{}{*}}.$$
Hence
$$\nab{\n{}{*}}{}A^{1,0}_{\n{}{*}}=\sqrt{-1}\frac{2s^{\Cc}}{m^2}A^{1,0}_{\n{}{*}}.$$
So we deduce that
$$R^{\Cc}(\Zh,\Wah)A^{1,0}_{\n{}{*}}=\sqrt{-1}\frac{2s^{\Cc}}{m^2}\omega(\Zh,\Wah)A^{1,0}_{\n{}{*}}$$
and
\begin{equation}\label{e64}
Q_{\Cc}^{1,0}=-\sqrt{-1}s^{\Cc}\dfrac{m+2}{m^2(m+1)}A^{1,0}_{\n{}{*}}.
\end{equation}
Hence (\ref{e63}) becomes
\begin{eqnarray}\label{e65}
{[C_{\Mc}^{1,1}(\Zh,\Wah)A^{1,0}_{\n{}{*}},A^{0,1}_{\n{}{*}}]}_1&=&\frac{s^{\Cc}}{m(m+1)}\Bigl((A^{1,0}_{\n{}{*}}\wedge A^{0,1}_{\n{}{*}})(\Zh,\Wah)\\
&-&\sqrt{-1}({[A^{1,0}_{\n{}{*}},A^{0,1}_{\n{}{*}}]}_{1}^{1,0}\wedge\omega^{0,1})(\Zh,\Wah)\
+\frac{2\sqrt{-1}}{m}\omega(\Zh,\Wah){[A^{1,0}_{\n{}{*}},A^{0,1}_{\n{}{*}}]}_{1}\Bigr)\nonumber.
\end{eqnarray}
By $C_{\Yc_{0}}=0$ together with the assumption Einstein-Fefferman, we obtain by (54) that $B_{\Cc_{0}}=0$. Also we have
\begin{equation}\label{e66}
\ds {[A^{1,0}_{\n{}{*}},A^{0,1}_{\n{}{*}}]}_1=-\dfrac{\sqrt{-1}}{m}(\tr\,{[A^{1,0}_{\n{}{*}},A^{0,1}_{\n{}{*}}]}_1)\omega.
\end{equation}
Substituting (\ref{e64}) and (\ref{e66}) in (\ref{e62}) and (\ref{e65}) together with the assumption $\del C_{\Yc_{0}}^{1,0}=0$, yields
\begin{eqnarray*}
{[C_{\Mc}^{1,1}(\Zh,\Wah)A^{1,0}_{\n{}{*}},A^{0,1}_{\n{}{*}}]}_1(\Xh,\Yah)&=&
-\frac{s^{\Cc}(m+2)}{m^2(m+1)}(A^{1,0}_{\n{}{*}}\wedge A^{0,1}_{\n{}{*}})(\Zh,\Wah)(\Xh,\Yah)\\
&-&\frac{s^{\Cc}(m+2)}{m^3(m+1)}(\tr\,{[A^{1,0}_{\n{}{*}},A^{0,1}_{\n{}{*}}]}_1)(\omega^{1,0}\wedge\omega^{0,1})(\Zh,\Wah)(\Xh,\Yah)\\
&-&\frac{2}{m}(\tr\,{[A^{1,0}_{\n{}{*}},A^{0,1}_{\n{}{*}}]}_1)C_{\Mc}^{1,1}(\Zh,\Wah)(\Xh,\Yah)\\
&=&\frac{s^{\Cc}}{m(m+1)}(A^{1,0}_{\n{}{*}}\wedge A^{0,1}_{\n{}{*}})(\Zh,\Wah)(\Xh,\Yah)\\
&-&\frac{s^{\Cc}}{m^2(m+1)}(\tr\,{[A^{1,0}_{\n{}{*}},A^{0,1}_{\n{}{*}}]}_1)(\omega^{1,0}\wedge\omega^{0,1})(\Zh,\Wah)(\Xh,\Yah)\\
&+&\frac{2s^{\Cc}}{m^3(m+1)}(\tr\,{[A^{1,0}_{\n{}{*}},A^{0,1}_{\n{}{*}}]}_1)(\omega\otimes\omega)(\Zh,\Wah)(\Xh,\Yah).
\end{eqnarray*}
We deduce that
\begin{eqnarray*}
 C_{\Mc}^{1,1}(\Zh,\Wah)(\Xh,\Yah)=&-&\frac{s^{\Cc}}{m}\Bigl(\frac{1}{m(m+1)}(\omega^{1,0}\wedge \omega^{0,1}+\omega\otimes\omega)(\Zh,\Wah)(\Xh,\Yah)\\
 &+&\frac{1}{\tr\,{[A^{1,0}_{\n{}{*}},A^{0,1}_{\n{}{*}}]}_1}A^{1,0}_{\n{}{*}}\wedge A^{0,1}_{\n{}{*}}\Bigr)(\Zh,\Wah)(\Xh,\Yah).\quad \Box
\end{eqnarray*}

\begin{Rem}
It follows from $$C_{\Mc}^{1,1}=R_{\Cc}^{1,1}+\frac{s^{\Cc}}{m(m+1)}\omega^{1,0}\wedge\omega^{0,1}+\frac{s^{\Cc}}{2m(m+1)}{(\omega\odot\omega)}^{1,1}$$ that $R^{\Cc}$ is given by:
\begin{equation}\label{e67}
R^{\Cc}=-\frac{s^{\Cc}}{m}\Bigl(\frac{1}{m}(\omega^{1,0}\wedge \omega^{0,1}+\omega\otimes\omega)+\frac{1}{\tr\,{[A^{1,0}_{\n{}{*}},A^{0,1}_{\n{}{*}}]}_1}A^{1,0}_{\n{}{*}}\wedge A^{0,1}_{\n{}{*}}\Bigr).
\end{equation}
\end{Rem}

Examples of Fefferman-Robinson manifolds with $C^{\Mc}$ given by (\ref{e58}) can be found in the class of locally symmetric Fefferman-Robinson manifolds which are special cases of locally almost Robinson Witt symmetric spaces (cf. \cite{Pe2}) and that are analogous of contact locally subsymmetric spaces (\cite{BFG}). Now, we give the definitions of these spaces.

\begin{Def} Let $(M^{2m+2},g,\Nc,\Rc,\Rcs)$ be an almost Robinson manifold. A diffeomorphism $\psi:M\to M$ is called an (complex) optical diffeomorphism if $d\psi(\Rc)\subset\Rc,d\psi(\Rcs)\subset\Rcs, d\psi(\Sc)\subset\Sc\;\mathrm{and}\;d\psi\circ J=J\circ d\psi$ (which is also equivalent to $d\psi(\Rc)\subset\Rc,d\psi(\Rcs)\subset\Rcs, d\psi(\So{1}{0})\subset\So{1}{0}$ and $d\psi(\So{0}{1})\subset\So{0}{1}$).
An optical isometry $\psi:M\to M$ is an isometric optical diffeomorphism ($\psi^{*}g=g$).
\end{Def}

\begin{Def} A (locally) symmetric Fefferman-Robinson manifold is a strongly geodetic almost Fefferman-Robinson manifold ($M^{2m+2},g,\Nc,\Rc,\Rcs)$ such that, for every point $x_0\in M$, there exists an optical (local) isometry $\psi$, called optical symmetry at $x_0$, satisfying $\psi(x_0)=x_0$ and $d\psi(x_0)_{/\Sc_{x_0}}=-id_{/\Sc_{x_0}}$.
\end{Def}

\begin{Prop} Let $(M^{2m+2},g,\Nc,\Rc,\Rcs,\psi)$ be a locally symmetric Fefferman-Robinson manifold endowed with its Chern-Robinson connection and let $(\n{}{},\n{}{*})$ be an adapted $g$-optical pairing. Then we have
\begin{enumerate}
\item $\Nc$ is integrable.
\item ${(\Lc_{\n{}{}}\sig{}{})}_{\Sc}={(\Lc_{\n{}{*}}\sig{}{})}_{\Sc}=0$, ${(\nab{\Sc}{}d\sig{}{})}_{\Sc}=\nabla_{\Sc}A_{\n{}{*}}=0$.
\item $\nab{\Sc}{}R=0$, ${(i_{\n{}{}}R)}_{\Sc}={(i_{\n{}{*}}R)}_{\Sc}=0$.
\end{enumerate}
\end{Prop}

\noindent{Proof.} Let $(M^{2m+2},g,\Nc,\Rc,\Rcs,\psi)$ be a locally symmetric Fefferman-Robinson manifold, $\psi:U_x\to U_x$ be an optical symmetry at $x\in M$ (with $U_x$ openset) and $(\n{}{},\n{}{*})$ be an adapted $g$-optical pairing. Then we have
$$d\psi(\n{}{})=\mu\n{}{},\;d\psi(\n{}{*})=\lambda\n{}{*},\;d\psi_{\Sc}(\So{1}{0})\subset\So{1}{0},\;d\psi_{\Sc}(\So{0}{1})\subset\So{0}{1}\;\mathrm{and}\;d\psi(x)_{/\Sc_{x}}=-id_{/\Sc_{x}},$$
together with $\lambda,\mu\in\Cc^{\infty}(U_x)$ satisfying $\lambda\mu=1$.
Since $d\psi(\n{}{})=\mu\n{}{}$ then $\psi^{*}\sig{}{*}=\mu\sig{}{*}$ and we have $$d(\psi^{*}\sig{}{*})=\psi^{*}(d\sig{}{*})=d
\mu\wedge\sig{}{*}+\mu d\sig{}{*}.$$
Since $d\psi(x)_{/\Sc_{x}}=-id_{/\Sc_{x}}$, then, for $X,Y\in \Sc_{x}$,
$$\psi^{*}(d\sig{}{*})(X,Y)=(d\sig{}{*})(d\psi(X),d\psi(Y))=d\sig{}{*}(X,Y)=\mu d\sig{}{*}(X,Y).$$
Hence $(\mu-1) d\sig{}{*}(X,Y)=(\mu-1)\omega(X,Y)=0$. Since $\omega$ is nondegenerate on $\Sc_{x}$, then $\mu=1$ and $\lambda=1$ at $x$.
Since $\psi$ is an optical local isometry, then using the formula
\begin{eqnarray*}
g((\nab{X}{}d\psi)(Y),d\psi(Z))&=&\frac{1}{2}\Bigl(g(T(d\psi(X),d\psi(Y)),d\psi(Z))-g(T(d\psi(X),d\psi(Z)),d\psi(Y))\\
                             &&-g(T(d\psi(Y),d\psi(Z)),d\psi(X))\\
                             &&-g(T(X,Y),Z)+g(T(X,Z),Y)+g(T(Y,Z),X)\Bigr),\quad X,Y,Z\in\Gamma(T^{\C}M),
\end{eqnarray*}
together with (\ref{e11}), we obtain that
$g((\nab{X}{}d\psi)(\Yh),d\psi(\Zah))=g((\nab{X}{}d\psi)(\n{}{}),d\psi(\n{}{*}))=0$. Also $\nab{X}{}d\psi=0$ and $\psi$ is affine for the Chern-Robinson connection $\nab{}{}$. Since $\psi$ is affine, then, for $X,Y,Z\in\Gamma(T^{\C}M)$,
$$T(d\psi(X),d\psi(Y))=d\psi(T(X,Y))\,\mathrm{and}\,
(\nab{d\psi(X)}{}T)(d\psi(Y),d\psi(Z))=d\psi((\nab{X}{}T)(Y,Z)).$$
Since $T(\Xh,\Yh)=-[\Xh,\Yh]^{0,1}+d\sig{}{}(\Xh,\Yh)\n{}{}$ and $d\psi(\Xh)=-\Xh,d\psi(\Yh)=-\Yh$ at $x$, we obtain
$$T(d\psi(\Xh),d\psi(\Yh))=T(\Xh,\Yh)=d\psi(T(\Xh,\Yh))=[\Xh,\Yh]^{0,1}+d\sig{}{}(\Xh,\Yh)\n{}{}.$$
Also $N(\Xh,\Yh)=[\Xh,\Yh]^{0,1}=0$ and $[\Gamma(\So{1}{0}),\Gamma(\So{1}{0})]\subset\Gamma(\So{1}{0})$. As $M$ is assumed strongly geodetic
$[\Gamma(\Rc^{\C}),\Gamma(\Nc)]\subset\Gamma(\Nc)$, we deduce that $\Nc$ is integrable.\\
From $d\psi(\n{}{})=\n{}{}$ and $d\psi(\n{}{*})=\n{}{*}$ at $x$ together with $T(\n{}{*},X)=(\Lc_{\n{}{*}}\sig{}{})(X)\n{}{}+\dfrac{1}{2}(\Lc_{\n{}{}}\sig{}{})(X)\n{}{*}+\tu{\Sc}(\n{}{*},X)$ for $X\in\Sc_{x}$, we deduce that :
$$T(d\psi(\n{}{*}),d\psi(X))=-T(\n{}{*},X)=d\psi(T(\n{}{*},X))=(\Lc_{\n{}{*}}\sig{}{})(X)\n{}{}+\frac{1}{2}(\Lc_{\n{}{}}\sig{}{})(X)\n{}{*}-\tu{\Sc}(\n{}{*},X).$$
Hence $(\Lc_{\n{}{}}\sig{}{})(X)=(\Lc_{\n{}{*}}\sig{}{})(X)=0$.\\
Now, for $X,Y,Z\in\Sc_{x}$,
$$(\nab{d\psi(X)}{}T)(d\psi(\n{}{*}),d\psi(Z))=(\nab{X}{}T)(\n{}{*},Z)=(\nab{X}{}\tu{\Sc})(\n{}{*},Z)=d\psi((\nab{X}{}\tu{\Sc})(\n{}{*},Z))=-(\nab{X}{}\tu{\Sc})(\n{}{*},Z).$$
We deduce that $(\nab{X}{}\tu{\Sc})(\n{}{*},Z)=0$ and $(\nab{X}{}A_{\n{}{*}})(Z,W)=0$ for $Z,W\in\Sc_{x}$.
\begin{eqnarray*}
 (\nab{d\psi(X)}{}T)(d\psi(Y),d\psi(Z))&=&-(\nab{X}{}T)(Y,Z)=-(\nab{X}{}d\sig{}{})(Y,Z)\n{}{}-d\sig{}{}(Y,Z)\nab{X}{}\n{}{}-\omega(Y,Z)\nab{X}{}\n{}{*}\\
 &=&d\psi((\nab{X}{}T)(Y,Z))=(\nab{X}{}d\sig{}{})(Y,Z)\n{}{}+d\sig{}{}(Y,Z)\nab{X}{}\n{}{}+\omega(Y,Z)\nab{X}{}\n{}{*}.
\end{eqnarray*}
Hence, $(\nab{X}{}d\sig{}{})(Y,Z)\n{}{}+d\sig{}{}(Y,Z)\nab{X}{}\n{}{}+\omega(Y,Z)\nab{X}{}\n{}{*}=0$. From $\omega$ nondegenerate on $\Sc_{x}$, we deduce that $g(\nab{X}{}\n{}{*},\n{}{})=0$.
Also,
$$\nab{X}{}\n{}{*}=g(\nab{X}{}\n{}{*},\n{}{})\n{}{*}=0,\;\nab{X}{}\n{}{}=-g(\nab{X}{}\n{}{*},\n{}{})\n{}{}=0\;\mathrm{and}\;(\nab{X}{}d\sig{}{})(Y,Z)=0.$$
It follows from (\ref{e12}) that $R(\n{}{(*)},Y)=0$, for any $Y\in\Sc$.\\
Since $\psi$ is affine, then, for $X,Y,Z,W\in\Gamma(T^{\C}M)$,
$$R(d\psi(X),d\psi(Y))d\psi(Z)=d\psi(R(X,Y)Z)\,\mathrm{and}\,
(\nab{d\psi(X)}{}R)(d\psi(Y),d\psi(Z))d\psi(W)=d\psi((\nab{X}{}R)(Y,Z)W).$$
Now, for $X,Y,Z,W\in\Sc_{x}$, we obtain :
\begin{eqnarray*}
(\nab{d\psi(X)}{}R)(d\psi(Y),d\psi(Z))d\psi(W)&=&(\nab{X}{}R)(Y,Z)W=d\psi((\nab{X}{}R)(Y,Z)W)=-(\nab{X}{}R)(Y,Z)W\\
(\nab{d\psi(X)}{}R)(d\psi(Y),d\psi(Z))d\psi(\n{}{(*)})&=&-(\nab{X}{}R)(Y,Z)\n{}{(*)}=d\psi((\nab{X}{}R)(Y,Z)\n{}{(*)})=(\nab{X}{}R)(Y,Z)\n{}{(*)}\\
(\nab{d\psi(X)}{}R)(d\psi(\n{}{}),d\psi(\n{}{*}))d\psi(W)&=&(\nab{X}{}R)(\n{}{},\n{}{*})W=d\psi((\nab{X}{}R)(\n{}{},\n{}{*})W)=-(\nab{X}{}R)(\n{}{},\n{}{*})W\\
(\nab{d\psi(X)}{}R)(d\psi(\n{}{}),d\psi(\n{}{*}))d\psi(\n{}{(*)})&=&-(\nab{X}{}R)(\n{}{},\n{}{*})\n{}{(*)}=d\psi((\nab{X}{}R)(\n{}{},\n{}{*})\n{}{(*)})=(\nab{X}{}R)(\n{}{},\n{}{*})\n{}{(*)}.
\end{eqnarray*}
Also $(\nab{X}{}R)(Y,Z)W=(\nab{X}{}R)(Y,Z)\n{}{(*)}=(\nab{X}{}R)(\n{}{},\n{}{*})W=(\nab{X}{}R)(\n{}{},\n{}{*})\n{}{(*)}=0$. We deduce that $\nab{X}{}R=0$.
$\Box$
\begin{Cor} Let $(M^{2m+2},g,\Nc,\Rc,\Rcs)$ be a locally symmetric Fefferman-Robinson manifold Einstein-Fefferman endowed with its Chern-Robinson connection, then the formulas (\ref{e58}) and (\ref{e67}) are valid.
\end{Cor}

\section{Appendix}

The proof of proposition 4.1 needs the two following lemmas

\begin{Le} Let $g_{\ld}=e^{2\ld}g$ and let $\nab{}{\ld}$ (resp $\nab{}{}$)
be a $g_{\ld}$ (resp $g$) connection, then for any $X,Y,Z\in\Gamma(TM)$, we have :
\begin{equation}\label{e68}
(\nab{X}{\ld}Y^{\flat_{\ld}})(Z)=e^{2\ld}\Bigl((\nab{X}{}Y^{\flat})(Z)+\Kc^{\ld}(X,Y,Z)-\Kc(X,Y,Z)+(d\ld(X)Y^{\flat}+d\ld(Y)X^{\flat}-g(X,Y)d\ld)(Z)\Bigr),
\end{equation}
with 
$$\Kc^{\ld}(X,Y,Z)=\frac{1}{2}\Bigl(g(T^{\ld}(X,Y),Z)-g(T^{\ld}(X,Z),Y)-g(T^{\ld}(Y,Z),X)\Bigr)\;\mathrm{and}\;\Kc=\Kc^{(\ld=0)}.$$
\end{Le}

\noindent{Proof.} 
By (\ref{e9}), we have the formula
\begin{equation}\label{e69}
(\nab{X}{\ld}Y^{\flat_{\ld}})(Z)=\frac{1}{2}\Bigl(dY^{\flat_{\ld}}(X,Z)+(\Lc_{Y}g_{\ld})(X,Z)\Bigr)+\Kc^{\nab{}{\ld}}(X,Y,Z),
\end{equation}
with $\ds\Kc^{\nab{}{\ld}}(X,Y,Z)=\frac{1}{2}\Bigl(g_{\ld}(T^{\ld}(X,Y),Z)-g_{\ld}(T^{\ld}(X,Z),Y)-g_{\ld}(T^{\ld}(Y,Z),X)\Bigr)$.\\ 
Since $g_{\ld}=e^{2\ld}g$, we have $Y^{\flat_{\ld}}=e^{2\ld}Y^{\flat}$. Hence
$$dY^{\flat_{\ld}}=e^{2\ld}(2d\lambda\wedge Y^{\flat}+dY^{\flat})\;\mathrm{and}\;\Lc_{Y}g_{\ld}=e^{2\ld}(2d\lambda(Y)g+\Lc_{Y}g).$$ We deduce that (\ref{e69}) becomes :
\begin{eqnarray}\label{e70}
(\nab{X}{\ld}Y^{\flat_{\ld}})(Z)&=&e^{2\ld}\Bigl[\frac{1}{2}\Bigl(dY^{\flat}(X,Z)+(\Lc_{Y}g)(X,Z)\Bigr)+d\ld(X)g(Y,Z)+d\ld(Y)g(X,Z)-d\ld(Z)g(X,Y)\nonumber\\
&+&\frac{1}{2}\Bigl(g(T^{\ld}(X,Y),Z)-g(T^{\ld}(X,Z),Y)-g(T^{\ld}(Y,Z),X)\Bigr)\Bigr].
\end{eqnarray}
The formula directly follows by substituting (\ref{e9}) in (\ref{e70}). $\Box$\\

\begin{Le}
Let $(M^{2m+2},g,\Nc,\Rc,\Rcs)$ be a geodetic almost Robinson manifold and $\ld\in\Cc_{\nu{}{}}^{\infty}(M)$. Then, for any $\Xhl,\Yhl\in\Gamma(\Sl{1}{0})$ and $\Xahl,\Yahl,\Zahl\in\Gamma(\Sl{0}{1})$, we have :
\begin{eqnarray}\label{e71}
\Kc(\Xhl,\Yhl,\Zahl)&=&-\frac{\sqrt{-1}}{2}d\omega(\Xhl,\Yhl,\Zahl)+\dlsl(\Xhl)\Bigl((\Lc_{\n{}{}}\omega)(\Yhl,\Zahl)-d\sig{}{*}(\Yhl,\Jl\Zahl)\Bigr)\nonumber\\
&-&\dlsl(\Yhl)d\sig{}{*}(\Xhl,\Jl\Zahl)\nonumber\\
&-&\dlsl(\Jl\Zahl)\Bigl(d\sig{}{*}(\Xhl,\Yhl)+(\Lc_{\n{}{}}g)(\Xhl,\Yhl)+\sqrt{-1}(\Lc_{\n{}{}}\omega)(\Xhl,\Yhl)\Bigr)
\end{eqnarray}
\begin{eqnarray}\label{e72}
\Kc(\Xahl,\Yhl,\Zahl)&=&-\frac{\sqrt{-1}}{2}d\omega(\Xahl,\Zahl,\Yhl)+\dlsl(\Xahl)\Bigl((\Lc_{\n{}{}}\omega)(\Yhl,\Zahl)+d\sig{}{*}(\Yhl,\Jl\Zahl)\Bigr)\nonumber\\
&+&\dlsl(\Jl\Yhl)\Bigl(d\sig{}{*}(\Xahl,\Zahl)+(\Lc_{\n{}{}}g)(\Xahl,\Zahl)-\sqrt{-1}(\Lc_{\n{}{}}\omega)(\Xahl,\Zahl)\Bigr)\nonumber\\
&+&\dlsl(\Zahl)d\sig{}{*}(\Yhl,\Jl\Xahl)
\end{eqnarray}
\begin{equation}\label{e73}
\Kc(\n{}{},\Yhl,\Zahl)=-\frac{1}{2}d\sig{}{*}(\Yhl,\Zahl)
\end{equation}
\begin{eqnarray}\label{e74}
\Kc(\n{}{*},\Yhl,\Zahl)&=&-\frac{1}{2}d\sig{}{}(\Yhl,\Zahl)+\dlsl(\Jl\Yhl)\Bigl((\Lc_{\n{}{*}}\sig{}{*})(\Zahl)-g([\n{}{},\n{}{*}]^{{(1,0)}_{\ld}}_{\scl},\Zahl)\Bigr)\nonumber\\
&-&\dlsl(\Jl\Zahl)\Bigl((\Lc_{\n{}{*}}\sig{}{*})(\Yhl)-g([\n{}{},\n{}{*}]^{{(0,1)}_{\ld}}_{\scl},\Yhl)\Bigr).
\end{eqnarray}
\end{Le}

\noindent{Proof.} For any $\Xhl,\Yhl\in\Gamma(\Sl{1}{0})$ and $\Xahl,\Yahl,\Zahl\in\Gamma(\Sl{0}{1})$, we have :
\begin{equation}\label{e75}
\Kc(\Xhl,\Yhl,\Zahl)=\frac{1}{2}\Bigl(g(T(\Xhl,\Yhl),\Zahl)-g(T(\Xhl,\Zahl),\Yhl)-g(T(\Yhl,\Zahl),\Xhl)\Bigr).
\end{equation}
Since $\Xhl=\Xh-2\sqrt{-1}\dls(\Xh)\nu{}{},\;\Yahl=\Yah+2\sqrt{-1}\dls(\Yah)\nu{}{}$, $\Xh\in\Gamma(\So{1}{0}),\;\Yah\in\Gamma(\So{0}{1})$ and $M$ is assumed to be geodetic (i.e. $(\Lc_{\n{}{}}\sig{}{*})(\Xh)=0$), then we obtain using (\ref{e9}) that :
\begin{eqnarray}\label{e76}
T(\Xhl,\Yhl)&=&\Tc^{1,0}(\Xh,\Yh)+2\sqrt{-1}\dls(\Yh)\Tc^{1,0}(\n{}{},\Xh)-2\sqrt{-1}\dls(\Xh)\Tc^{1,0}(\n{}{},\Yh)\nonumber\\
&-&[\Xh,\Yh]^{0,1}+2\sqrt{-1}\dls(\Xh)[\n{}{},\Yh]^{0,1}-2\sqrt{-1}\dls(\Yh)[\n{}{},\Xh]^{0,1}\nonumber\\
&+&\Bigl(d\sig{}{}(\Xh,\Yh)+\sqrt{-1}\dls(\Xh)(\Lc_{\Yh}g)(\n{}{},\n{}{*})-\sqrt{-1}\dls(\Yh)(\Lc_{\Xh}g)(\n{}{},\n{}{*})\Bigr)\n{}{}\nonumber\\
&+&d\sig{}{*}(\Xh,\Yh)\n{}{*}\nonumber\\
T(\Xhl,\Yahl)&=&2\sqrt{-1}\dls(\Xh)[\n{}{},\Yah]^{1,0}-2\sqrt{-1}\dls(\Yah)\Tc^{1,0}(\n{}{},\Xh)\nonumber\\
&+&2\sqrt{-1}\dls(\Yah)[\n{}{},\Xh]^{0,1}-2\sqrt{-1}\dls(\Xh)\Tc^{0,1}(\n{}{},\Yah)\nonumber\\
&+&\Bigl(d\sig{}{}(\Xh,\Yah)+\sqrt{-1}\dls(\Xh)(\Lc_{\Yh}g)(\n{}{},\n{}{*})+\sqrt{-1}\dls(\Yh)(\Lc_{\Xh}g)(\n{}{},\n{}{*})\Bigr)\n{}{}\nonumber\\
&+&d\sig{}{*}(\Xh,\Yah)\n{}{*}.
\end{eqnarray}
Substituting (\ref{e76}) in (\ref{e75}), we obtain
\begin{eqnarray*}
\Kc(\Xhl,\Yhl,\Zahl)&=&\sqrt{-1}\Bigl(-\frac{1}{2}d\omega(\Xh,\Yh,\Zah)\\
&+&\dls(\Yh)(\Lc_{\n{}{}}g)(\Xh,\Zah)+\dls(\Zah)(\Lc_{\n{}{}}g)(\Xh,\Yh)\\
&+&\dls(\Xh)d\sig{}{*}(\Yh,\Zah)+\dls(\Yh)d\sig{}{*}(\Xh,\Zah)+\dls(\Zah)d\sig{}{*}(\Xh,\Yh)\Bigr).
\end{eqnarray*}
Now, we have also 
\begin{eqnarray*}
d\omega(\Xhl,\Yhl,\Zahl)&=&d\omega(\Xh,\Yh,\Zah)+2\sqrt{-1}\Bigl(\dls(\Zah)(\Lc_{\n{}{}}\omega)(\Xh,\Yh)\\
&+&\dls(\Yh)(\Lc_{\n{}{}}\omega)(\Xhl,\Zahl)-\dls(\Xh)(\Lc_{\n{}{}}\omega)(\Yh,\Zah)\Bigr). 
\end{eqnarray*}
Since $d\ld(\nu{}{})=0$, then $\dlsl(\Xhl)=\dls(\Xh)$ and $\dlsl(\Jl\Xhl)=\dls(J\Xh)=\sqrt{-1}\dls(\Xh)$. 
Hence, we obtain, together with the assumption $(\Lc_{\n{}{}}\sig{}{*})(\Xh)=0$, that :
\begin{eqnarray*}
\Kc(\Xhl,\Yhl,\Zahl)&=&-\frac{\sqrt{-1}}{2}d\omega(\Xhl,\Yhl,\Zahl)+\dlsl(\Xhl)\Bigl((\Lc_{\n{}{}}\omega)(\Yhl,\Zahl)-d\sig{}{*}(\Yhl,\Jl\Zahl)\Bigr)\\
&-&\dlsl(\Yhl)d\sig{}{*}(\Xhl,\Jl\Zahl)\\
&-&\dlsl(\Jl\Zahl)\Bigl(d\sig{}{*}(\Xhl,\Yhl)+(\Lc_{\n{}{}}g)(\Xhl,\Yhl)+\sqrt{-1}(\Lc_{\n{}{}}\omega)(\Xhl,\Yhl)\Bigr).
\end{eqnarray*}
In the same way, 
\begin{eqnarray*}
\Kc(\Xahl,\Yhl,\Zahl)&=&\frac{1}{2}\Bigl(g(T(\Xahl,\Yhl),\Zahl)-g(T(\Xahl,\Zahl),\Yhl)-g(T(\Yhl,\Zahl),\Xahl)\Bigr)\\
&=&-\frac{\sqrt{-1}}{2}d\omega(\Xahl,\Zahl,\Yhl)+\dlsl(\Xahl)\Bigl((\Lc_{\n{}{}}\omega)(\Yhl,\Zahl)+d\sig{}{*}(\Yhl,\Jl\Zahl)\Bigr)\\
&+&\dlsl(\Jl\Yhl)\Bigl(d\sig{}{*}(\Xahl,\Zahl)+(\Lc_{\n{}{}}g)(\Xahl,\Zahl)-\sqrt{-1}(\Lc_{\n{}{}}\omega)(\Xahl,\Zahl)\Bigr)\\
&+&\dlsl(\Zahl)d\sig{}{*}(\Yhl,\Jl\Xahl).
\end{eqnarray*}
Now, we have 
\begin{eqnarray}\label{e77}
T(\n{}{},\Yhl)&=&T(\n{}{},\Yh)=\Tc^{1,0}(\n{}{},\Yh)-[\n{}{},\Yh]^{0,1}-\frac{1}{2}(\Lc_{\Yh}g)(\n{}{},\n{}{*})\n{}{}\nonumber\\
T(\n{}{*},\Yhl)&=&T(\n{}{*},\Yh)-2\sqrt{-1}\dls(\Yh)T(\n{}{*},\n{}{})\nonumber\\
&=&\Tc^{1,0}(\n{}{*},\Yh)-[\n{}{*},\Yh]^{0,1}-2\sqrt{-1}\dls(\Yh)[\n{}{},\n{}{*}]_{\Sc}\nonumber\\
&+&(\Lc_{\n{}{*}}\sig{}{})(\Yh)\n{}{}-\frac{1}{2}(\Lc_{\Yh}g)(\n{}{},\n{}{*})\n{}{*}.
\end{eqnarray}
It directly follows from (\ref{e76}) and (\ref{e77}) that :
\begin{eqnarray*}
\Kc(\n{}{},\Yhl,\Zahl)&=&\frac{1}{2}\Bigl(g(T(\n{}{},\Yhl),\Zahl)-g(T(\n{}{},\Zahl),\Yhl)
-g(T(\Yhl,\Zahl),\n{}{})\Bigr)\\
&=&-\frac{1}{2}d\sig{}{*}(\Yh,\Zah)=-\frac{1}{2}d\sig{}{*}(\Yhl,\Zahl).
\end{eqnarray*}
Now by (\ref{e76}) and (\ref{e77}), we have
\begin{eqnarray*}
\Kc(\n{}{*},\Yhl,\Zahl)&=&\frac{1}{2}\Bigl(g(T(\n{}{*},\Yhl),\Zahl)-g(T(\n{}{*},\Zahl),\Yhl)
-g(T(\Yhl,\Zahl),\n{}{*})\Bigr)\\
&=&-\frac{1}{2}d\sig{}{}(\Yh,\Zah)-\sqrt{-1}\dls(\Yh)\Bigl((\Lc_{\Zah}g)(\n{}{},\n{}{*})+g([\n{}{},\n{}{*}]^{1,0}_{\Sc},\Zah)\Bigr)\\
&-&\sqrt{-1}\dls(\Zah)\Bigl((\Lc_{\Yh}g)(\n{}{},\n{}{*})+g([\n{}{},\n{}{*}]^{0,1}_{\Sc},\Yh)\Bigr).
\end{eqnarray*}
Since $$(\Lc_{\Yh}g)(\n{}{},\n{}{*})=-(\Lc_{\n{}{}}\sig{}{})(\Yh)-(\Lc_{\n{}{*}}\sig{}{*})(\Yh),$$
then
\begin{eqnarray*}
\Kc(\n{}{*},\Yhl,\Zahl)&=&-\frac{1}{2}d\sig{}{}(\Yh,\Zah)+\sqrt{-1}\dls(\Yh)(\Lc_{\n{}{}}\sig{}{})(\Zah)
+\sqrt{-1}\dls(\Zah)(\Lc_{\n{}{}}\sig{}{})(\Yh)\\
&+&\sqrt{-1}\dls(\Yh)\Bigl((\Lc_{\n{}{*}}\sig{}{*})(\Zah)-g([\n{}{},\n{}{*}]^{1,0}_{\Sc},\Zah)\Bigr)\\
&+&\sqrt{-1}\dls(\Zah)\Bigl((\Lc_{\n{}{*}}\sig{}{*})(\Yh)-g([\n{}{},\n{}{*}]^{0,1}_{\Sc},\Yh)\Bigr))\\
&=&-\frac{1}{2}d\sig{}{}(\Yhl,\Zahl)+\dlsl(\Jl\Yhl)\Bigl((\Lc_{\n{}{*}}\sig{}{*})(\Zahl)-g([\n{}{},\n{}{*}]^{{(1,0)}_{\ld}}_{\scl},\Zahl)\Bigr)\\
&-&\dlsl(\Jl\Zahl)\Bigl((\Lc_{\n{}{*}}\sig{}{*})(\Yhl)-g([\n{}{},\n{}{*}]^{{(0,1)}_{\ld}}_{\scl},\Yhl)\Bigr).
\end{eqnarray*}
The result follows. $\Box$\\

\noindent{Proof of the proposition 4.1} We have by (\ref{e68}) :
\begin{eqnarray}\label{e78}
(\nab{\Xhl}{\ld}{\Yhl}^{\flat_{\ld}})(\Zahl)&=&e^{2\ld}\Bigl((\nab{\Xhl}{}{\Yhl}^{\flat})(\Zahl)+\Kc^{\ld}(\Xhl,\Yhl,\Zahl)-\Kc(\Xhl,\Yhl,\Zahl)\nonumber\\
&+&\bigl(\dlsl(\Xhl){\Yhl}^{\flat}+\dlsl(\Yhl){\Xhl}^{\flat}-g(\Xhl,\Yhl)\dlsl\bigr)(\Zahl)\Bigr)\nonumber\\
&=&e^{2\ld}(\nab{\Xhl}{}{\Yhl}^{\flat})(\Zahl)+e^{2\ld}\Kc^{\ld}(\Xhl,\Yhl,\Zahl)-e^{2\ld}\Kc(\Xhl,\Yhl,\Zahl)\nonumber\\
&+&\dlsl(\Xhl)\gl(\Yhl,\Zahl)+\dlsl(\Yhl)\gl(\Xhl,\Zahl).
\end{eqnarray}

Using (\ref{e8}), we have :
\begin{eqnarray*}
e^{2\ld}\Kc^{\ld}(\Xhl,\Yhl,\Zahl)&=&\frac{1}{2}\Bigl(\gl(T^{\ld}(\Xhl,\Yhl),\Zahl)-\gl(T^{\ld}(\Xhl,\Zahl),\Yhl)-\gl(T^{\ld}(\Yhl,\Zahl),\Xhl)\Bigr)\\ 
&=&\frac{1}{2}\gl(T^{\ld}(\Xhl,\Yhl),\Zahl)=-\frac{\sqrt{-1}}{2}d\omega_{\ld}(\Xhl,\Yhl,\Zahl)\\
&=&-\frac{\sqrt{-1}}{2}e^{2\ld}d\omega(\Xhl,\Yhl,\Zahl)+\dlsl(\Xhl)\Bigl(\gl(\Yhl,\Zahl)-d\sig{}{*}(\Yhl,\Jl\Zahl)\Bigr)\\
&-&\dlsl(\Yhl)\Bigl(\gl(\Xhl,\Zahl)-d\sig{}{*}(\Xhl,\Jl\Zahl)\Bigr)\\
&-&\dlsl(\Jl\Zahl)d\sig{}{*}(\Xhl,\Yhl).
\end{eqnarray*}

Now substituting the last formula and (\ref{e71}) in (\ref{e78}) yields the formula for
$(\nab{\Xhl}{\ld}{\Yhl}^{\flat_{\ld}})(\Zahl)$.
Using (\ref{e72}) and (\ref{e73}), we similary prove formulas for $(\nab{\Xahl}{\ld}{\Yhl}^{\flat_{\ld}})(\Zahl)$ and
$(\nab{\n{\ld}{}}{\ld}{\Yhl}^{\flat_{\ld}})(\Zahl)$.

Now, we have  
\begin{eqnarray}\label{e79}
(\nab{\n{\ld}{*}}{\ld}{\Yhl}^{\flat_{\ld}})(\Zahl)&=&e^{2\ld}\Bigl((\nab{\n{\ld}{*}}{}{\Yhl}^{\flat})(\Zahl)+\Kc^{\ld}(\n{\ld}{*},\Yhl,\Zahl)-\Kc(\n{\ld}{*},\Yhl,\Zahl)\nonumber\\
&+&\bigl(d\ld(\n{\ld}{*}){\Yhl}^{\flat}+\dlsl(\Yhl){\n{\ld}{*}}^{\flat}-g(\n{\ld}{*},\Yhl)\dlsl\bigr)(\Zahl)\Bigr)\nonumber\\
&=&e^{2\ld}(\nab{\n{\ld}{*}}{}{\Yhl}^{\flat})(\Zahl)+e^{2\ld}\Kc^{\ld}(\n{\ld}{*},\Yhl,\Zahl)-e^{2\ld}\Kc(\n{\ld}{*},\Yhl,\Zahl)\nonumber\\
&+&d\ld(\n{\ld}{*})\gl(\Yhl,\Zahl).
\end{eqnarray}

Using (\ref{e8}), we have :
\begin{eqnarray}\label{e80}
e^{2\ld}\Kc^{\ld}(\n{\ld}{*},\Yhl,\Zahl)&=&\frac{1}{2}\Bigl(\gl(T^{\ld}(\n{\ld}{*},\Yhl),\Zahl)-\gl(T^{\ld}(\n{\ld}{*},\Zahl),\Yhl)-\gl(T^{\ld}(\Yhl,\Zahl),\n{\ld}{*})\Bigr)\nonumber\\ 
&=&-\frac{1}{2}\sig{\ld}{}(T^{\ld}(\Yhl,\Zahl))=-\frac{1}{2}d\sig{\ld}{}(\Yhl,\Zahl).
\end{eqnarray}

Now, we have 
\begin{eqnarray*}
\n{\ld}{*}&=&e^{-2\ld}\bigl(\n{}{*}-2J{(d\ld)}^{\sharp}_{\Sc}-2\dlsn\n{}{}\bigr)=e^{-2\ld}\n{}{*}-2\Jl\dlsd{\sharp_\ld}+2e^{-2\ld}\dlsn\n{}{}\\
&=&e^{-2\ld}\bigl(\n{}{*}-2\sqrt{-1}\dlsd{\sharp\,{(1,0)}_{\ld}}+2\sqrt{-1}\dlsd{\sharp\,{(0,1)}_{\ld}}+2\dlsn\n{}{}\bigr).
\end{eqnarray*}
Hence
\begin{eqnarray*}
e^{2\ld}\Kc(\n{\ld}{*},\Yhl,\Zahl)&=&\Kc(\n{}{*},\Yhl,\Zahl)-2\sqrt{-1}\Kc(\dlsd{\sharp\,{(1,0)}_{\ld}},\Yhl,\Zahl)+2\sqrt{-1}\Kc(\dlsd{\sharp\,{(0,1)}_{\ld}},\Yhl,\Zahl)\\
&+&2\dlsn\Kc(\n{}{},\Yhl,\Zahl). 
\end{eqnarray*}
Using (\ref{e71}),(\ref{e72}),(\ref{e73}) and (\ref{e74}), we obtain :
\begin{eqnarray}\label{e81}
e^{2\ld}\Kc(\n{\ld}{*},\Yhl,\Zahl)&=&-\frac{1}{2}d\sig{}{}(\Yhl,\Zahl)+\dlsn d\sig{}{*}(\Yhl,\Zahl)-d\omega(\dlsd{\sharp},\Yhl,\Zahl)\nonumber\\
&+&\dlsl(\Jl\Yhl)\Bigl((\Lc_{\n{}{*}}\sig{}{*})(\Zahl)-2\sqrt{-1}d\sig{}{*}(\dlsd{\sharp\,{(1,0)}_{\ld}},\Zahl)+2\sqrt{-1}d\sig{}{*}(\dlsd{\sharp\,{(0,1)}_{\ld}},\Zahl))\nonumber\\
&-&g([\n{}{},\n{}{*}]^{{(1,0)}_{\ld}}_{\scl},\Zahl)+2\sqrt{-1}(\Lc_{\n{}{}}g)(\dlsd{\sharp\,{(0,1)}_{\ld}},\Zahl)+2(\Lc_{\n{}{}}\omega)(\dlsd{\sharp\,{(0,1)}_{\ld}},\Zahl)\Bigr)\nonumber\\
&-&\dlsl(\Jl\Zahl)\Bigl((\Lc_{\n{}{*}}\sig{}{*})(\Yhl)-2\sqrt{-1}d\sig{}{*}(\dlsd{\sharp\,{(1,0)}_{\ld}},\Yhl)+2\sqrt{-1}d\sig{}{*}(\dlsd{\sharp\,{(0,1)}_{\ld}},\Yhl))\nonumber\\
&-&g([\n{}{},\n{}{*}]^{{(0,1)}_{\ld}}_{\scl},\Yhl)-2\sqrt{-1}(\Lc_{\n{}{}}g)(\dlsd{\sharp\,{(1,0)}_{\ld}},\Yhl)+2(\Lc_{\n{}{}}\omega)(\dlsd{\sharp\,{(1,0)}_{\ld}},\Yhl)\Bigr).
\end{eqnarray}
Since 
\begin{eqnarray*}
g([\n{}{},\dlsd{\sharp\,{(1,0)}_{\ld}}]^{{(0,1)}_{\ld}}_{\scl},\Yhl)&=&-\frac{1}{2}\Bigl((\Lc_{\n{}{}}g)(\dlsd{\sharp\,{(1,0)}_{\ld}},\Yhl)+\sqrt{-1}(\Lc_{\n{}{}}\omega)(\dlsd{\sharp\,{(1,0)}_{\ld}},\Yhl)\Bigr)\\
g([\n{}{},\dlsd{\sharp\,{(0,1)}_{\ld}}]^{{(1,0)}_{\ld}}_{\scl},\Zahl)&=&-\frac{1}{2}\Bigl((\Lc_{\n{}{}}g)(\dlsd{\sharp\,{(0,1)}_{\ld}},\Zahl)-\sqrt{-1}(\Lc_{\n{}{}}\omega)(\dlsd{\sharp\,{(0,1)}_{\ld}},\Zahl)\Bigr)
\end{eqnarray*}
and 
$$(\Lc_{\dlsd{\sharp}}\omega)(\Yhl,\Zahl)=d\omega(\dlsd{\sharp},\Yhl,\Zahl)-(dJ^{*}\dls)(\Yhl,\Zahl),$$
then (\ref{e81}) becomes
\begin{eqnarray}\label{e82}
e^{2\ld}\Kc(\n{\ld}{*},\Yhl,\Zahl)&=&-\frac{1}{2}\Bigl(d\sig{}{}(\Yhl,\Zahl)+2(dJ^{*}\dls)(\Yhl,\Zahl)-2\dlsn d\sig{}{*}(\Yhl,\Zahl)\Bigr)\nonumber\\
&+&\dlsl(\Jl\Yhl)\Bigl(d\sig{}{*}(e^{2\ld}\n{\ld}{*},\Zahl)-g([\n{}{},\n{}{*}+4\sqrt{-1}\dlsd{\sharp\,{(0,1)}_{\ld}}]^{{(1,0)}_{\ld}}_{\scl},\Zahl)\Bigr)\nonumber\\
&-&\dlsl(\Jl\Zahl)\Bigl(d\sig{}{*}(e^{2\ld}\n{\ld}{*},\Yhl)-g([\n{}{},\n{}{*}-4\sqrt{-1}\dlsd{\sharp\,{(1,0)}_{\ld}}]^{{(0,1)}_{\ld}}_{\scl},\Yhl)\Bigr)\nonumber\\
&-&(\Lc_{\dlsd{\sharp}}\omega)(\Yhl,\Zahl).
\end{eqnarray}
Now, we have 
$$d\sig{\ld}{}=d\sig{}{}+2dJ^{*}\dls-2d(\dlsn)\wedge\sig{}{*}-2\dlsn d\sig{}{*},$$
$$\n{}{*}+4\sqrt{-1}\dlsd{\sharp\,{(0,1)}_{\ld}}=e^{2\ld}\n{\ld}{*}+2\sqrt{-1}\dlsd{\sharp}-2\dlsn\n{}{}$$
$$\n{}{*}-4\sqrt{-1}\dlsd{\sharp\,{(1,0)}_{\ld}}=e^{2\ld}\n{\ld}{*}-2\sqrt{-1}\dlsd{\sharp}-2\dlsn\n{}{}$$
and
$$d\sig{}{*}(e^{2\ld}\n{\ld}{*},\Zahl)=d\sig{\ld}{*}(\n{\ld}{*},\Zahl)+2\dlsl(\Zahl),\quad d\sig{}{*}(e^{2\ld}\n{\ld}{*},\Yhl)=d\sig{\ld}{*}(\n{\ld}{*},\Yhl)+2\dlsl(\Yhl).$$
Also (\ref{e82}) becomes
\begin{eqnarray}\label{e83}
e^{2\ld}\Kc(\n{\ld}{*},\Yhl,\Zahl)&=&-\frac{1}{2}d\sig{\ld}{}(\Yhl,\Zahl)
+2\Bigl(\dlsl(\Jl\Yhl)\dlsl(\Zahl)-\dlsl(\Jl\Zahl)\dlsl(\Yhl)\Bigr)\nonumber\\
&+&\dlsl(\Jl\Yhl)\Bigl((\Lc_{\n{\ld}{*}}\sig{\ld}{*})(\Zahl)-\gl([\n{}{},\n{\ld}{*}]^{{(1,0)}_{\ld}}_{\scl},\Zahl)\Bigr)\nonumber\\
&-&\dlsl(\Jl\Zahl)\Bigl((\Lc_{\n{\ld}{*}}\sig{\ld}{*})(\Yhl)-\gl([\n{}{},\n{\ld}{*}]^{{(0,1)}_{\ld}}_{\scl},\Yhl)\Bigr)\nonumber\\
&-&2\sqrt{-1}\Bigl(\dlsl(\Jl\Yhl)g([\n{}{},\dlsd{\sharp}]^{{(1,0)}_{\ld}}_{\scl},\Zahl)
+\dlsl(\Jl\Zahl)g([\n{}{},\dlsd{\sharp}]^{{(0,1)}_{\ld}}_{\scl},\Yhl)\Bigr)\nonumber\\
&-&(\Lc_{\dlsd{\sharp}}\omega)(\Yhl,\Zahl)\nonumber\\
&=&-\frac{1}{2}d\sig{\ld}{}(\Yhl,\Zahl)
+2\Bigl(\dlsl(\Jl\Yhl)\dlsl(\Zahl)-\dlsl(\Jl\Zahl)\dlsl(\Yhl)\Bigr)\nonumber\\
&+&\dlsl(\Jl\Yhl)\Bigl((\Lc_{\n{\ld}{*}}\sig{\ld}{*})(\Zahl)-\gl([\n{}{},\n{\ld}{*}]^{{(1,0)}_{\ld}}_{\scl},\Zahl)\Bigr)\nonumber\\
&-&\dlsl(\Jl\Zahl)\Bigl((\Lc_{\n{\ld}{*}}\sig{\ld}{*})(\Yhl)-\gl([\n{}{},\n{\ld}{*}]^{{(0,1)}_{\ld}}_{\scl},\Yhl)\Bigr)\nonumber\\
&-&2\Bigl(\dlsl(\Jl\Yhl)(\Lc_{\dlsd{\sharp}}\omega)(\n{}{},\Zahl)
-\dlsl(\Jl\Zahl)(\Lc_{\dlsd{\sharp}}\omega)(\n{}{},\Yhl)\Bigr)\nonumber\\
&-&(\Lc_{\dlsd{\sharp}}\omega)(\Yhl,\Zahl).
\end{eqnarray}
Now the result follows by substituting (\ref{e80}) and (\ref{e83}) in (\ref{e79}). $\Box$\\

\noindent{Proof of the corollary 4.1.} First since $(M^{2m+2},g,\Nc,\Rc,\Rcs)$ is a strongly geodetic almost Fefferman-Robinson manifold then $(M^{2m+2},\gl,\Nc_{\ld},\Rc_{\ld},\Rcs_{\ld})$ is also a strongly geodetic almost Fefferman-Robinson manifold. Hence we have by proposition 2.1 and lemma 4.1 that $\Lc_{\n{}{}}\ol=0$, $[\n{}{},\n{\ld}{*}]_{\scl}=0$, ${(\Lc_{\n{}{}}\gl)}_{\scl}=0$ and 
$$\gl(\Xhl,\Zahl)=d\sig{\ld}{*}(\Xhl,\Jl\Zahl),\;d\sig{\ld}{*}(\Xhl,\Yhl)=d\sig{\ld}{*}(\n{\ld}{},\Zahl)=d\sig{\ld}{*}(\n{\ld}{*},\Zahl)=0.$$
Now, since $\dlsd{\sharp\,{(1,0)}_{\ld}}={(d\ld)}^{\sharp\,1,0}_{\Sc}-\sqrt{-1}\dlsn\n{}{}$ then
\begin{equation}\label{e84}
\dlsl(\Zahl)=g(\dlsd{\sharp\,{(1,0)}_{\ld}},\Zahl)=-e^{-2\ld}\gl(\dlsd{\sharp\,{(1,0)}_{\ld}},\Zahl).
\end{equation}
Now, using (\ref{e84}), then the two first equations in (\ref{e20}) becomes :
\begin{eqnarray*}
\gl(\nab{\Xhl}{\ld}\Yhl,\Zahl)&=&\gl({\bigl(\nab{\Xhl}{}\Yhl\bigr)}^{{(1,0)}_{\ld}}_{\scl}+2\bigl(\dlsl(\Xhl)\Yhl+\dlsl(\Yhl)\Xhl\bigr),\Zahl)\\
\gl(\nab{\Xahl}{\ld}\Yhl,\Zahl)&=&\gl({\bigl(\nab{\Xahl}{}\Yhl\bigr)}^{{(1,0)}_{\ld}}_{\scl}-2g(\Xahl,\Yhl)\dlsd{\sharp\,{(1,0)}_{\ld}},\Zahl).
\end{eqnarray*}
Now, we have 
\begin{eqnarray*}
(\Lc_{\dlsd{\sharp}}\omega)(\Yhl,\Zahl)&=&d\omega(\dlsd{\sharp},\Yhl,\Zahl)-(dJ^{*}\dls)(\Yhl,\Zahl)\\
(\Lc_{\dlsd{\sharp}}\omega)(\n{}{},\Zahl)&=&d\omega(\dlsd{\sharp},\n{}{},\Zahl)-(dJ^{*}\dls)(\n{}{},\Zahl)\\
(\Lc_{\dlsd{\sharp}}\omega)(\n{}{},\Yhl)&=&d\omega(\dlsd{\sharp},\n{}{},\Yhl)-(dJ^{*}\dls)(\n{}{},\Yhl).
\end{eqnarray*}
Since $\omega=d\sig{}{*}$, then $d\omega=0$, also to obtain the fourth equation, we have 
to calculate $(dJ^{*}\dls)(\Yhl,\Zahl)$, $(dJ^{*}\dls)(\n{}{},\Zahl)$ and $(dJ^{*}\dls)(\n{}{},\Yhl)$.
\\ First, recall the following formula :
$$d\alpha(X,Y)=(\nab{X}{}\alpha)(Y)-(\nab{Y}{}\alpha)(X)+\alpha(T(X,Y)).$$
Now, applying the previous formula to $\dls=d\ld-\sig{}{*}\otimes d\ld(\n{}{*})$ ($d\ld(\n{}{})=0$), we obtain for $X,Y\in\Gamma(\mathrm{Ker}\,\sig{}{*})$ that :
\begin{equation}\label{e85}
(\nab{X}{}\dls)(Y)-(\nab{Y}{}\dls)(X)=-\dls(T(X,Y))-d\sig{}{*}(X,Y)d\ld(\n{}{*}).
\end{equation}
Also, using (\ref{e85}), we have
\begin{eqnarray*}
(dJ^{*}\dls)(\Yhl,\Zahl)&=&(\nab{\Yhl}{}J^{*}\dls)(\Zahl)-(\nab{\Zahl}{}J^{*}\dls)(\Yhl)+\dls(JT(\Yhl,\Zahl))\\
&=&-\sqrt{-1}((\nab{\Yhl}{}\dls)(\Zahl)+(\nab{\Zahl}{}\dls)(\Yhl))+\dls(JT(\Yhl,\Zahl))\\
&=&-2\sqrt{-1}(\nab{\Yhl}{}\dls)(\Zahl)-\sqrt{-1}\dls\bigl(T(\Yhl,\Zahl)+\sqrt{-1}JT(\Yhl,\Zahl)\bigr)\\
&-&\sqrt{-1}d\sig{}{*}(\Yhl,\Zahl)d\ld(\n{}{*})\\
&=&-2\sqrt{-1}\bigl((\nab{\Yhl}{}\dls)(\Zahl)+\dls(T^{0,1}(\Yhl,\Zahl)\bigr)+g(\Yhl,\Zahl)d\ld(\n{}{*})
\end{eqnarray*}
\begin{eqnarray*}
(dJ^{*}\dls)(\n{}{},\Zahl)&=&(\nab{\n{}{}}{}J^{*}\dls)(\Zahl)-(\nab{\Zahl}{}J^{*}\dls)(\n{}{})+\dls(JT(\n{}{},\Zahl))\\
&=&-\sqrt{-1}(\nab{\n{}{}}{}\dls)(\Zahl)+\dls(JT(\n{}{},\Zahl))\\
&=&\sqrt{-1}\dls\bigl(T(\n{}{},\Zahl)-\sqrt{-1}JT(\n{}{},\Zahl)\bigr)+\sqrt{-1}(\Lc_{\n{}{}}\sig{}{*})(\Zahl)d\ld(\n{}{*})\\
&=&2\sqrt{-1}\dls(T^{1,0}(\n{}{},\Zahl))
\end{eqnarray*}
and
$$
(dJ^{*}\dls)(\n{}{},\Yhl)=-2\sqrt{-1}\dls(T^{0,1}(\n{}{},\Yhl)).
$$
Now using (\ref{e76}) and (\ref{e77}), we have
$$\dls(T^{0,1}(\Yhl,\Zahl))=2\sqrt{-1}\dls(\Zah)\dls([\n{}{},\Yh]^{0,1})=0$$
and
$$\dls(T^{1,0}(\n{}{},\Zahl))=-\dls([\n{}{},\Zah]^{1,0})=0,\quad \dls(T^{0,1}(\n{}{},\Yhl))=-\dls([\n{}{},\Yh]^{0,1})=0.$$

We deduce that : 
\begin{eqnarray}\label{e86}
(dJ^{*}\dls)(\Yhl,\Zahl)&=&g(-2\nab{\Yhl}{}\Jl\dlsd{\sharp\,{(1,0)}_{\ld}}+d\ld(\n{}{*})\Yhl,\Zahl)\nonumber\\
(dJ^{*}\dls)(\n{}{},\Zahl)&=&(dJ^{*}\dls)(\n{}{},\Yhl)=0.
\end{eqnarray}
Substituting (\ref{e86}) in the last equation of (\ref{e20}), we obtain :
$$
\gl(\nab{\n{\ld}{*}}{\ld}\Yhl,\Zahl)=\gl({\bigl(\nab{\n{\ld}{*}}{}\Yhl\bigr)}^{{(1,0)}_{\ld}}_{\scl}
+2e^{-2\ld}\Bigl(\nab{\Yhl}{}\Jl\dlsd{\sharp\,{(1,0)}_{\ld}}-2\dlsl(\Yhl)\Jl\dlsd{\sharp\,{(1,0)}_{\ld}}\Bigr),\Zahl).
$$
$\Box$

\end{document}